\documentclass[preprint,12pt]{elsarticle}

\usepackage{amssymb}
\usepackage{amsmath}
\usepackage{amsthm}
\usepackage{mathtools}
\usepackage{dsfont}
\usepackage[capitalise]{cleveref}
\usepackage[printonlyused,withpage]{acronym}
\usepackage{caption}
\usepackage{subfig}
\usepackage{algorithm}
\usepackage{algpseudocode}
\usepackage{subfig}
\usepackage{multirow}

\usepackage{color}
\usepackage{booktabs}

\usepackage{framed} 
\usepackage{multicol} 
\usepackage{nomencl} 
\makenomenclature
 \biboptions{sort&compress}
\newtheorem{theorem}{Theorem}

\newcommand{\gf}{\Omega}
\newcommand{\gb}{\varGamma}
\newcommand{\bu}{\boldsymbol{u}}
\newcommand{\hbu}{\boldsymbol{\hat{u}}}
\newcommand{\obu}{\boldsymbol{\overline{u}}}
\newcommand{\ubu}{\boldsymbol{\underline{u}}}
\newcommand{\reu}{\boldsymbol{\mathrm{\tilde{u}}}}
\newcommand{\ru}{\boldsymbol{\mathrm{u}}}
\newcommand{\rv}{\boldsymbol{\mathrm{v}}}

\newcommand{\bv}{\boldsymbol{v}}
\newcommand{\s}{\mathfrak{s}}
\renewcommand{\u}{{u}}
\renewcommand{\v}{{v}}
\renewcommand{\t}{t}
\newcommand{\bx}{\boldsymbol{x}}
\newcommand{\x}{{x}}
\newcommand{\y}{{y}}
\newcommand{\p}{p}
\newcommand{\op}{\overline{p}}
\newcommand{\hp}{\hat{p}}
\newcommand{\rp}{\mathrm{p}}
\newcommand{\rep}{\mathrm{\tilde{p}}}
\newcommand{\q}{\mathrm{q}}
\renewcommand{\Re}{\mathrm{Re}}
\renewcommand{\k}{k}
\newcommand{\rF}{\mathrm{F}}
\newcommand{\reF}{\mathrm{\overline{F}}}
\newcommand{\rG}{\mathrm{G}}
\newcommand{\gbn}{\varGamma_{\mathrm{N}}}
\newcommand{\gbd}{\varGamma_{\mathrm{D}}}
\newcommand{\gd}{{G_D}}
\newcommand{\gn}{{G_N}}
\newcommand{\f}{{f}}
\newcommand{\strain}{\epsilon}
\newcommand{\rf}{\mathrm{f}}
\newcommand{\rg}{\mathrm{g}}
\newcommand{\n}{\boldsymbol {n}}

\renewcommand{\P}{\mathcal{Q}}
\newcommand{\Z}{\mathcal{Z}}
\newcommand{\Pk}[1]{\mathds{P}_{#1}}
\newcommand{\Th}{\mathcal{T}_h}
\newcommand{\K}{\mathrm{K}}
\newcommand{\Nv}{\mathrm{N_v}}
\newcommand{\Ne}{\mathrm{N_e}}
\newcommand{\bev}{\boldsymbol{\psi}}
\newcommand{\bep}{{\phi}}
\newcommand{\clv}{\mathit{v}}
\newcommand{\clq}{\mathit{q}}
\newcommand{\lin}{\boldsymbol{\ell}}
\newcommand{\lins}{\mathrm{L}}

\newcommand{\mA}{\mathsf{A}}
\newcommand{\mK}{\mathsf{K}}
\newcommand{\mB}{\mathsf{B}}
\newcommand{\mU}{\mathsf{U}}
\newcommand{\mP}{\mathsf{P}}

\newcommand{\mF}{\mathsf{F}}
\newcommand{\X}{\mathsf{X}}

\newcommand{\Res}{\mathcal{R}es}
\newcommand{{\m}}{\mathcal{N}}
\newcommand{{\mn}}{n}
\newcommand{\h}{\mathrm{h}}
\newcommand{\sT}{\mathds{T}}

\newcommand{\up}[1]{$^{\mathrm{#1}}$}
\newcommand{\tol}{\mathsf{tol}}

\renewcommand{\t}{\mathrm{t}}
\newcommand{\I}{\mathbf{I}}

\newcommand{\V}{\mathcal{V}}
\newcommand{\J}{\mathcal{J}}
\newcommand{\C}[1]{\mathrm{C}_{#1}}
\newcommand{\T}[1]{\mathbf{T}_{#1}}

\renewcommand{\d}{\mathrm{d}}

\renewcommand{\b}{\mathrm{b}}

\newcommand{\w}{\mathrm{w}}

\newcommand{\Vr}{\mathring{\V}}
\newcommand{\Vrh}{\Vr_h}

\newcommand\dg\mathfrak

\definecolor{myred}{RGB}{255,0,0}

\journal{Mathematics and Computers in Simulation}

\begin{document}

\begin{frontmatter}



\title{Numerical Integration of Navier-Stokes Equations by Time Series Expansion and Stabilized FEM}

\author[label1]{Ahmad Deeb}
\ead{ahmad.deeb@ku.ac.ae}
\author[label1,label2]{Denys Dutykh}
\ead{denys.dutykh@ku.ac.ae}

\affiliation[label1]{organization={Khalifa University of Science and Technology},
            addressline={PO Box 127788},
            city={Abu Dhabi},
            country={United Arab Emirates}}

\affiliation[label2]{organization={Causal Dynamics},
            addressline={Pty LTD},
            city={Perth},
            country={Australia}}



\begin{abstract}
This manuscript introduces an advanced numerical approach for the integration of incompressible \ac{NS} equations using a \ac{TSE} method within a \ac{FEM} framework. The technique is enhanced by a novel stabilization strategy, incorporating a
\ac{DSR} technique, which significantly augments the computational efficiency of the algorithm. The stabilization mechanism is meticulously designed to improve the stability and validity of computed series terms, enabling the application of the \ac{FS} algorithm for series resummation. This approach is pivotal in addressing the challenges associated with the accurate and stable numerical solution of \ac{NS} equations, which are critical in \ac{CFD} applications. The manuscript elaborates on the variational formulation of Stokes problem and present convergence analysis of the method using the \ac{LBB} condition. It is followed by the \ac{NS} equations and the implementation details of the stabilization technique, underscored by numerical tests on laminar flow past a cylinder, showcasing the method's efficacy and potential for broad applicability in fluid dynamics simulations. The results of the stabilization indicate a substantial enhancement in computational stability and accuracy, offering a promising avenue for future research in the field.
\end{abstract}


\begin{highlights}
\item Convergence analysis of Stokes problem solution in the TSE-FEM framework.
\item Factorial Series algorithm for the incompressible Navier-Stokes equation in FEM.
\item Stabilization technique to improve space modes computation of the TSE.
\item The formulation is linked  to the Navier Stokes alpha model.
\item The methodology is validated on the flow past a cylinder test case.

\end{highlights}


\end{frontmatter}




\section{Introduction}

\label{sec1}
\acf{CFD} solvers of \acf{NS} equations exist abundantly via multitude of software that are based on different frameworks of space and time discretization of continuous mechanics models \cite{CFD_sagot,CFD-book-09,book-CFSI-12,Archives-1,SPH-archives-15,NM-TF-15,ASPH-5,phase-NS-hydro-23,IBM-CFD-24}. For the space discretization, the \ac{FEM} is one of the framework that has multiple convergent numerical schemes tailored to the \ac{NS} equations \cite{NumScheme-INS-18}. It is also known to deal successfully with complex geometries and having the ability to increase the precision of the space approximation via increasing the polynomial degrees of \acp{FE} \cite{taylor-hood-73,FEM-HO-NS-80, Archives-4,FSI-06,LoggEtal2012,FEM-CFD-14,FEM-FD-14}.
After assembling the \ac{NS} equations for a given spatial framework, the semi-discrete problem is represented by an \ac{IVP} as follows:
\begin{equation}
\label{dyn_pb}
 \frac{\d \y}{\d \t} = \f(\t,\y), \quad \text{with} \quad \y(\t_0) = \y_0,
\end{equation}
Thus, the numerical integration of the above system is to be constructed. This has been widely developed for the incompressible \ac{NS} equations \cite{NS-IVP-64,TIS-NS-CM-94,TIS-NS-CM-16,TIS-NS-01,Assessment-TIS-IMEX-23,TIS-NS-24}. Different schemes have been developed such as the one using the fractional time step \cite{TIS-NS-01,TIS-NS-11,TIS-NS-FTS-22}, time Galerkin-collocation method \cite{Nitsche-NS-21}, the one based on \ac{ETD} \cite{TIS-NS-20}, on \ac{BDF} \cite{TIS-NS-LMS-16,TIS-NS-BDF-20,TIS-NS-BDF-24,TIS-NS-IMEX-BDF-24}, on Linear-Multi Step methods \cite{Layton-22,Layton-23} and \ac{RK} methods \cite{TIS-RK-NS-19} among others. In additions, stabilization techniques of \ac{NS} equations \cite{stab-FEM-NS-92,stab-FEM-NS-01,stab-FEM-NS-15,stab-FEA-NS-20,stab-FEM-NS-23} were developed also in the framework of \ac{FEM}, such as the pressure gradient projection method \cite{stab-FEM-NS-PG-00,Press-stab-stokes-12}, the variational multi-scale method for high Reynolds numbers \cite{stab-VMS-FEM-NS-10} and for adaptive anisotropic meshing \cite{stab-VMS-FEM-NS-13}, the local Projection method \cite{stab-FEM-localProj-16}, the Subgrid \cite{stab-FEM-subgrid-19} and the \ac{SUPG} methods \cite{TIS-FEM-SUPG-NS-97,stab-SUPG-FEM-03,satb-FEM-SUPG-NS-21}.

A new time integration strategy using the \acf{TSE} and \acf{DSR} has been developed by Razafindralandy \cite{dina-thesis} in the framework of \ac{FDM} for the \ac{NS} equations and was transposed to the framework of \ac{FEM} by Deeb \cite{deeb-thesis}. Given the differential equation of a dynamical problem \cref{dyn_pb},
the solution could be represented with the \ac{TSE} in the vicinity of the initial time $\t_0$ as follows:
\begin{equation}
\label{tse}
 \y(\t) = \sum\limits_{\k=0}^{\infty} \y_\k (\t-\t_0)^\k.
\end{equation}
Here, $\y_\k$ are the terms of the series obtained by a recurrence formula:
\begin{equation}
 \y_{\k+1} = \frac{1}{(\k+1)!} \left.\frac{\d^\k \f}{\d \t^\k} \right\vert_{\t=\t_0},
\end{equation}
that is generated after inserting formula \eqref{tse} in \cref{dyn_pb} and equating terms of $(\t - \t_0)^\k$. When the series is convergent, the partial sum is used to approximate the solution for a finite convergence radius $\tau$. However, and when the series is divergent, resummation techniques as the \ac{BL} resummation \cite{borel-1901} are applied to produce function whose asymptotic expansion is the original series. Numerically, this technique was transformed into numerical method such the \ac{BPL} algorithm which was firstly introduced in \cite{dina-2012} and then emphasized in several works \cite{ahmad_bpl_2014,ahmad_comp_bpl_sfg_2015,ahmad_icnpaa_2016}.
The procedure has proven its effectiveness for solving stiff and non-stiff dynamical systems \cite{ahmad_robust_integrators_2019,tayeh-21,DEEB_2022_bpl}.

However, algorithms based on \ac{TSE} present instabilities when using high-order \ac{FEM}. In a recent published paper, Deeb and Dutykh \cite{deeb:stab-serie} have enhanced the computational efficiency of these algorithms, in the framework of higher-order \ac{FEM}, by adding an artificial diffusion term on the left-hand side of the recurrence formulas. Considered as a stabilization technique, it improved the computational rank validity of spatial modes of the series for the one dimensional diffusion and Burgers equation.
It was shown that the value of the diffusion coefficient depends on the rank of the term and is to be found by minimizing the condition number of the new linear system.

In this paper, we extend this algorithm to solve numerically the \ac{NS} equations. After inserting the time series expansion as a representative form of the solution in the equations, the recurrence formula that generates the computation of spatial modes, for both velocity vector $\bu$ and pressure $\p$, is obtained. In the weak form of the recurrence formula case, subsequent terms are obtained by solving a linear system encapsulating a mixed weak formulation of the problem.
A convergence analysis of the numerical method built in the framework of \ac{TSE}-\ac{FEM} has been made and a condition of stability for this technique is established.
Without the stabilization technique, only the term relative to the first order of the velocity expansion is valid. However, applying the stabilization leads to computing valid terms for higher rank, and having a larger stability condition.

The outline of this paper is as follows. \cref{sec2} will present the outcome on the recurrence formula generating terms of the series in the case of \ac{NS} equations and the computation of the terms in the weak formulation. \Cref{sec3} presents a convergence analysis study of the numerical method when applied for the Stokes problem using the \ac{LBB} conditions.
\cref{sec4} presents the time numerical integration technique based on divergent series resummation and Factorial series.
\cref{sec5} presents the Taylor-Green Vortex problem as a study case. \Cref{sec6} presents the stabilization technique to compute the terms of the series in the mixed formulation, ending with the exhibition of Algorithm integrating \ac{NS} equations with \acf{FS}. \cref{sec7} exhibits the performance of \cref{ALgo1} in computing terms and applying the \ac{TSE} algorithm for solving the flow past a cylinder benchmark.
The main conclusions and potential perspectives are outlined in \Cref{sec8} and \ref{sec9}, respectively.

\section{\ac{NS} equations and \ac{TSE} technique}
\label{sec2}

\subsection{Main equation}

We consider an open domain $\gf \subset \mathds{R}^2$ where $\bx = (\x,\y) \in \gf$ is the space variable, $\t$ is the time variable spanning in $ [0,\T{}]$. The incompressible \ac{NS} equations are presented below in their dimensionless form:
\begin{equation}
\label{NS_Eq}
\left\lbrace
\begin{array}{ccl}
 \displaystyle \frac{\partial \bu}{\partial \t} + \left(\bu\cdot\nabla \right)\bu  + \nabla \p &=& \displaystyle \frac{1}{\Re}\Delta \bu   \\
 \nabla\cdot \bu &=& 0
 \end{array}
 \right., \quad
 \forall (\t,\bx) \in [0,\T{}]\times \gf,
\end{equation}
where $\bu \coloneqq (\u,\v)^\top$ is the velocity vector field, $^\top$ is the transpose operator, $\p$ is the scalar pressure, $\Re \coloneqq \rho U_\infty L/\nu$ is the dimensionless Reynolds number,
$\nabla \coloneqq (\frac{\partial }{\partial \x},\frac{\partial }{\partial \y})^\top$ is the gradient operator, $\cdot$ is the symbol scalar product in $\mathds{R}^2$ such that $\nabla\cdot \bu \equiv \frac{\partial \u}{\partial \x}+\frac{\partial \v }{\partial y}$ is the divergence operator of a given velocity vector field $\bu$ and $\Delta  \coloneqq \nabla \cdot \nabla = \frac{\partial^2}{\partial \x^2} + \frac{\partial^2}{\partial \y ^2} $ is the Laplace operator. We end with defining the scalar product element $\bu\cdot \nabla \coloneqq \u\frac{\partial}{\partial \x} + \v\frac{\partial}{\partial \y}$.

The first equation encapsulates the momentum equation and the second one represents the continuity equation or the incompressibility condition. Here $\rho$ is the constant density of the fluid, $\nu$ is the kinematic diffusion parameter, $U_\infty$ and $L$ are respectively the characteristic speed and length.

The system is completed with a given initial condition $\bu(0,\bx)\coloneqq \ru_0(\bx)$, and the following boundary conditions:
\begin{equation}
 \label{bc_NS}
 \begin{array}{cccc}
 \bu(\t,\bx) &=&  \gd(\t,\bx) & \forall \bx \in \gbd ,\\
  \sigma(\bx) \cdot\n \coloneqq\left( \displaystyle\frac{2}{\Re}\strain(\bu)  - \p \I \right)\cdot \n &=&  \gn(\t,\bx) & \forall \bx \in \gbn,
\end{array}
  \end{equation}
with $\sigma(\bx)$ is the stress tensor, $\strain(\bu) \coloneqq (\nabla \bu + (\nabla \bu)^\top)/2$ is the rate of the strain tensor, $\n$ is the outward normal vector to the boundary of $\gf$ and $\I$ is the identity tensor of rank two. Having decomposed the boundary set of $\gf$, denoted by $\partial \gf \coloneqq \gb$, into $\gb= \gbn \cup \gbd$, the first equation in \cref{bc_NS} represents the Dirichlet boundary condition and the second one represents the slip condition. Note that $\gbn \cap \gbd = \emptyset$.

\subsection{\acf{TSE} of the solution}

The solution $\bu$ and $p$ are written in their \ac{TSE} forms, denoted by $\hbu$ and $\hp$ respectively, as follows:
\begin{align}
 \label{up_TSE}
 \hbu(\t,\bx) &\coloneqq \sum\limits_{\k=0}^{\infty} \ru_{\k}(\bx) \,\t^{\k}, & \hp(\t,\bx) &\coloneqq \sum\limits_{\k=0}^{\infty} \rp_{\k}(\bx) \,\t^{\k},
\end{align}
where $\ru_{\k}(\bx)$ and $\rp_{\k}(\bx)$ are spatial functions defined in the fluid domain $\gf$ and represent respectively the $k$\up{th} terms of the series $\hbu$ and $\hp$.
Using relations in \cref{up_TSE}, we can write their partial derivatives as follows:
\begin{align}
\label{TSE_up}
 \frac{\partial \hbu}{\partial \t}(\t,\bx)& =\sum\limits_{\k=1}^{\infty} \k\ru_{\k}(\bx) \,\t^{\k-1}, & \nabla\hp(\t,\bx) &= \sum\limits_{\k=0}^{\infty} \nabla\rp_{\k}(\bx) \,\t^{\k}.
\end{align}
The convective term $(\hbu\cdot\nabla)\hbu$ is obtained in its formal series as follows using the Cauchy product of two power series:
\begin{align}
\label{TSE_convective_Term}
(\hbu(\t,\bx)\cdot\nabla)\,\hbu(\t,\bx)&=\sum\limits_{\k=0}^{\infty} \underbrace{\left(\sum\limits_{m=0}^{\k} \left(\ru_{m}(\bx)\cdot \nabla\right)\ru_{\k-m}\right)}\limits_{w_\k}\t^\k = \sum\limits_{\k=0}^{\infty}  w_\k(\bx)\, \t^\k,
\end{align}
and the Laplace and the divergence terms are expressed in power series as follows:
\begin{align}
\label{TSE_div_laplace}
 \nabla\cdot \hbu\,(\t,\bx) & = \sum\limits_{\k=0}^{\infty} \nabla \cdot \ru_\k(\bx)\, \t^\k , & \Delta \hbu \,(\t,\bx) & = \sum\limits_{\k=0}^{\infty} \Delta \cdot \ru_\k(\bx)\, \t^\k.
\end{align}
We insert elements of \cref{TSE_div_laplace,TSE_convective_Term,TSE_up} in \cref{NS_Eq} and
assembles terms of $\t^k$ together in each side. We obtain:
\begin{equation}
\label{NS_Eq_series}
\left\lbrace
\begin{array}{ccl}
  \sum\limits_{\k=0}^{\infty} \Big[(\k+1)\ru_{\k+1}(\bx) +  w_\k(\bx) + \nabla\rp_{\k}(\bx) \Big]\,\t^{\k} &=&
 \displaystyle \sum\limits_{\k=0}^{\infty} \left[\frac{1}{\Re}\Delta \cdot \ru_\k(\bx) \right]\, \t^\k   \\
 \sum\limits_{\k=0}^{\infty} [\nabla \cdot \ru_\k(\bx)]\, \t^\k &=& 0,
 \end{array}
 \right..
\end{equation}

Equating terms of $\t^\k$ on both sides of the system \eqref{NS_Eq_series}, we find naturally the following recurrence formula that, if having the initial condition $\ru_0$, generates $\ru_{\k+1}$ and $\rp_{\k}$ for $\k\geqslant 0$:
\begin{equation}
\label{NS_Eqk}
\left\lbrace
\begin{array}{lcl}
(\k+1)\,\ru_{\k+1} + \nabla \rp_{\k} &=&\rF_{\k}  \\
 \nabla\cdot \ru_{\k+1} &=& 0
 \end{array}
 \right.,\quad
 \forall \bx \in \gf, \forall \k \in \mathds{N},
\end{equation}
where the right hand side $\rF_{\k}$ is given as follows:
\begin{equation}
 \label{rhs_NS_Eqk}
 \rF_{\k} \coloneqq -\displaystyle\sum\limits_{m=0}^{\k} \left(\ru_m\cdot\nabla \right)\ru_{\k-m} +  \frac{1}{\Re}\Delta \ru_{\k}.
\end{equation}
This is promising as this process transform the non-linear \ac{NS} \cref{NS_Eq} into a cascade of linear equations. This linearisation is numerically more efficient when comparing with implicit numerical methods as the \ac{BDF} schemes in the framework of \ac{FEM}. The efficiency comes from the fact that the Mass matrix will be assembled and inverted only once during the whole simulation, while in implicit methods the mass matrix has to be assembled and inverted at every time step. This reduces the cost time computation and promise to accelerate the computational process.

As for the boundary conditions, functions $\gd$ and $\gn$ are written in their \ac{TSE} forms:
\begin{align}
 \gn(\t,\bx) & =\sum\limits_{\k=0}^{\infty} \rf_{\k}(\bx) \,\t^{\k}, &
 \gd(\t,\bx) & =\sum\limits_{\k=0}^{\infty} \rg_{\k}(\bx) \,\t^{\k},
\end{align}
thus the system \cref{rhs_NS_Eqk} is completed by adding the following boundary conditions in the determination of $\ru_{\k+1}$  and $\rp_\k$:
\begin{equation}
 \label{bck_NS}
 \begin{array}{cccc}
 \ru_{\k+1}(\bx) &=&  \rg_{\k+1}(\bx), & \forall \bx \in \gbd ,\\
  \left( \frac{2}{\Re}\strain(\ru_\k)  - \rp_\k \I \right)\cdot \n &=&  \rf_\k(\bx), & \forall \bx \in \gbn,
\end{array}
, \quad \forall \k \in \mathds{N}
  \end{equation}
To this end, we can see that the governing equations \cref{NS_Eqk} contain the term $\nabla\rp_k$ which means that the equation is solved by finding $\rp_k$ up to an arbitrary function of time. To fix this degree of freedom, we add the following condition:
\begin{equation}
 \label{cond_pk}
 \int_\gf \rp_\k \, \d\bx = 0.
\end{equation}
If one has the initial condition, \emph{i.e.} the term $\ru_0$ relative the $0$\up{th} rank of $\hbu$, the first term $\ru_1$ is computed by solving System \eqref{NS_Eqk} for $\k=0$ with the associated boundary conditions and so on till computing $(\ru_{\k+1},\rp_\k)$ up to any required rank $\m$. Next, we present the variational formulation of System \eqref{NS_Eqk}, preparing the ground for the formulation of \ac{FEM}.

\subsection{Mixed formulation}

Defining $\V \coloneqq \left( H^1 (\gf)\right)^2$ the space containing velocity field vector $\bv$ of two components defined on domain $\gf$, and consider the space $\V_\k$:
$$\V_\k = \left\lbrace \bv \in \V  \,\Big\vert \left. \bv\right\vert_{\gbd} = \rg_\k \right\rbrace.$$
We denote by $\bev_\ell \coloneqq (\psi_\ell,\psi_\ell)$  a basis' element of $\V_\k$ with $\psi_\ell \in H^1 (\gf)$. Consider also the space $\P$ such that:
$$\P \coloneqq \left\lbrace\q \in L^2(\gf) \left\vert \int_\gf \q \,\d\bx = 0 \,\wedge \,\left(\nabla \q\cdot \n\right)\Big\vert_{\gbd} = 0\right.\right\rbrace,$$
where $\wedge$ is the logical and, and the space $\P$
containing scalar functions $\q$ defined on $\gf$. 
We denote by $\bep_\ell$ a basis element of this space. For any function $\bv \in \V_\k$ and $\q \in \P$, they are represented by linear combination of basis elements as follows:
\begin{align}
 \bv(\bx) &\coloneqq \sum\limits_{\ell} \clv_{\ell} \bev_\ell(\bx),&
 \q(\bx) &\coloneqq \sum\limits_{\ell} \clq_{\ell} \phi_\ell(\bx),
\end{align}
where $\clv_{\ell} \in \mathds{R}^2$ and $\clq_{\ell}\in\mathds{R}$ . To establish the weak formulation of System \eqref{NS_Eqk}, we multiply it by the couple of test functions $(\bv, \q) \in \V_\k\times \P$ and integrate over the domain $\gf$. We get the following mixed formulation of the problem. We are seeking a weak solution $(\ru_{\k+1},\rp_k) \in \V_\k \times \P$ such that:
\begin{equation}
 \label{NS_eqk_vf}
 \left\lbrace
 \begin{array}{rcl}
  a_{\k+1}(\ru_{\k+1},\bv) +b(\bv,\rp_\k) &=& \lin_\k(\bv)\\
  b(\ru_{\k+1},\q) &=& 0
 \end{array}
 \right., \quad
 \forall \,(\bv,\q)\in \V_\k \times \P,
\end{equation}
where the forms $a_{\k+1}(\cdot,\cdot)$, $b(\cdot,\cdot)$ and $\lin_\k(\cdot)$ are defined as follows:
\begin{align}
 a_{\k}(\bu,\bv) &\coloneqq \k \int_\gf \bu\cdot\bv \,\d\bx, &
 b(\bv,\q) & \coloneqq -\int_\gf \q\, \nabla\cdot\bv \,\d\bx,
 &
 \lin_\k(\bv) &\coloneqq \displaystyle\int_\gf \rF_\k\cdot \bv \, \d\bx.
\end{align}
The existence and uniqueness of the solution of \cref{NS_eqk_vf} could be studied throughout the use of the Brezzi's conditions \cite{brezzi-75}.
Note that, if $a_\k(\bu,\bv) \equiv \int_{\gf} \nabla \bu :\nabla\bv\,\d\bx$, we found the steady Stokes problem such that $A:B$ is the second order tensors contraction operator.
The first Brezzi's condition demands to ensure the so-called $\inf-\sup$ condition of the bilinear form $b$, stating that there exists a constant $\beta>0$ such that:
\begin{equation}
 \label{inf-sup-cond}
 \inf\limits_{0\neq \q \in \P} \sup\limits_{0\neq \bv \in \V_\k} \frac{b(\bv,\q)}{\|\bv\|_{\V_\k}\,\|\q\|_{\P}} \geqslant \beta.
\end{equation}
The second condition states the coercivity of bilinear operator $a_\k$:
\begin{equation}
 \label{coercive_a}
 a_\k(\bu,\bu) \geqslant \delta_\k \|\bu\|^2_{\Z}, \quad \forall \bu \in \Z.
\end{equation}
This condition is verified too here for $\delta_\k=\k \delta$, $\delta$ is the coercivity constant of $a_1$, and $\Z  = \{\bu \in \V_\k \vert b(\bu,\q)=0\,, \forall \q \in \P\}$.

\subsection{\ac{FEM} framework}
We denote by $\Pk{r}[{\bx}]$ the set of scalar polynomials defined on $\gf$ of degree $r$, and by $\Th$ a triangulation of the domain $\gf$. We denote by $\K_i, i\in\{1,\ldots \Ne\}$ the simplices that forms this triangulation, and by $\bx^j, j\in\{1,\ldots,\Nv\},$ the vertices of the resulting simplices forming the set of nodes on which the solution will be approximated. We define the set of functions $\bep_\ell: \gf\mapsto \mathds{R}$ such that:
\begin{align}
 \label{basisp-FE}
 \bep_\ell\vert_{\K_i} \in \Pk{r}, &&  \int_\gf \bep_\ell \cdot \bep_{\ell^{\prime}}\,\d\bx &\equiv \delta_{\ell}^{\ell^{\prime}},
 & \bep_\ell(\bx^j)&= \delta_{\ell}^{j}, & \forall \ell,\ell^{\prime},i,j\in\{1,\ldots \Nv\},
\end{align}
where $\delta_{\ell}^{\ell^{\prime}}$ is the standard Kronecker symbol. The set of vector-functions polynomials $\bev_\ell: \gf\mapsto \mathds{R}^2$ of degree $s$ verifies also the following  conditions:
\begin{align}
 \label{basisv-FE}
 \bev_\ell\vert_{\K_i} \in (\Pk{s}[\bx])^2, &&  \int_\gf \bev_\ell \cdot \bev_{\ell^{\prime}}\,\d\bx &\equiv \delta_{\ell}^{\ell^{\prime}},
 & \bev_\ell(\bx^j)&\equiv \boldsymbol{\delta}_{\ell}^{j}, & \forall \ell,\ell^{\prime},i,j\in\{1,\ldots \Nv\},
\end{align}
where $\boldsymbol{\delta}_{\ell}^{j} \in \mathds{R}^2$ taking at each components the value ${\delta}_{\ell}^{j}$ at the nodes $\left\lbrace \bx^j\right\rbrace_{j=1}^{\Nv}$.

For velocity field functions, the subspace $\V^s_{h,\k}\subseteq \V_\k$ of vector functions that are continuous and locally polynomials of order $s$ on the triangulation is considered as finite space, while the subspace  $\P^r_h\subseteq \P$ of continuous functions and locally polynomials of order $r$ is considered for the pressure element. We seek to find approximation $(\ru_{\k+1,h}, \rp_{\k,h}) \in \V^s_{h,\k}\times\P^r_h$ for $(\ru_{\k+1}, \rp_\k)$
Here $h$ stands for the mesh size that we define as the highest inscribed circle of all simplices of the given mesh, $r$ and $s$ are the degrees of space approximations of the pressure and velocity, respectively. To ensure the $\inf-\sup$ condition, the degree $s$ should be strictly greater than $r$ ($s>r$). Verifying also the coercivity conditions of $a(\cdot,\cdot)$ in the following space:
$$\Z_h \coloneqq \{ \bu \in \V^s_{h,\k} \vert b(\bu,\q) = 0 \, \forall \q_h \in  \P^r_h\}$$
will complete the Brezzi's conditions in order to converge to a unique weak solution when the size of mesh $h\rightarrow 0$.

We return now to the system \cref{NS_Eqk}. Consider a given triangulation $\Th$ and the finite elements spaces are constructed. To find the approximations $(\ru_{\k+1,h}, \rp_{\k,h}) \in \V^s_{h,\k}\times \P^r_h$, we multiply System \eqref{NS_Eqk} by their basis elements $(\bev_{\ell}, \bep_{\ell}) \in \V^s_{h,\k}\times \P^r_h$, integrate over $\Th$ and assemble the discrete system. The following linear system is obtained:
\begin{equation}
 \label{mix_form_FE}
 \left(
 \begin{array}{cl}
  \mA_{\k+1} & \mB^\top \\ \mB&0
 \end{array}
 \right)
  \left(
 \begin{array}{cc}
  \mU_{\k+1,h} \\ \mP_{\k,h}
 \end{array}
 \right) =
   \left(
 \begin{array}{cc}
  \mF_{\k,h} \\ 0
 \end{array}
 \right) ,
\end{equation}
where $\mA_{\k+1}$ and $\mB$ are matrices defined as follows:
\begin{align}
 \left(\mA_{\k}\right)_{\ell,\ell^{\prime}} &\coloneqq \k \int_{\gf} \bev_{\ell} \cdot \bev_{\ell^{\prime}} \d\bx , &
 \left(\mB\right)_{\ell,j}  &\coloneqq \int_{\gf}\, \bep_{j}\, \nabla\cdot \bev_{\ell} \, \d\bx.
\end{align}
Here, $\mU_{\k+1,h}$ and $\mP_{\k,h}$ are vectors representing the approximation to the solutions $\ru_{\k+1}$ and $\rp_\k$ on the nodes $\bx^j$ of the mesh $\Th$:
\begin{equation}
 \left(\mU_{\k,h}\right)_j \coloneqq \ru_{\k,h}(\bx^j)  , \quad
 \left(\mP_{\k,h}\right)_j \coloneqq \rp_{\k,h}(\bx^j),
\end{equation}
thus the functions $\ru_{\k,h}$ and $\rp_{\k,h}$ defined as below:
\begin{align}
 \ru_{\k,h}(\bx) & = \sum\limits_\ell\left(\mU_{\k,h}\right)_\ell\bev_{\ell}(\bx),&
\rp_{\k,h}(\bx) & = \sum\limits_\ell\left(\mP_{\k,h}\right)_\ell\bep_{\ell}(\bx),
\end{align}
and $\mF_{\k,h}$ is the vector representing the projection of $\rF_\k$ on the subspace $\V_{\k,h}^s$ via its finite elements basis.


\section{Convergence analysis for the Stokes problem}
\label{sec3}
In this section, we will present a preliminary result for convergence of the numerical method in the framework of \ac{TSE} and \ac{FEM} for the case of unsteady Stokes problem:
\begin{equation}
\label{S_Eq}
\left\lbrace
\begin{array}{ccl}
 \displaystyle \frac{\partial \bu}{\partial \t} + \nabla \p &=& \displaystyle \frac{1}{\Re}\Delta \bu   \\
 \nabla\cdot \bu &=& 0
 \end{array}
 \right., \quad
 \forall (\t,\bx) \in [0,\T{}]\times \gf.
\end{equation}
For the boundary conditions, consider the case of homogeneous boundary condition of the velocity. If it is not the case, a variable change could be applied. The modals space $\ru_\k(\bx)$ of the \ac{TSE}, $\hbu(\t,\bx)$, are generated by the following recurrence formulas:
\begin{equation}
\label{S_Eqk}
\left\lbrace
\begin{array}{lcl}
(\k+1)\,\ru_{\k+1} + \nabla \rp_{\k} &=&\frac{1}{\Re} \Delta \ru_\k  \\
 \nabla\cdot \ru_{\k+1} &=& 0
 \end{array}
 \right.,\quad
 \forall \bx \in \gf, \forall \k \in \mathds{N}.
\end{equation}

\subsection{Variational formulation and Brezzi's condition}
For all $(\bv, \q) \in \Vr\times \P$, $\Vr$ is the space of functions that vanish on boundary $\partial \gf$, the variational form is written as follows:
\begin{equation}
 \label{S_eqk_vf}
 \left\lbrace
 \begin{array}{rcl}
  a_{\k+1}(\ru_{\k+1},\bv) +b(\bv,\rp_\k) &=& \lins_{\ru_{\k}}(\bv)\\
  b(\ru_{\k+1},\q) &=& 0
 \end{array}
 \right., \quad
 \forall \,(\bv,\q)\in \Vr \times \P,
\end{equation}
with $a_\k$ and $b$ are the bilinear operators defined above and $\lins_\k$  is the following linear operator obtained after using the divergence theorem and the fact the $\bv\in \Vr$:
\begin{align}
 \lins_{\ru_{\k}} (\bv)= \frac{1}{\Re}\int_\gf \nabla\ru_\k : \nabla \bv \,\d \bx.
\end{align}
Having $a_\k$ a continuous coercive bilinear operator with coercivity constant $\delta_\k = \delta \k$, $b$ a bilinear continuous operator fulfilling the inf-sup condition for the desired space $\Vr\times \P$ and $\lins_{\ru_{\k}}$ is a linear continuous operator, we can apply the Brezzi's conditions to prove the existence and uniqueness of solution $(\ru_{\k+1}, \rp_\k) \in \Vr\times \P$ for the mixed formulation \eqref{S_eqk_vf}. In addition, from Brezzi's conditions we have the following estimation:
\begin{align}
 \label{est_S_vp}
 \begin{aligned}
 \left\| \ru_{\k+1}\right\|_{\Vr} &\leqslant \frac{1}{\delta_{\k+1}} \left\| \lins_{\ru_{\k}}\right\|_{\left(\Vr\right)^\prime} \leqslant \frac{1}{\delta_{\k+1}\Re} \left\| \ru_\k\right\|_{\Vr},\\
 \left\| \rp_{\k}\right\|_{\P} &\leqslant \frac{2M_{\k+1}}{\delta_{\k+1}\beta}\left\| \lins_{\ru_{\k}}\right\|_{\left(\Vr\right)^\prime} \leqslant \frac{2M_{\k+1}}{\delta_{\k+1}\beta\Re} \left\| \ru_\k\right\|_{\Vr},
 \end{aligned}
\end{align}
with $\delta_\k$ is the coercivity constant of $a_\k$, $M_\k$ is its continuity constant and $\beta$ is the inf-sup constant condition. The proof of inequalities \eqref{est_S_vp} is a classical result of Brezzi (See the theorem in \ref{append1} or \cite{brenner-10}).
Using the formula of $a_\k$, we have $\delta_\k = \delta \cdot\k$ and $M_\k = M\cdot \k$ where $\delta>0$ and $M>0$ are the coercivity and continuity constants of $a_1$. This give us the following bounds:
\begin{align}
\label{bounds_S}
\begin{aligned}
\left\| \ru_{\k}\right\|_{\Vr} &\leqslant \frac{1}{\delta^\k\Re^\k \k!}   \left\| \ru_0\right\|_{\Vr}, \quad\forall \k>1,\\
\left\| \rp_{\k}\right\|_{\P} &\leqslant \frac{2M}{\beta\delta^{\k+1}\Re^{\k+1}(\k+1)!} \left\| \ru_0\right\|_{\Vr},\quad \forall \k>0.
\end{aligned}
\end{align}
This means that the norm of spatial modes are bounded by the norm of the initial condition, the coercivity constant and the Reynolds number. The series
\begin{equation}
 \sum\limits_{\k=0}^{\infty} \left\|\ru_\k(\bx)\right\|_{\Vr} \t^\k,
\end{equation}
is an entire function in the complex plane.

\subsection{Error of \ac{TSE}-\ac{FEM} approximation}
In this section, we will estimate the error of the numerical approximations $(\ru_{\k+1,h} ,\rp_{\k,h}) \in \Vrh^s \times \P_h^r$ such that the subspace $\Vrh^s \subset \Vr$ contains piecewise polynomials of order $s$ and $\P_h^r\subset \P$ contains piecewise polynomials of order $r$  approximating, respectively, the velocity $\ru_{\k+1}$ and  pressure $\rp_\k$ for the problem \eqref{S_eqk_vf}. Note that we should have $r<s$ to ensure the inf-sup Brezzi's condition.

Numerically, the spatial mode of rank $\k+1$ is approximated using the last modes that already have been approximated in the finite space, namely $\ru_{\k,h}$. For the Stokes problem, the approximation $(\ru_{\k+1,h},\rp_{\k,h})$ is the solution of the following mixed problem:
\begin{equation}
\left\lbrace
 \begin{array}{rcl}
  a_{\k+1}(\ru_{\k+1,h},\bv_h) +b(\bv_h,\rp_{\k,h}) &=& \lins_{\ru_{\k,h}}(\bv_h) \\
  b(\ru_{\k+1,h},\q_h) &=& 0
 \end{array}
 \right., \quad
 \forall \,(\bv_h,\q_h)\in \Vrh^s \times \P_h^r.
\end{equation}
By adding and subtracting the exact $\ru_\k$ in $\lins_{\ru_{\k,h}}$, we have:
\begin{align}
 \begin{aligned}
\lins_{\ru_{\k,h}}(\bv_h)  & = \lins_{\ru_{\k}}(\bv_h) + \lins_{\ru_{\k,h}-\ru_{\k}}(\bv_h).
\end{aligned}
\end{align}
First, we neglect the numerical error in $\ru_{\k,h}$ and  we use the exact solution of $\ru_\k$ in the approximation of $\ru_{\k+1}$, denoted by $\reu_{\k+1,h}$, \emph{i.e.} $ \lins_{\ru_{\k,h}-\ru_{\k}} \equiv 0$. Using a similar way of Céa's Lemma \cite{Cea-64},  we can have the following estimation:
\begin{align}
 \label{error_Vh}
 \begin{aligned}
 \left\| \reu_{\k+1} - \ru_{\k+1,h} \right\|_{\Vr} & \leqslant \frac{4M_{\k+1} M_b}{\delta_{\k+1} \beta_h} E_{\ru_{\k+1}} + \frac{M_b}{\delta_{\k+1}} E_{\rp_\k} ,
 \\
 \left\| \rep_{\k} - \rp_{\k,h} \right\|_{\P}& \leqslant \frac{4M^2_{\k+1} M_b}{\delta_{\k+1} \beta^2_h} E_{\ru_{\k+1}} + \frac{3M_{\k+1}M_b}{\delta_{\k+1}\beta_h} E_{\rp_\k},
 \end{aligned}
\end{align}
with $M_b$ is the continuity constant of $b(\cdot,\cdot)$ and $E_{\ru_{\k+1}}$ and $E_{\rp_\k}$ are the best approximation errors of $\ru_{\k+1}$ and $\rp_{k}$ in $\Vrh^s$ and $\P_h$ respectively:
\begin{align}
\label{Eu_Deltax}
\begin{aligned}
E_u &= \inf_{u_I\in \Vrh^s}\|u-u_I\|_{\Vr}\approx O(h^s),& E_p &= \inf_{p_I\in \P_h^r}\|p-p_I\|_\P\approx O(h^r).
\end{aligned}\end{align}
The proof of bounded error in \cref{error_Vh} when using the exact solution of $\ru_\k$ in computing $\reu_{\k+1,h}$ is presented in \ref{append21}. Now, we consider again the error propagation in $\ru_{\k,h}$. By continuity of $\lins:\Vr \mapsto (\Vr)^\prime$ and using the error bounds in \eqref{error_Vh}, we have
\begin{equation}
\left\| \lins_{\ru_{\k,h}-\ru_{\k}}\right\|_{(\Vr)^\prime} \leqslant \frac{1}{\Re}\left\| \reu_{\k,h}-\ru_{\k}\right\|_{\Vr} \leqslant \frac{4M M_b}{\Re\,\delta \beta_h} O(h^s) + \frac{M_b}{\k\Re\,\delta} O(h^r)
\end{equation}
Thus, we repeat the same process by incorporating $\lins_{}$ to find the following bound error:
\begin{equation}
\label{error_Vhc}
 \begin{aligned}
  \left\|\ru_{\k+1}-\ru_{\k+1,h} \right\|_{\Vr} &\leqslant
  \frac{4MM_b}{\delta\beta_h}\left[1 + \frac{1}{(\k+1)\delta\Re} \right]O(h^s) \\
  &\quad+ \frac{M_b}{(\k+1)\delta} \left[ 1+\frac{1}{\k\delta \Re}\right] O(h^r)\\
  \left\|\rp_{\k}-\rp_{\k,h} \right\|_{\Vr} & \leqslant \frac{4M^2M_b}{\delta\beta_h^2}\left[\k+1 + \frac{2}{\delta\Re} \right]O(h^s) \\
  &\quad+ \frac{2MM_b}{\delta\beta_h} \left[ 1+\frac{1}{\k\delta \Re}\right] O(h^r)\\
 \end{aligned}
\end{equation}
To see the full proof, please refer to \ref{append2}.

\subsection{Error Convergence}

For a fixed value of $\t = \Delta \t$, we have $\bu(\Delta \t,\bx) \in \Vr$. The approximation in the framework of \ac{TSE}-\ac{FEM} given by the following partial sum also belongs to $ \Vr$:
\begin{equation}
\label{TSE-FEM}
  \sum\limits_{\k=0}^{\m} \ru_{\k,h}(\bx)\, \Delta \t^\k \in \Vr.
\end{equation}
Thus the error between the exact and the one built in the framework of \ac{TSE}-\ac{FEM} is estimated as follows:
\begin{align}
 \begin{aligned}
  e_\m(\Delta \t,h)
  &\coloneqq \left\| \bu(\Delta \t,\bx) - \sum\limits_{\k=0}^{\m} \ru_{\k,h}(\bx)\, \Delta\t^\k \right\|_{\Vr}\\
 & = \left\| \sum\limits_{\k=0}^{\infty} \ru_\k(\bx) \Delta\t^\k - \sum\limits_{\k=0}^{\m} \ru_{\k,h}(\bx)\, \Delta\t^\k \right\|_{\Vr}
  \end{aligned}
\end{align}
which is smaller than:
\begin{align}
 \begin{aligned}
  &\leqslant\left\| \ru_{0} - \ru_{0,h} \right\|_{\Vr} + \sum\limits_{\k=1}^{\m} \left\| \ru_{\k+1} - \ru_{\k+1,h} \right\|_{\Vr} \Delta\t^\k \\
  & \quad+\sum\limits_{\k=\m+1}^{\infty} \left\|\ru_\k(\bx)\right\|_{\Vr} \Delta\t^\k
   \end{aligned}
\end{align}
The first term in the last inequality represents the error of the projection of the initial condition to the finite space which is in this case is of order $O(h^s)$, or the accumulation of error when reaching such an instant time $\t$. Consider we are at $\t=0$, using inequalities \eqref{bounds_S} and \eqref{error_Vhc} we obtain:
\begin{align}
\label{e_Nhtau}
 \begin{aligned}
   e_\m(\Delta \t,h) &\leqslant O(h^s) +\sum\limits_{\k=1}^{\m}\left(\frac{4M M_b}{\delta \beta_h} \cdot O(h^s) + \frac{M_b}{\delta({\k+1})} \cdot O(h^r) \right)\cdot  \Delta\t^k \\
  &\quad  +\sum\limits_{\k=1}^{\m}\left(\frac{4M M_b}{(\k+1)\delta^2 \beta_h\Re} \cdot O(h^s) + \frac{M_b}{\k({\k+1})\delta^2\Re} \cdot O(h^r) \right)\cdot  \Delta\t^k\\
   &\qquad +   \left\| \ru_0\right\|_{\Vr}
   \sum\limits_{\k=\m+1}^{\infty}\frac{\Delta \t^\k}{\delta^\k\Re^\k \k!}
 \end{aligned}
\end{align}
where the second term in the right hand side is a geometrical series that is estimated as follows:
\begin{equation}
 e_\m^1(\Delta \t,h) = \sum\limits_{\k=1}^{\m}\frac{4M M_b}{\delta \beta_h} \cdot O(h^s)\cdot  \Delta\t^k = \frac{4M M_b}{\delta \beta_h}\cdot O(h^s)\cdot \Delta \t\cdot\left(\frac{1-\Delta \t^{\m}}{1-\Delta \t}\right).
\end{equation}
Having the fact that $\beta_h \xrightarrow[h\rightarrow 0]{} \beta >0$  and $\displaystyle \frac{M}{\delta} = \kappa(\mA_1) \sim O(\frac{1}{h})$ (in the discrete form), we can conclude that its limit, for a fixed $\m$, w.r.t $\Delta \t$ and $h$ is:
\begin{equation}
\lim\limits_{h \rightarrow 0} \lim\limits_{\Delta \t \rightarrow 0} e_\m^1(\Delta \t,h) =  \lim\limits_{h \rightarrow 0}   \frac{4 M_b}{ \beta_h} \cdot O(h^{s-1}) \cdot \lim\limits_{\Delta \t \rightarrow 0}\left[ \Delta \t \cdot \frac{1-\Delta \t^{\m}}{1-\Delta \t} \right] = 0.
\end{equation}
The second term in \cref{e_Nhtau} could be approximated as follows:
\begin{equation}
 e_\m^2(\Delta \t,h) = \sum\limits_{\k=1}^{\m} \frac{M_b}{\delta({\k+1})} \cdot O(h^r)\cdot  \Delta\t^k = \frac{M_b \,O(h^r)}{\delta}\sum\limits_{\k=1}^{\m} \frac{\Delta \t^\k}{\k+1}
\end{equation}
where
\begin{equation}
 \sum\limits_{\k=1}^{\m} \frac{\Delta \t^\k}{\k+1}\leqslant  -\frac{\ln(1-\Delta\t)}{\Delta \t} -1 \xrightarrow[\Delta \t \rightarrow 0]{} 0,
\end{equation}
thus
\begin{equation}
 \lim\limits_{h \rightarrow 0} \lim\limits_{\Delta \t \rightarrow 0} e_\m^2(\Delta \t,h) \leqslant  \lim\limits_{h \rightarrow 0} \frac{M_b \,O(h^r)}{\delta}\cdot \lim\limits_{\Delta \t \rightarrow 0} \left[-\frac{\ln(1-\Delta\t)}{\Delta \t} -1\right]= 0.
\end{equation}
The third term in \cref{e_Nhtau} follows the same estimation as the second to obtain:
\begin{equation}
 e_\m^3(\Delta \t,h) = \sum\limits_{\k=1}^{\m} \frac{4M M_bO(h^s)}{(\k+1)\beta_h\delta^2 \Re} \Delta \t^\k \leqslant \frac{4M M_bO(h^s)}{\beta_h\delta^2 \Re}\left[-\frac{\ln(1-\Delta\t)}{\Delta \t} -1\right],
\end{equation}
which has the following limits when $\Delta \t \rightarrow 0$:
\begin{equation}
\lim\limits_{\Delta \t \rightarrow 0} e_\m^3(\Delta \t,h) \leqslant \frac{4M M_bO(h^s)}{\beta_h\delta^2 \Re} \sim \frac{4 M_bO(h^{s-1})}{\beta_h\delta \Re} \xrightarrow[h\rightarrow 0]{}0.
\end{equation}
The fourth term in the inequality has the following limit:
\begin{equation}
 \lim\limits_{h \rightarrow 0} \lim\limits_{\Delta \t \rightarrow 0}\sum\limits_{\k=1}^\m \left(\frac{M_b}{\k({\k+1})\delta^2\Re} \cdot O(h^r) \right)\cdot  \Delta\t^k = 0.
\end{equation}
For a bounded initial condition, the last term in \cref{e_Nhtau} has the following sum:
\begin{equation}
 \left\| \ru_0\right\|_{\Vr}
   \sum\limits_{\k=\m+1}^{\infty}\frac{\Delta \t^\k}{\delta^\k\Re^\k \k!} = \left\| \ru_0\right\|_{\Vr} \exp\left({\frac{\Delta \t}{\delta \Re}}\right) \left[1- \frac{\Gamma(\m+1,\frac{\Delta \t}{\delta \Re})}{\Gamma(\m+1)}\right],
\end{equation}
with $\Gamma$ is the Gamma function verifying $\Gamma(\m+1) = \m \Gamma(\m)$. We have:
\begin{equation}
\Gamma(\m+1,\frac{\Delta \t}{\delta \Re})\xrightarrow[\Delta \t \rightarrow 0]{}\Gamma(\m+1),
\end{equation}
and finally for a fixed $\m$ we have:
\begin{equation}
 \lim\limits_{h \rightarrow 0} \lim\limits_{\Delta \t \rightarrow 0} e_\m(\Delta \t,h) = 0,
\end{equation}
which completes study of convergence of the \ac{TSE}-\ac{FEM} numerical scheme in the Stokes problem.

\subsection{Discussions}
In the weak formulation of the Stokes problem, the product $\delta \cdot\Re$ is important in how the partial sum is efficient in approximating the weak solution. If $\delta \cdot\Re >1$, \emph{i.e.} the coercivity constant of $a_1$ is bigger than the inverse of the Reynolds number ($\delta > \frac{1}{\Re}$), the partial sum
\begin{equation}
  \sum\limits_{\k=0}^{\m} \ru_\k(\bx) \t^\k
\end{equation}
does not present any problem in approximating the solution $\bu(\t,\bx)$ up to any order $\m$ for a fixed $\t$. However, if $\delta < \frac{1}{\Re}$, the partial sum will not be efficient at the beginning, \emph{i.e.} the error of the partial sum with the entire function will increase up to $\m\approx \cfrac{\t}{\delta \cdot \Re}$ and then decreases for $\k>\m$.

On the other hand and in finite space, the dual space of $\Vrh$ is itself meaning that the application $\ru_h\mapsto a_1(\ru_h,\cdot) \in \Vrh$, thus the continuity constant $M$ of $a_1$ verifies the following inequalities:
\begin{align}
 \lambda_{\max}&\leqslant  M ,&\lambda_{\min}&\geqslant \delta
\end{align}
where $\lambda_{\max}$ and $\lambda_{\min}$ are the largest and smallest eigenvalues of $\mA_1$. This gives us their relation of  coercivity and continuity constants with the condition number $\kappa(\mA_1)$ of matrix $\mA_1$ as follows:
\begin{equation}
 \kappa(\mA_1) = \frac{\lambda_{\max}}{\lambda_{\min}}\leqslant \frac{M}{\delta} .
\end{equation}

By checking the error $e_\m(\Delta \t, h)$ in \cref{e_Nhtau}, you can notice that the following elements plays major role in the convergence error
\begin{itemize}
 \item The error $O(h^s)$ of the velocity and $O(h^r)$ of the pressure.
 \item The ratio $\displaystyle\frac{M}{\delta}$ in $e_\m^1(\Delta \t,h)$ which is represented by the condition number $\kappa(\mA_1)$ of the Mass matrix in the discrete form. The smaller the ratio is, the smaller is its contribution in the global error; \emph{i.e.} having a bigger coercivity constant implies that the error is smaller.
 \item The ratio $\displaystyle\frac{M_b}{\beta_h}$ represents how the inf-sup condition is fulfilled in the discrete form. the bigger the constant $\beta_h>0$ we find verifying the inf-sup condition, the smaller the error of the first part is.
\end{itemize}

In the next section, we will state the resummation technique for the time series expansion, namely the Borel-Laplace summation. The sum is presented effectively using the factorial series for what it is suitable for the functions rather than the \ac{BPL} algorithm which is for moment tailored for pointwise evaluation.

\section{Resummation technique for divergent series}
\label{sec4}
In this section, we will start by presenting the \ac{BL} resummation for series solution of \ac{PDE}. Then the \ac{FS} algorithm is presented as he has the advantage for dealing with elements in Banach spaces. We end by presenting the process of numerical flow in the concept of continuation.

\subsection{Borel Laplace Resummation}

In the past work \cite{DEEB_2022_bpl}, the \ac{BL} resummation technique has been used in the framework of \ac{ODE} to develop numerical integrators. When having \ac{PDE}, we consider the cases of solutions $\bu(\t,z)$ such that $z$ is the complex variable, and coefficients of the series, $\ru_\k(z)$ are elements in a Banach space $\mathbb{E}_r$ of functions that are holomorphic in the disc $D_r$ or radius $r$ such that there exists constants $\mathrm{A}$ and $\mathrm{C}$:
$$\left\| \ru_\k(z) \right\|_{\mathbb{E}r} < \mathrm{C}\mathrm{A}^\k \k!, \quad \forall \k>0.$$
In this case, we have a Gevrey series of order $1$.
For more details, please refer to \cite{balser-99}.
Having the series
$$\hbu(t,z) = \sum\limits_{k=0}^{\infty}\ru_\k(z) \t^\k,$$
we define its Borel transform with respect to its variable $\t$ by
$$ \big(\mathcal{B}\hbu(t,z) \big)(\xi,z) =\sum_{\k=0}^{\infty}\cfrac{\ru_{\k+1}(z)}{\k!}\ \xi^\k.$$
If the above series converges absolutely (which is the case if we have a Gevrey series of order $1$) for $|z|<r$ and $|t|<R$, for $R>0$, to a function that could be continued by $P(\xi,z)$ with respect to $\xi$ into the sector $S_{\gamma} = \left\lbrace \xi: |\arg \xi|< \gamma\right\rbrace$. Moreover, for evert $\gamma_1<\gamma$ there exists constants $C$ and $K$ such that:
$$\left\|P(\xi,.)\right\|_{\mathbb{E}_r} \leqslant C \exp(K|\xi|),\quad \forall \xi \in S_{\gamma_1},$$
then the Laplace transform of $P(\xi,z)$,
$$\ru_0(z) + \int_0^{\infty} P(\xi,z)\, e^{-\xi/\t} \d\xi,$$
integrating along the positive real axis is the Borel sum of the formal series $\hbu(\t,z)$. This is an holomorphic function over the domain $S_{\pi}\times D_r$.
We resume these steps in the below diagram.
\captionsetup[table]{name=Diagram,skip=1ex, labelfont=bf}
\begin{table}[!ht]
\caption{Schematic of \ac{BL} resummation for }
\centering  {\sf\small
	\begin{tabular}{clc}
		$\displaystyle\hbu(\t,z) =\sum_{\k=0}^{\infty}\ru_\k(z) \t^\k$& $\underleftarrow{\hspace{0.5cm}\sim_{\mathbf{V}}\hspace{0.5cm}}$   &\hspace{-1cm}$\mathcal{S} \hbu(\t,z) \coloneqq\ru_0(z)+\displaystyle
	\int_{0}^{\infty} P(\xi,z) e^{-\xi/\t}\operatorname d\xi$
	\\\\
	 $\left.\begin{array}{c}\\\text{Borel}\\\\\end{array}\right\downarrow$&
	 &$\left\uparrow \begin{array}{c}\\\text{Laplace} \\\\\end{array}\right.$
	\\\\
	$\displaystyle\mathcal{B}\hbu(\xi,z) =
	\sum_{\k=0}^{\infty}\cfrac{\ru_{\k+1}(z)}{\k!}\ \xi^\k$
	&$\overrightarrow{\hspace{.5cm}\text{Prolongation}\hspace{.5cm}}$
	&\hspace{-1cm}$P(\xi,z)$
	\end{tabular}}
	\label{schema-bl}
\end{table}

When having \ac{PDE}, the \ac{BPL} algorithm could be developed in the framework of point wise once the space is discretized. This is because the Padé-approximant represents a challenge when dealing with function space. However, point values are not well suited. This is why we will use an effective technique to represent the \ac{BL} using the \ac{FS}.

\subsection{Effective Borel resummation using Factorial series}
Here, we will develop the process of resummation technique for divergent series using \ac{BL} resummation. Instead of using the \ac{BPL} algorithm \cite{ahmad_bpl_2014,DEEB_2022_bpl}, we will employ the one based on the \ac{FS}. Note that Delabaere \emph{et al.} \cite{Delabaere-07} have shown, in the case of having scalar coefficients $\ru_\k$, how the \ac{BL} resummation is obtained via \ac{FS} using symbolic mathematics:
\begin{equation}
 \label{FS-BL}
 \ru_0 + \int_0^\infty \left(\sum\limits_{\k=0}^{\infty} \frac{\ru_{\k+1}}{\k!} \xi^\k\right)\, e^{-\xi/\t} \d\xi = \ru_0 + \sum\limits_{\k=0}^{\infty} \frac{\k! \rv_k \t^{\k+1}}{(1+\t)\ldots(1+\k\t)},
 \end{equation}
 where $\rv_\k$ are computed using Stirling algorithm:
 \begin{equation}
  \label{vn_stirling}
  \rv_\k = \frac{1}{\k!} \sum\limits_{p=1}^{\k+1}| \s(\k,p-1)|\ru_p
 \end{equation}
 and $\s(\k,p)$ are the first order Stirling coefficients.
 Both algorithms, the \ac{BPL} and the one using \ac{FS} have been compared in \cite{ahmad_comp_bpl_sfg_2015} on a series of \ac{ODE}s problems. One of the advantages of \ac{FS} over the \ac{BPL} is numerical stability as it does not involves poles as the Padé approximants do. Another one is the possibility of applying the process to coefficients $\ru_\k(\bx)$ that are space functions. Recall that the poles of the effective sum using \ac{FS} are the set of points $\{-\frac{1}{\k},  \k =1,2,\ldots\}$ which lies in the negative part of the real axis, thus the sum does not present any singularity as the time is a positive variable.

 \subsection{\ac{FS} Algorithm}
In practice, we use a finite sum of the series presented in the right side of \cref{FS-BL} to approximate the Borel sum.
As $\k$ increases, Stirling coefficients increase exponentially which alter numerically the process. Instead, we use the following identity for Stirling coefficients:
\begin{equation}
 \label{stirling_id}
 |\s(\k+1,p)| = \k|\s(\k,p)| + |\s(\k,p-1)|, \quad \forall  1 \leqslant p \leqslant \k +1,
\end{equation}
thus, we can compute $\rv_\k$ using the algorithm in \cite{Thomann-90}:
\begin{equation}
 \label{w_k}
\w_{\k+1} \coloneqq  \frac{\k! \rv_k \t^{\k+1}}{(1+\t)\ldots(1+\k\t)} = \frac{1}{(1+\t)\ldots(1+\k\t)} \sum\limits_{p=0}^{\k}|\s(\k,p)|\w_{p+1}^{(1)},
\end{equation}
with $\w_{p+1}^{(1)} = \ru_p \t^{p+1}$, and the $\w_{\k+1}$ are computed in a recursive way such that:
\begin{equation}
\label{digram_thoman}
\left\lbrace
\begin{array}{llll}
 \w^{(1)}_{p+1}&=&\ru_p \t^{p+1}, &p\geq0 \\ \\
 \w^{(j+1)}_{p+1}&=&\cfrac{{(j-1)}\w_{p+1}^{(j)}+\frac{1}{\t} \w_{p}^{(j)}}{\frac{1}{\t}+{j}}, &p\geq j\geq1\\\\
 \w_{\k+1}&=&\w_{\k+1}^{(\k+1)}, & \k\geq0.
\end{array}\right.\end{equation}
In this case the Borel sum could be written:
\begin{equation}
 \label{Borel-FS}
 \mathcal{S} \hbu(\t) = \ru_0 + \sum\limits_{\k=1}^{\infty} \w_{\k}.
\end{equation}
In practice, we approximate the Borel sum by the following finite sum:
\begin{equation}
 \label{finite_FS_sum}
 \mathcal{I}^\m(\t) \coloneqq  \ru_0 + \sum\limits_{\k=1}^{\m} \w_{\k}
\end{equation}

\subsection{Continuation process-numerical flow}
From \cref{digram_thoman} we build the numerical flow of the \ac{FS} scheme.
Having a set of equally spaced discrete points $\t_\mn$ with a time step $\tau$, we denote by $\bu_\mn$- without $\bx$ for simplicity- the approximation of $\bu(\t_\mn,\cdot)$ using a given numerical flow. Starting from $\t_0 = 0$ with $\ru_0\equiv \bu(0,\cdot)$, we reach instant $\t_{\mn+1}$ to compute $\bu_{\mn+1}$ following these steps:

\begin{enumerate}
\item Consider the approximation $\bu_n$ at instant $\t_n$ as the rank $0$ of the time series expansion:$$\ru_0\gets \bu_\mn$$ .
\item Apply the variable change $\t - \t_{\mn}$.
\item Compute $\ru_k$  for $\k= 1,\ldots, \m$.
\item Apply following numerical flow
\begin{equation}
 \label{map-BPL}
 \bu_\mn\mapsto  \bu_{\mn+1}\coloneqq\Psi_{\tau}^{[\m]}\big(\bu_\mn\big) =\mathcal{I}^\m(\tau) ,
\end{equation}
producing an approximation of $\bu(\t_{\mn+1} \equiv \t_\mn +\tau,\cdot)$.
\item $\mn \gets \mn+1$  and repeat from step 1.
\end{enumerate}

The above steps are repeated until reaching $\mn$ such that $\t_{\mn} \geqslant \T{}$. We denote that the second step, consisting of the computation of the terms of series depends on the problem itself and on the spatial framework used for the space discretization technique. For this moment on, we use \cref{mix_form_FE} in the framework of \ac{FEM}. Next, we will assess the proposed numerical method in solving the Taylor Green Vortex problem.
\captionsetup[table]{name=Table,skip=1ex}

\section{Taylor Green vortex test}
\label{sec5}

In this section, we consider a case study to assess the computation of the series terms, comparing them to the known exact closed-form analytical solution. The two dimensional Taylor Green vortex is represented by the following solution of the velocity field in the geometrical domain $\gf = [0,1]^2 \subset \mathds{R}^2$:
\begin{equation}
 \label{NS_TG_velo}
 \bu(\t,\bx) = \left( -\sin(\pi\x) \cos(\pi\y)\exp(-\frac{2\pi^2\t}{\Re}) , \cos(\pi\x) \sin(\pi\y) \exp(-\frac{2\pi^2\t}{\Re}) \right)^\top.
\end{equation}
Thus, the pressure will be given by the following function:
\begin{equation}
 \label{NS_TG-pres}
 \p(\t,\bx) = \frac{1}{4} \big(\cos(2\pi\x) +  \cos(2\pi\y)\big)\cdot\exp\Big(-\frac{4\pi^2\t}{\Re}\Big).
\end{equation}
The exact form of the  velocity vector and pressure spatial modes are given below:
\begin{align}
 \label{NS_TG_TSE}
 \ru_{\k}(\bx)&= \left( -\sin(\pi\x) \cos(\pi\y)\frac{(-2\pi^2)^\k}{\k!}, \cos(\pi\x) \sin(\pi\y) \frac{(-2\pi^2)^\k}{\k!} \right)^\top, \\
 \rp_\k(\bx) & = \frac{1}{4} \big(\cos(2\pi\x) +  \cos(2\pi\y)\big)\cdot \frac{(-4\pi^2)^\k}{\k!},
\end{align}
which can be directly obtained from \eqref{NS_TG_velo} and \eqref{NS_TG-pres} by expanding the $\exp(\cdot)$ function on a series.
The Taylor-Hood element is used for the \ac{FEMF}. It approximates locally the velocity field by second-order Lagrange polynomials ($\Pk{2}$) and first-order ones ($\Pk{1}$) for the pressure.
After the triangulation of $\gf$, the mixed formulation is adopted to solve the system $(\mU_{\k+1},\mP_{\k})^\top$ within the FEniCS project \cite{BarattaEtal2023,ScroggsEtal2022,LoggEtal2012}.

\begin{figure}[!ht]
 \begin{center}
  \includegraphics[scale=0.24]{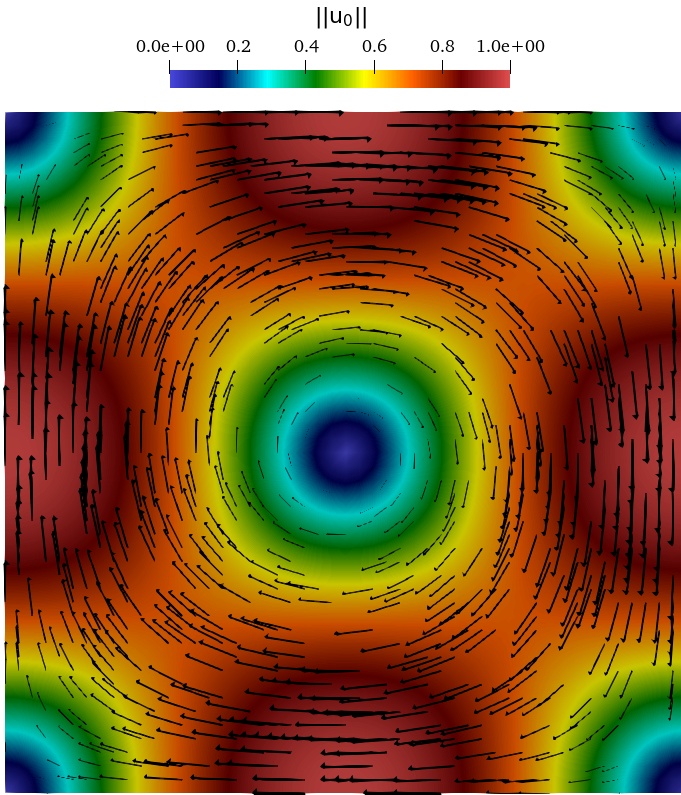}
 \end{center}
\caption{Streamlines of the magnitude of the initial condition vector field of the Taylor-Green vortex.}
\label{Fig1}
\end{figure}

The streamlines along with the magnitude of the initial condition $\ru_0 = (\u_0,\v_0)^\top$ of the velocity field vector, defined by:
\begin{equation}
 \|\ru_0\| = \sqrt{\u_0^2 + \v_0^2}
\end{equation}
are plotted in \cref{Fig1}.
Note that the exact solution of the $\k$\up{th} vector field has the same pattern as the $0$\up{th} term, but with a constant multiplication depending on $\k$. The mesh over the square $\gf = [0,1]\times[0,1]$ is constructed: the $\x$ and $\y$ axes are decomposed into a finite number of points that are uniformly spaced. Then, small squares are drawn in this square, to split every square into triangles forming the set of simplices of our triangulation.
The space of approximations $\V_{h,\k}^2$ and $\P_{h,\k}^1$ are constructed and the weak formulation is assembled using their basis elements to create the linear system \cref{mix_form_FE}.

\begin{figure}[!ht]
 \begin{center}
  \includegraphics[scale=0.24]{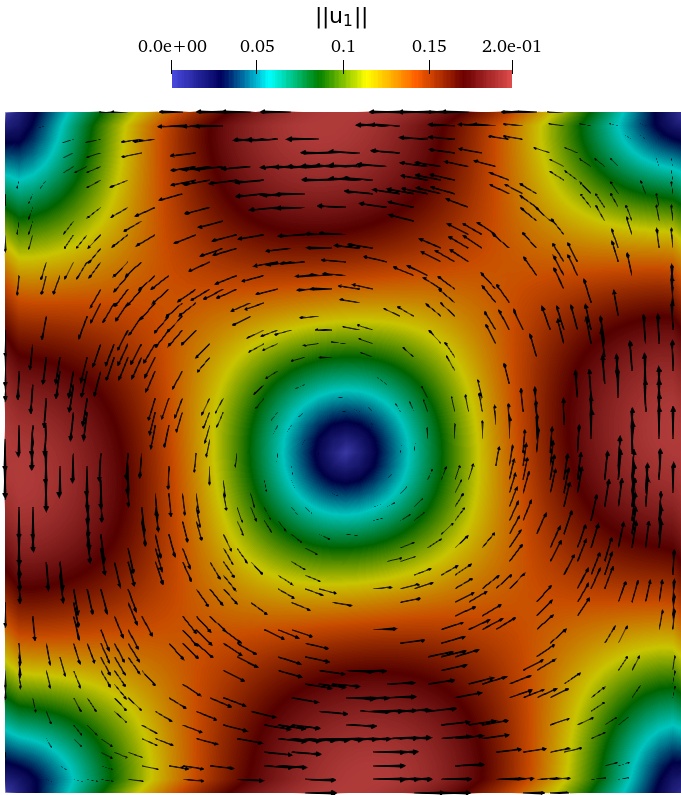}
  \includegraphics[scale=0.24]{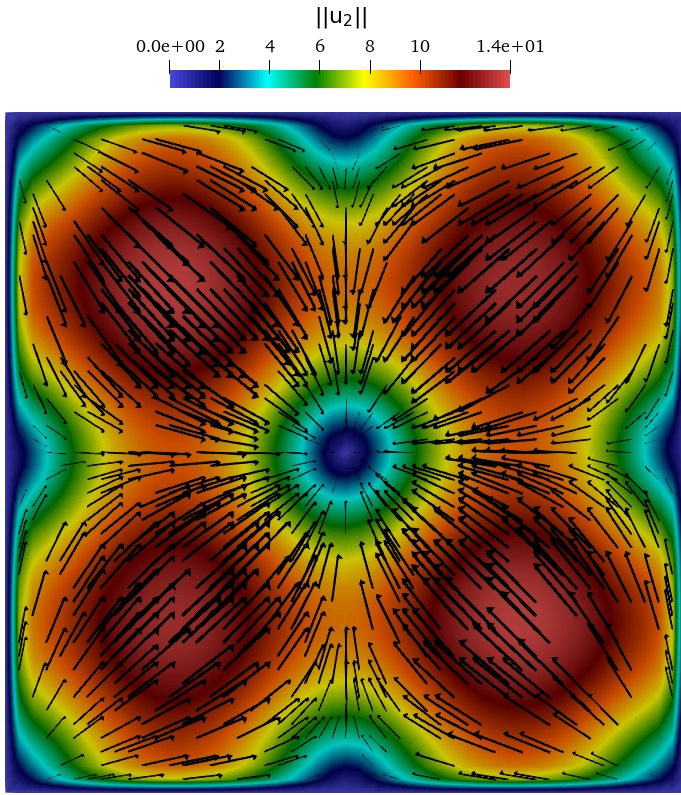}
 \end{center}
\caption{Streamlines of the magnitude of the Taylor-Green's first and second terms approximation in the \ac{FEMF} by means of the mixed formulation.}
\label{Fig2}
\end{figure}

After solving \cref{mix_form_FE} for $\k\in\{1,2\}$, \cref{Fig2} shows the streamlines of the approximation of the first and the second terms of the series in the \ac{FEMF} using the Taylor-Hood element. We can see that $\ru_1$ has a good approximation as the vector field kept almost the same streamline pattern, while the approximation of the vector field $\ru_2$ does not show this pattern, with presenting high values of its magnitude. This implies that the second term of the series is not valid and could not be considered in the computation for higher rank vector fields.

\begin{figure}[!ht]
 \begin{center}
  \includegraphics[scale=0.4]{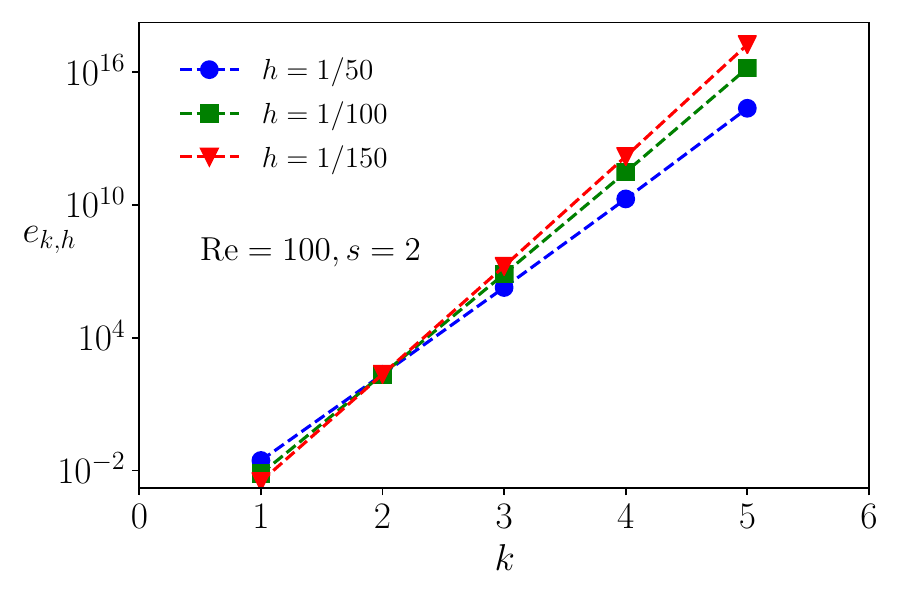}
  \includegraphics[scale=0.4]{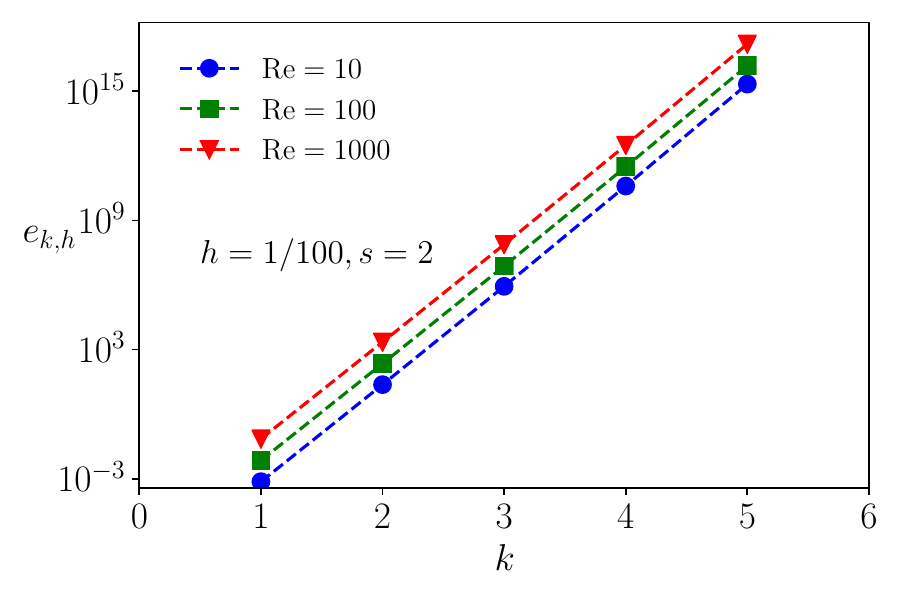}
 \end{center}
\caption{Error of the approximations of Taylor-Green terms for different mesh size $h$ (left panel) and for different Reynolds numbers (right panel),  for a fixed polynomial order $s$.}
\label{Fig3}
\end{figure}

The following error is evaluated to assess the quality of the approximation.
\begin{equation}
 e_{\k,h} \coloneqq \frac{ \displaystyle \int_\gf (\ru_{\k} -\ru_{\k,h})\cdot(\ru_{\k} -\ru_{\k,h}) \,\d\bx }{ \displaystyle \int_\gf \ru_\k \cdot \ru_\k \d\bx}. 
\end{equation}
This error is plotted in \cref{Fig3} for different simulation parameters.
First, the Reynolds number is set to $\Re=100$, and different mesh sizes were considered: $h\in\{1/50,1/100,1/150\} $. The left panel in \cref{Fig3} shows the errors $e_{\k,h}$ for the above features.
We see that the logarithmic error decreases when $h$ decreases for the first term of the approximation, but increases when $h$ decreases for $\k>2$. This aligns with the explanation and analysis given in \cite{deeb:stab-serie}. Namely, it was presented in that work that the error increases with $\k$ and when the size of the mesh decreases. This is related to the condition number of the mass matrix $\mA_{\k}$. The right panel of \cref{Fig3} shows the error for a fixed value of the mesh size and different values of Reynolds number, remarking that increasing its value makes the quality of approximations worse.

\begin{figure}[!ht]
 \begin{center}
\includegraphics[scale=0.46]{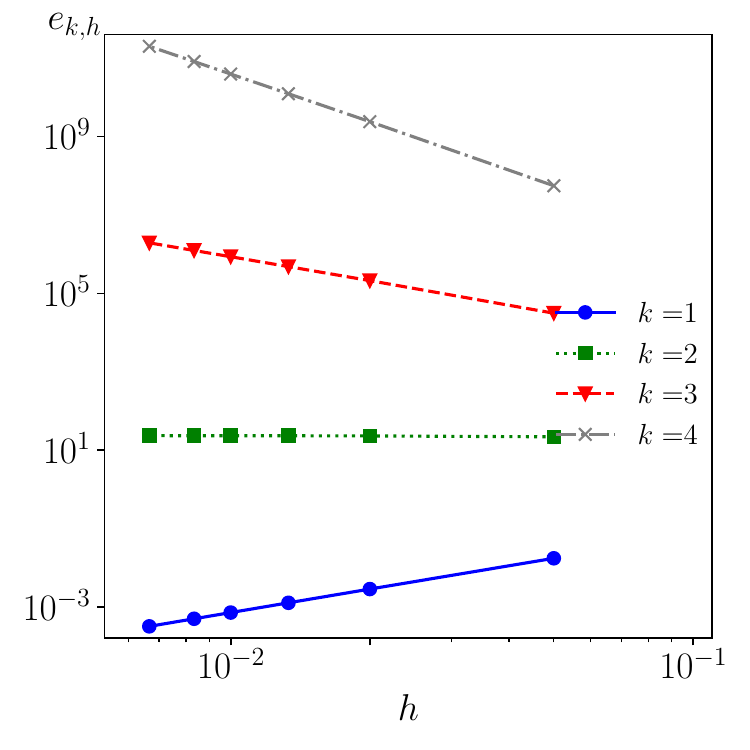}
\includegraphics[scale=0.46]{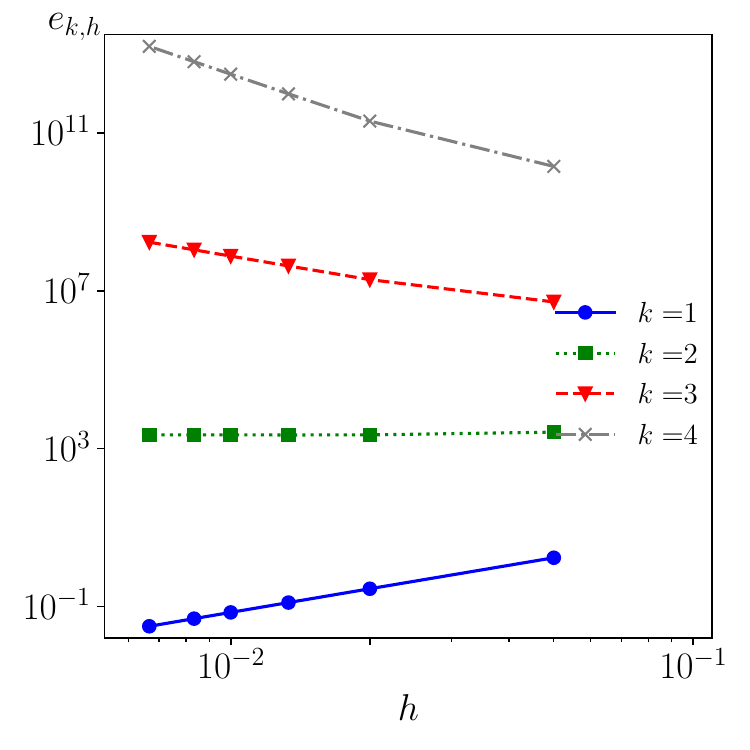}
 \end{center}
\caption{Evolution of the error for the first four terms of the Taylor Green time series solution with respect to the mesh size $h$ with two values of Reynolds number: $\Re = 1$ (left panel) and $\Re=1000$ (right panel).}
\label{Fig31}
\end{figure}

The evolution of the error $e_{\k,h}$  with respect to the mesh size $h$ and for a fixed value of Reynolds number is presented in \cref{Fig31}. We can see that the error $e_{1,h}\rightarrow 0$ when $h\rightarrow 0$ for the first term, while the error fails to converge to zero for the second and thus the third  and fourth also. We can see also that the $\Re$ number only affects the absolute value of the error $e_{\k,h}$ and not its rate of evolution with respect to $h$. This prove the necessity to stabilize the computation of the series terms to enable their use in the algorithm of \ac{FS}.

\begin{figure}[!ht]
 \begin{center}
  \includegraphics[scale=0.44]{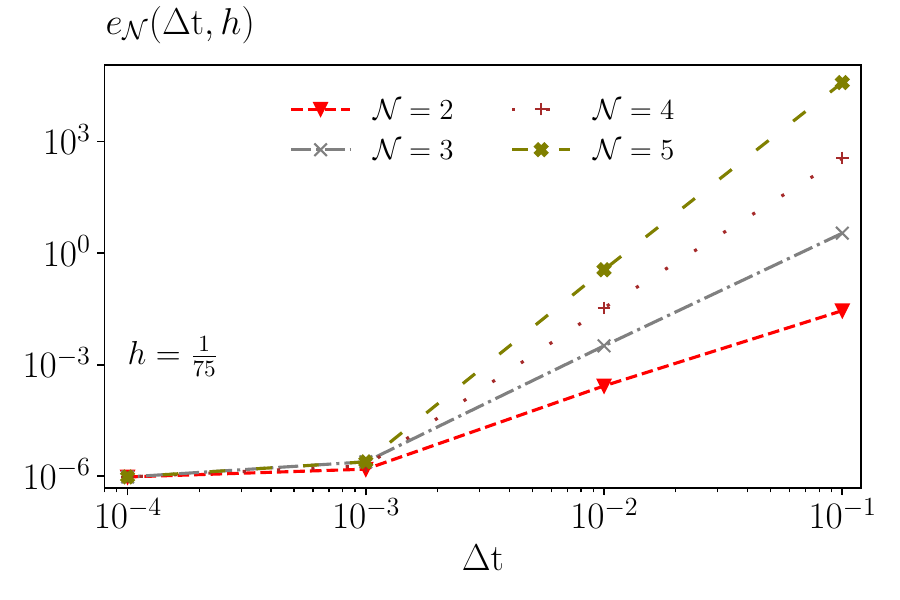}
  \includegraphics[scale=0.44]{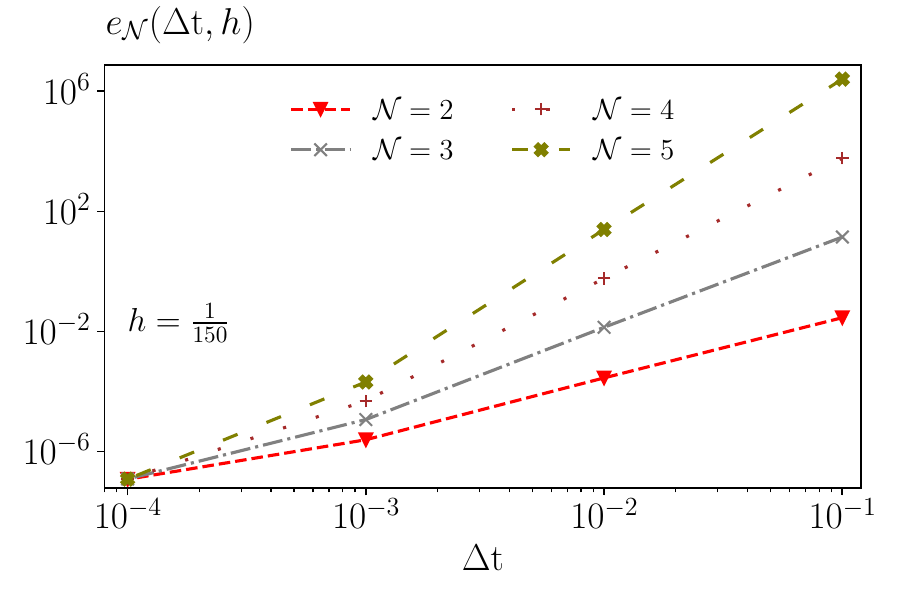}
 \end{center}
\caption{Plot of the error of the \ac{TSE}-\ac{FEM} when $h = 1/75$ (left panel) and $h=1/150$ (right panel) and $\Delta \t \rightarrow 0$ and for different truncation rank $\m = {2,3,4,5}$.}
\label{Fig32}
\end{figure}

We present now a graphical convergence analysis of the numerical method depending on the time step $\Delta \t$ and mesh size $h$. The Reynolds number is taken $\Re=100$. The plots in \cref{Fig32} show the evolution of the error w.r.t $\Delta \t$, for different truncation rank $\m$ and for two mesh size $h=1/75, 1/150$. As a first remark, the error increases with $\m$ which was already remarked in the plots of error in $\ru_\k$. But the error converges when $\Delta \t \rightarrow 0$ for the two mesh sizes $h$ with a bigger error when $h$ is small but both reaches the same value of error when $\Delta \t$ tend to 0: The convergence is faster for a smaller $h$. This means that the stability condition of the numerical scheme imposes a smaller time step for  smaller mesh size.

\begin{figure}[!ht]
 \begin{center}
  \includegraphics[scale=0.44]{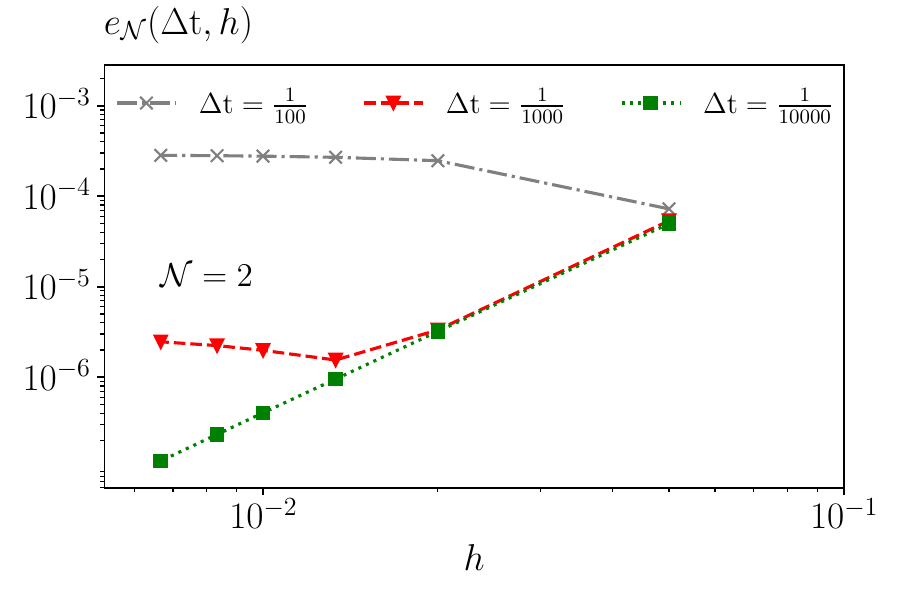}
  \includegraphics[scale=0.44]{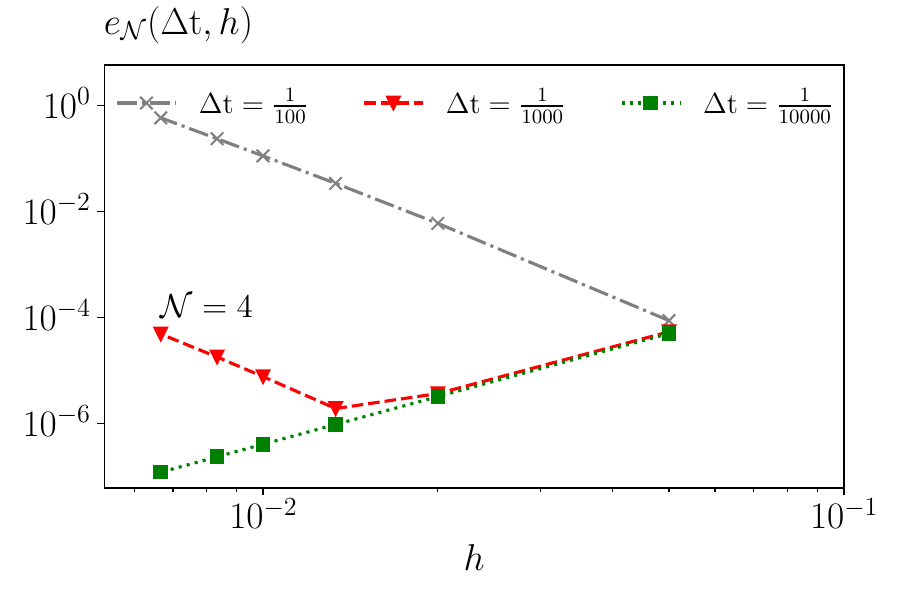}
 \end{center}
\caption{Plot of the error of the \ac{TSE}-\ac{FEM} when $\Delta \t$ is fixed and $\h \rightarrow 0$ and for two truncation rank: $\m=2$ in the left panel and $\m=4$ in the right panel.}
\label{Fig34}
\end{figure}

This can be seen in \cref{Fig34}, where the error does not converge at all w.r.t $h$ when $\Delta \t = 1/00$ and $\m=4$, but it converges first before it diverges w.r.t. $h$ when $\Delta \t= 1/1000$ and converges completely for $\Delta \t = 1/10000$. This behaviour is lowered for a small value of $\m=2$: We can see that the error for $\Delta \t=1/100$  diverges before it reach a stabilisation of order $\Delta \t^2$, and it diverges slower for $\Delta \t=1/1000$ and when $\m=2$ than the case of $\m=4$ for the same step size.

\begin{figure}[!ht]
 \begin{center}
  \includegraphics[scale=0.44]{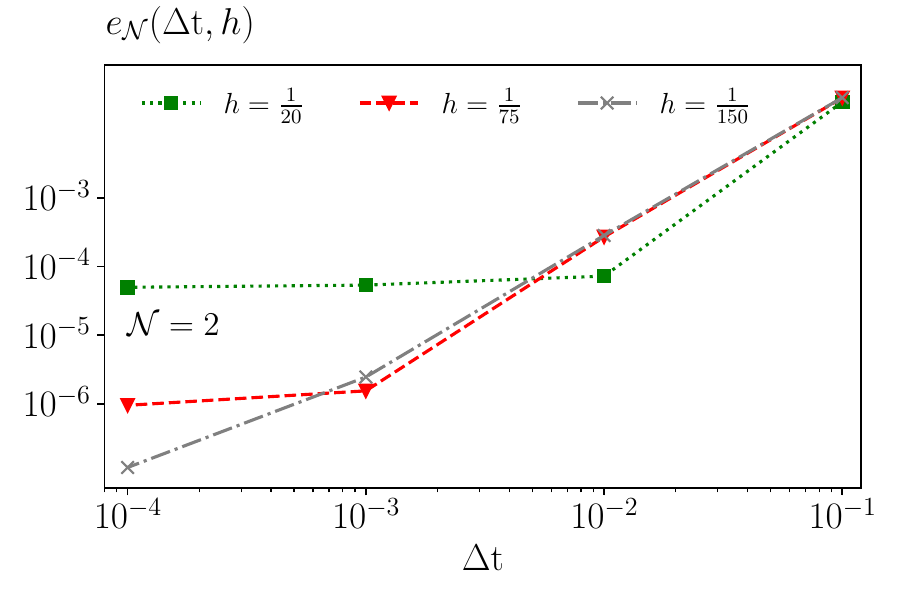}
  \includegraphics[scale=0.44]{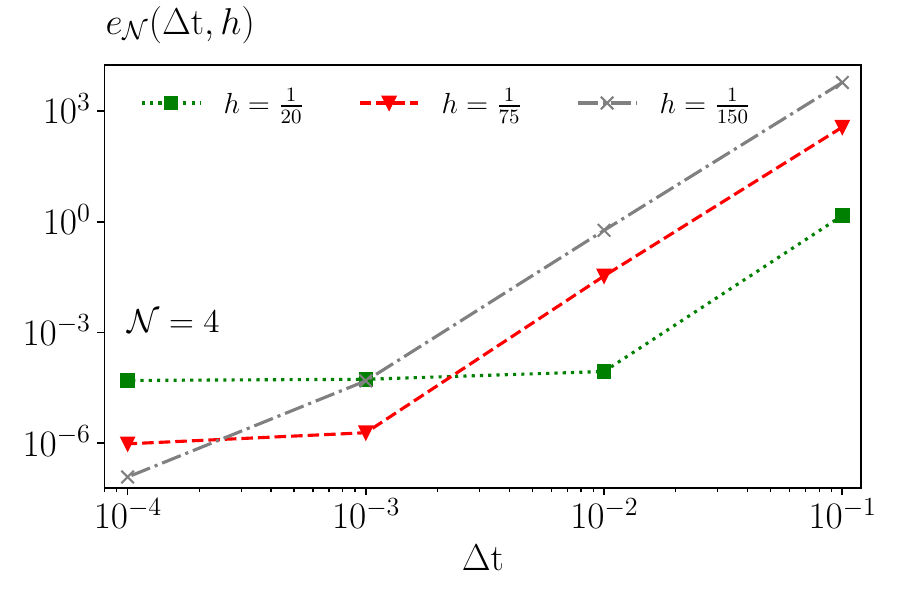}
 \end{center}
\caption{Plot of the convergence error w.r.t $\Delta \t$ for different mesh size $h=\{1/20,1/75,1/150\}$ and for two value of truncation rank: $\m=2$ in the left panel and $\m=4$ in the right panel.}
\label{Fig33}
\end{figure}
We end this section by showing in \cref{Fig33} the error of the approximation $e_\m(\Delta \t,h)$ w.r.t $\Delta \t$ and a set of mesh size values. We can see that the error converges for all cases before reaching a plate limited by the term related to the error in the mesh size. We can see also that when $\m$ increases, the error increases when $h$ decreases for a range of values of $\Delta \t$ ($\Delta \t> \epsilon$). This is related to the stability condition that need to be identified.

In the next section, we extend the stabilization technique, developed in \cite{deeb:stab-serie} for parabolic equations to the \ac{NS} equations in order to improve the computational validity of the terms vector fields for higher rank $\k$, thus the stability of the proposed numerical method.

\section{Stabilization technique}
\label{sec6}
In a recent paper \cite{deeb:stab-serie}, a stabilization technique was developed to improve the computation of the higher rank terms of the series in higher-order \ac{FEMF}. Regarding the problem of \ac{NS}, the idea is extended as follows: instead of solving System \ref{NS_Eqk}, an alternative one is solved by adding an artificial diffusion to the left-hand side of the first equation:
\begin{equation}
 \label{NS_Eqk_stab}
 \left\lbrace
\begin{array}{lcl}
(\k+1)\,\left(\ru_{\k+1} - \alpha_k \Delta \ru_{k+1}\right) + \nabla \rp_{\k} &=&\rF_{\k}  \\
 \nabla\cdot \ru_{\k+1} &=& 0
 \end{array}
 \right.,\quad
 \forall \bx \in \gf,\quad \forall \k \in \mathds{N},
\end{equation}
where $\alpha_k$ is a positive real coefficient, depending on the mesh size and on the term rank $\k$. It was shown in \cite{deeb:stab-serie} that if $\ru_{\k+1}$ is the minimum of a given functional as below:
\begin{equation}
 \label{minimum_original}
 (\ru_{\k+1},\rp_\k) \coloneqq \min\limits_{(\bv,\q)\, \in \V_\k\times \P} \J_\k(\bv,\q)
\end{equation}
then, adding the artificial diffusion creates the new minimization problem where the total variation of the unknown $\ru_{\k}$ is considered:
\begin{equation}
 \label{minimum_new}
 (\ru_{\k+1},\rp_\k )
 \coloneqq \min\limits_{(\bv,\q)\, \in \V_\k\times \P} \left(\J_\k(\bv,\q) + (\k+1)\frac{\alpha_{\k}}{2} \int_\gf \lvert \nabla\bv\rvert^2 \,\d\bx \right).
\end{equation}
This will filter out short wavelengths in $\ru_{\k+1}$ and minimize their amplification in the computation for $\ru_{\k+2}$.

In the next section, we show the relation between this stabilization and the \ac{NSa} model \cite{nsalpha-03}.

\subsection{Comparison with \ac{NSa} model}

If the coefficient of the stabilization is taken $\alpha_\k \coloneqq \alpha_0$, then the momentum balance in \cref{NS_Eq} is replaced by the following regularized equation:
\begin{equation}
 \label{NSfull-stab}
 \partial_t \ubu + \left(\bu\cdot\nabla \right)\bu  + \nabla \p = \frac{1}{\Re}\Delta \bu,
\end{equation}
with $\ubu \coloneqq (1- \alpha_0 \Delta)\bu$. This is the inverse of the regularization operator $(1- \alpha_0 \Delta)^{-1}$. Apply the latter on \cref{NSfull-stab} and use the  fact that it commutes with the time and space derivatives, we obtain the following system:
\begin{equation}
 \partial_t \bu + (\obu\cdot \nabla) \obu + \nabla \op = \frac{1}{\Re}\Delta \obu - \nabla\cdot \sT,
\end{equation}
where $\obu \coloneqq (1- \alpha_0 \Delta)^{-1}\bu$ is the filtered operation and $\sT \coloneqq \overline{\bu\otimes \bu} - \obu\otimes \obu$ is the so-called sub-grid scale tensor or Reynolds stress tensor.

The model we obtain in \cref{NSfull-stab} is different than the one defined by Leray \cite{leray-34} and given by $\obu \coloneqq \phi_\varepsilon\ast \bu$, where $\ast$ is the convolution operation. Guermond \emph{et al.} \cite{nsalpha-03} interprets the Leray regularization as the \ac{NSa} model
\begin{equation}
\label{NSa-model}
 \partial_t \bu + \left(\obu\cdot\nabla \right)\bu  + \nabla \p = \frac{1}{\Re}\Delta \bu \quad \text{with}\quad  \nabla\cdot \bu = 0 .
\end{equation}
Here, either the Helmholtz filter that is given by $\obu \coloneqq(1- \varepsilon^2 \Delta)^{-1}\bu $ or the Leray regularization lead to the above model.

\subsection{Determination of the artificial diffusion coefficient}

After constructing the \ac{FEMF}, the choice of the coefficient $\alpha_\k$ is considered to result from a recurrence formula as follows:
\begin{equation}
 \label{recurrence_alpha_k}
 \alpha_\k = \C{h}\times \alpha_{\k-1},
 \end{equation}
where $\alpha_0$ is taken as the $\arg\min$ of the condition number of the following matrix:
\begin{equation}
\label{argmin_alhpa0}
\alpha_0 \coloneqq \arg\min\limits_{\alpha>0} \kappa\left(\mA_{1} + \alpha\mK\right),
\end{equation}
with $\mK$ is the stiffness matrix defined  as:
\begin{equation}
\label{stiff_matrix}
 (\mK)_{\ell,\ell^{\prime}}   \coloneqq \int_{\gf} \nabla\bev_{\ell} : \nabla \bev_{\ell^{\prime}} \d\bx,
\end{equation}
and $\C{h}$ is to be chosen. If $\C{h}\coloneqq 1$, we found the \ac{NSa} model. Another choice of $\C{h}$ could be considered: $\C{h} = c^{m}$, where $c,m \in \mathds{R}^{+*}$.
Note that $\kappa(\mA) \coloneqq {\|\mA\|}\cdot{\|\mA^{-1}\|}$ is the condition number of the matrix $\mA$. 

To show the motivation behind \cref{argmin_alhpa0}, consider  solving the following linear system:
\begin{equation}
 \label{linear_sys}
 \mA \X = \b.
\end{equation}
If the right hand side $\b$ is perturbed by $\delta \b$, which is the case when $\ru_0$ is approximated by $\ru_{0,h}$ to be used in the computation of $\ru_{1,h}$, the solution  of the new problem is affected by $\delta \X$ as follows:
\begin{equation}
 \label{new_linear_sys}
 \mA (\X+ \delta \X)= \b +\delta\b,
\end{equation}
where the norm of $\delta \X$ is traditionally bounded as:
\begin{equation}
 \frac{\|\delta \X\|}{\|\X\|} \leqslant \kappa(\mA)\frac{\|\delta \b\|}{\|\b\|}.
\end{equation}
Thus, the reason of choosing $\alpha_0$ to be the argmin of the condition number $\kappa(\mA_{1} + \alpha \mK)$ is to minimize the amplification of the error which is produced in $\ru_{0}$ due to the discretization. Note that the computational cost of computing the condition number is high. However it is an offline computation for what we compute it only once at the beginning of the computation. In addition, we can use an estimation of the condition number based on the 1-norm condition estimator of Hager \cite{hager-84}. It was shown in \cite[Section 5.2]{deeb:stab-serie} that the proposition of the above minimization problem contributes in satisfying the \ac{DMP} \cite{DMP_Farago} in the case of heat equation.


\subsection{Stabilization algorithm}

To summarize the technique of stabilization for computing the terms of the series generated by \cref{NS_Eqk}, the following algorithm presents the steps to follow in the case of a given \ac{FEMF} adopted for \ac{NS} equations.

\begin{algorithm}[htp]
 \caption{Stabilization technique for solving problem \cref{NS_Eqk}.\label{ALgo1}}
 \begin{algorithmic}
  \Require $\ru_0$, $\Th$, $\V_{\k,h}^s$, $\P_{h}^r$
  \State $ \mU_{0,h} \gets$ project $\ru_0$ to $\V_{\k,h}^s$
  \State Assemble $\mA_1$, $\mK$
  \State Find $\alpha_0$ solution to \cref{argmin_alhpa0}
 \For{$\k\gets 0,...\m $}
 \State $\beta_{\k} \gets (\k+1)\alpha_{\k}$
 \State Solve the weak form of System \cref{NS_Eqk_stab} with boundary conditions:
 \begin{equation}
   \label{mix_form_FEstab}
 \left(
 \begin{array}{cl}
  (\mA_{k+1} + \beta_\k \mK) & \mB^\top \\ \mB&0
 \end{array}
 \right)
  \left(
 \begin{array}{cc}
  \mU_{\k+1,h} \\ \mP_{\k,h}
 \end{array}
 \right) =
   \left(
 \begin{array}{cc}
  \mF_{\k,h} \\ 0
 \end{array}
 \right) ,
 \end{equation}
\State $\alpha_{\k+1} \gets \C{h} \alpha_{\k}$
 \EndFor
 \end{algorithmic}
\end{algorithm}

In this paper, the above algorithm is programmed within the FEniCS project \cite{LoggEtal2012} under the Python programming language. It will be applied to compare the approximation of vector field terms $\ru_{\k,h}$, in the Taylor-Green problem, with stabilization to those without stabilization.
The first step consists of finding the recurrence formula for $\alpha_\k$, namely $\alpha_0$ and $\C{h}$.
Hereafter and for simplicity, we approximate $\alpha_0$ by $h^m$, though $m$ is to be determined. We plot in  the evolution of the condition number of the matrix $\mA_1+ h^m \mK$ for a range of values of $m \in [0,6]$.

\begin{figure}[!ht]
 \begin{center}
  \includegraphics[scale=0.65]{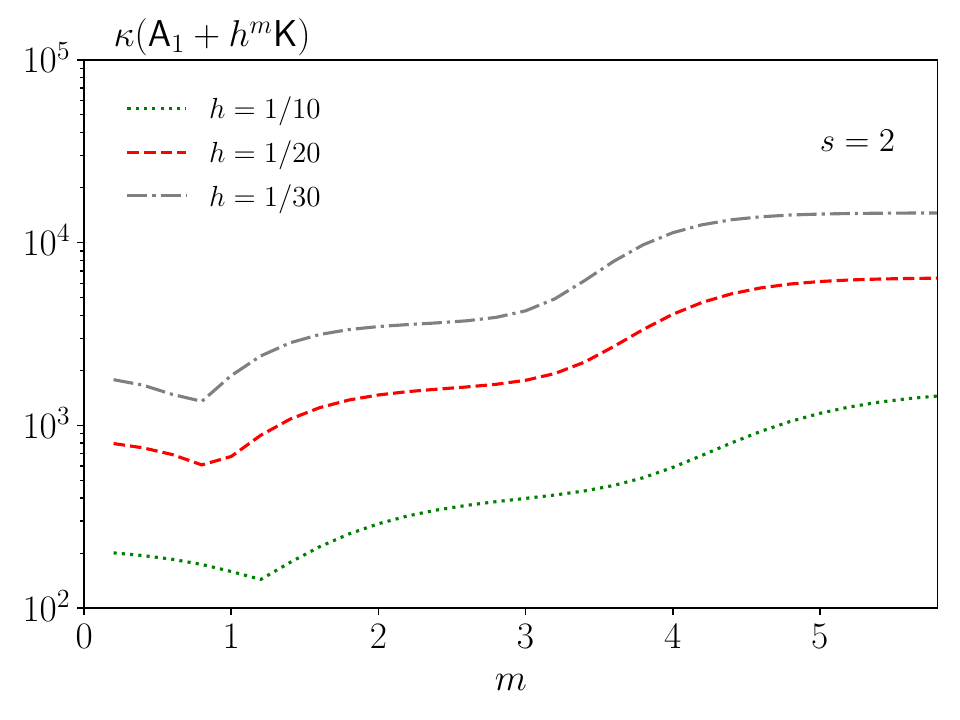}
 \end{center}
\caption{Condition number of the matrix $\mA_1 + h^m\mK$ for a range of values of $m\in ]0,6[$ and three different triangulations.}
\label{Fig5}
\end{figure}

\begin{figure}[!ht]
 \begin{center}
  \includegraphics[scale=0.4]{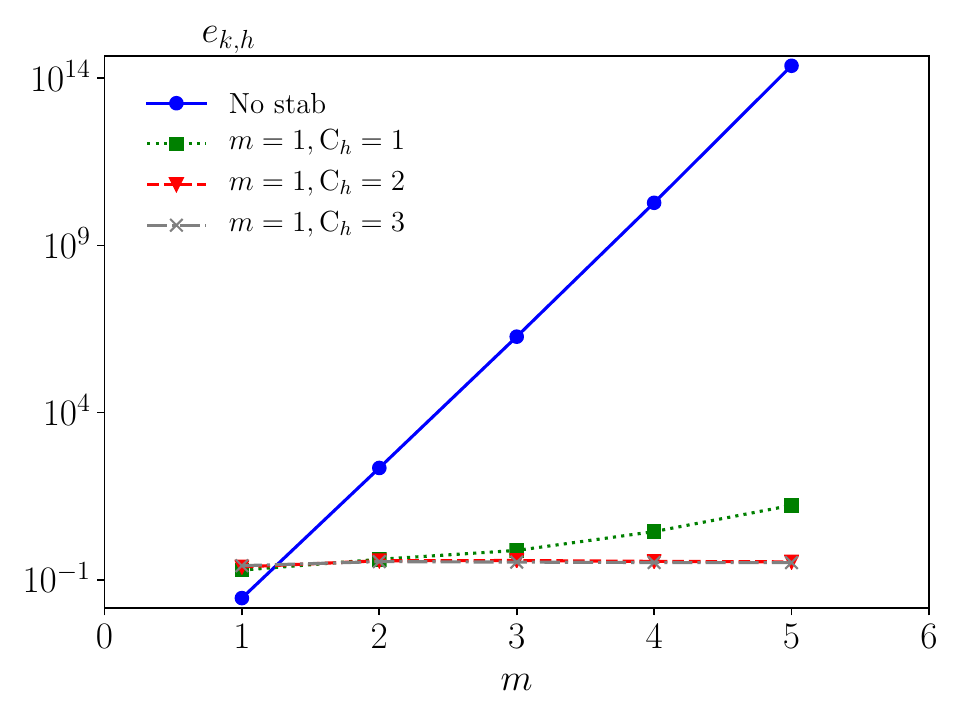}
    \includegraphics[scale=0.4]{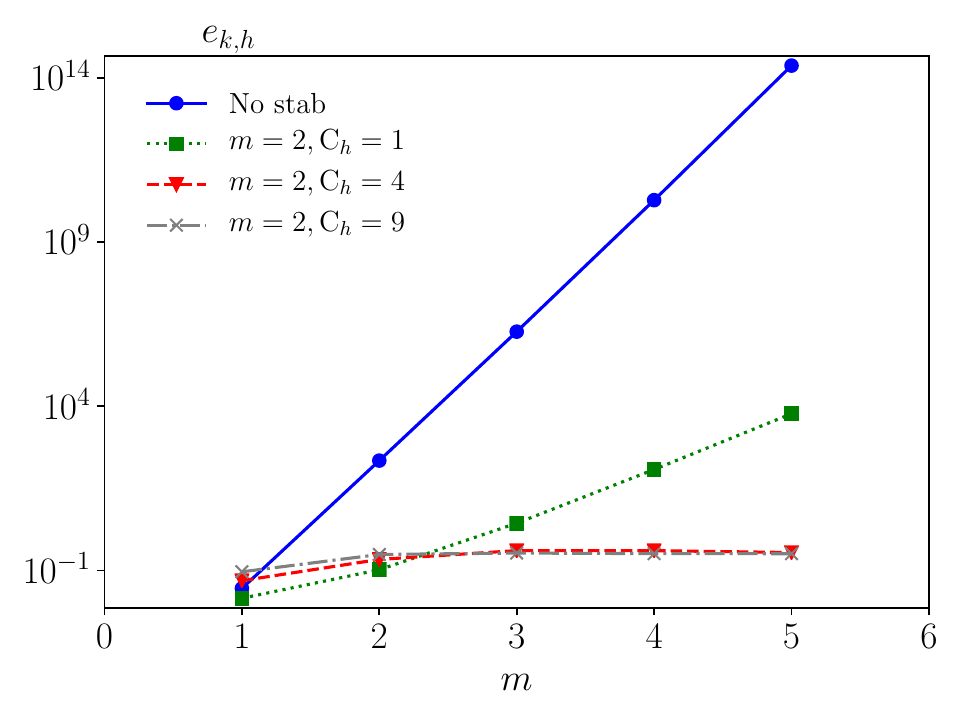}
 \end{center}
\caption{Relative error of the vector-field terms for different features of stabilization.}
\label{Fig6}
\end{figure}

\cref{Fig5} shows that there is a unique value on the considered interval that minimizes the condition number of the matrix. In addition, the minimum decreases when the size of mesh decreases too. If we take into consideration the error between the exact vector field term and its \ac{FEM} approximation, we tend to increase the value of $m$. To understand the choice of $\alpha_0$ and $\C{h}$, \cref{Fig6} presents the error of the first five vector-fields of the series with different features of stabilization.
We remark from these two panels, that the stabilization truly enhances the computation of terms of the series as the error was reduced drastically. This confirms the purpose of the stabilization. We remark also that for a fixed value of $\C{h}$, the greater the value of $\alpha_0 \coloneqq h^m$, the more the terms of the series are stable. We remark also that the errors present the same asymptotic when $\k \rightarrow \infty$ and when $\C{h}$ increases.

\begin{figure}[!ht]
 \begin{center}
  \includegraphics[scale=0.4]{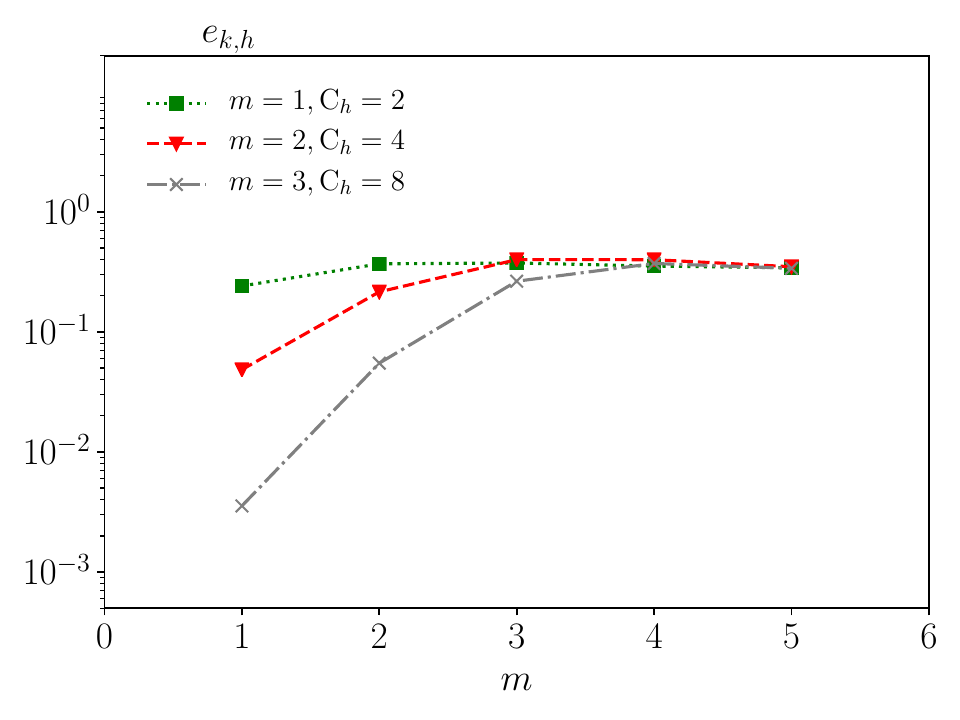}
    \includegraphics[scale=0.4]{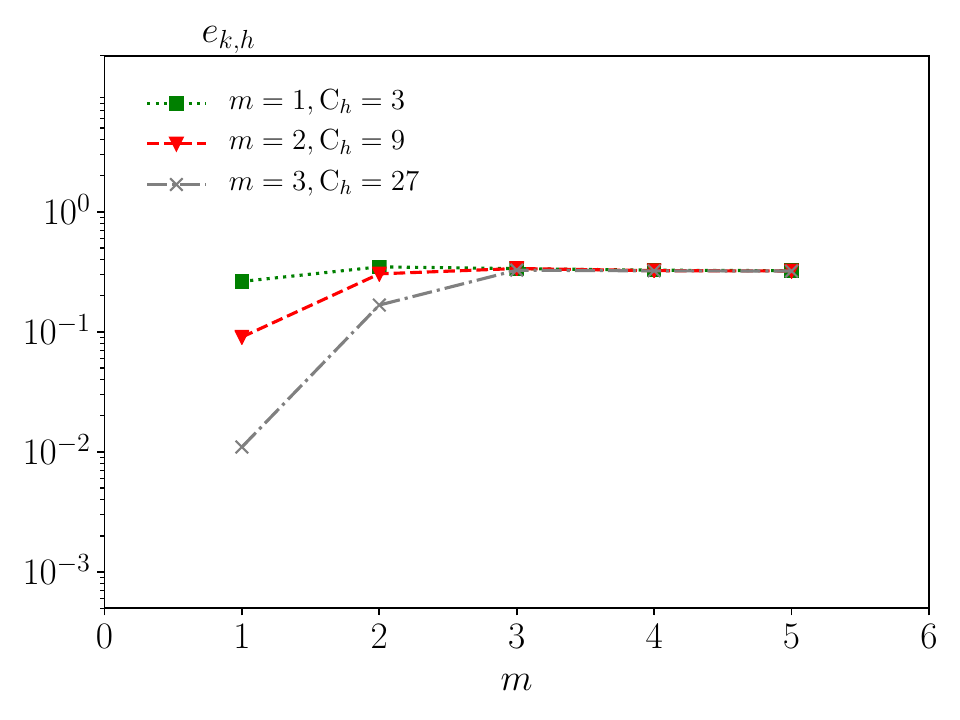}
 \end{center}
\caption{Relative error of the vector-field terms for different features of the stabilization.}
\label{Fig61}
\end{figure}

However, approximations present higher-errors in the first vector-field when high values of stabilization coefficients $\alpha_0\coloneqq h^m$ or $\C{h}$ are considered.
This is expected as the new vector-field terms are solutions to adapted stabilized problems with addition of $\alpha_\k \Delta \ru_{\k+1}$: The magnitude of the difference between the original problem and the stabilized one is of order $\alpha_\k$.

To check the impact of the choice $\C{h}$ on the stabilization technique, \cref{Fig6} shows that the computation of the terms of the series, when $\C{h} = 2^m$ were more stable than the case when $\C{h}=1$. The stability is not improved for higher values of $\C{h}$ as shown in \cref{Fig61}. Nevertheless, when $\C{h}$ increases for a fixed $\alpha_0$, the error in the first rank vector-field $\ru_{\k,h}$ increases. This aligns with an optimal choice to be taken, where the optimality will comprise between the stability of the computation of the terms and their errors.

\subsection{Convergence graphical study}

\begin{figure}[!ht]
 \begin{center}
\includegraphics[scale=0.35]{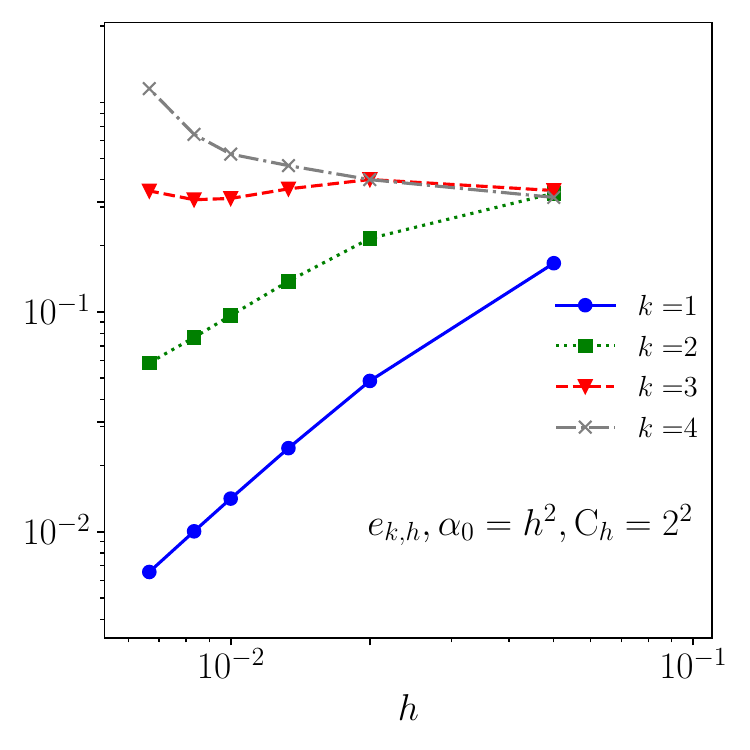}
\includegraphics[scale=0.35]{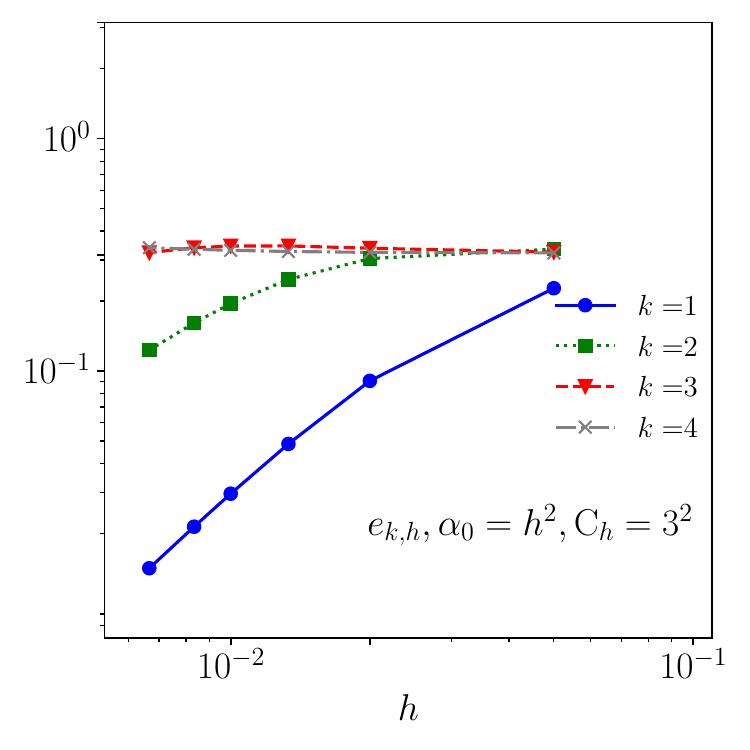}
\includegraphics[scale=0.35]{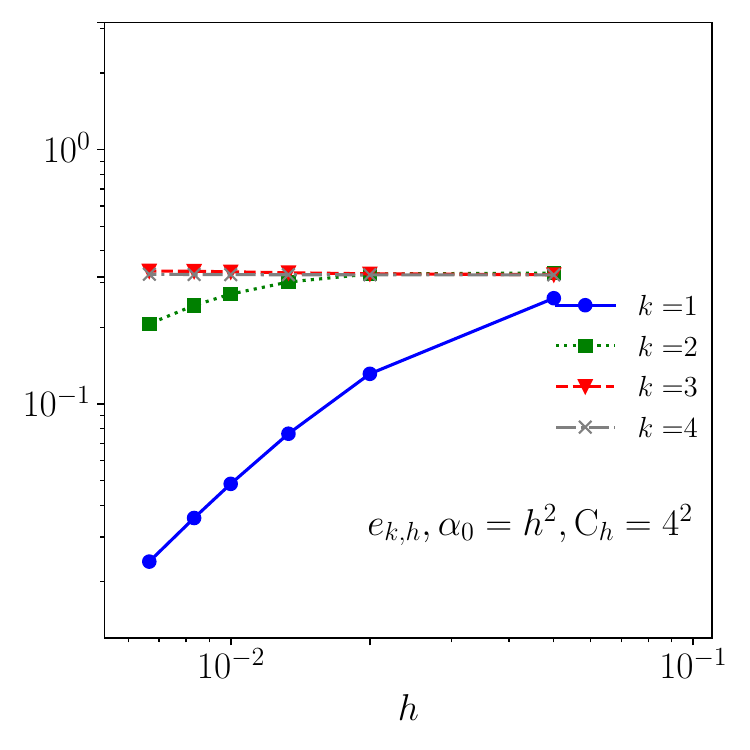}

\includegraphics[scale=0.35]{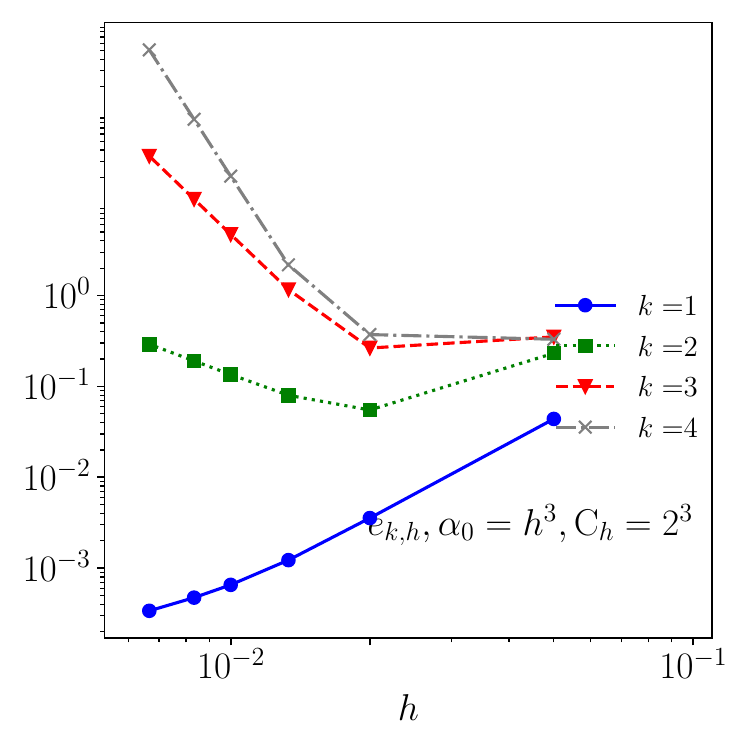}
\includegraphics[scale=0.35]{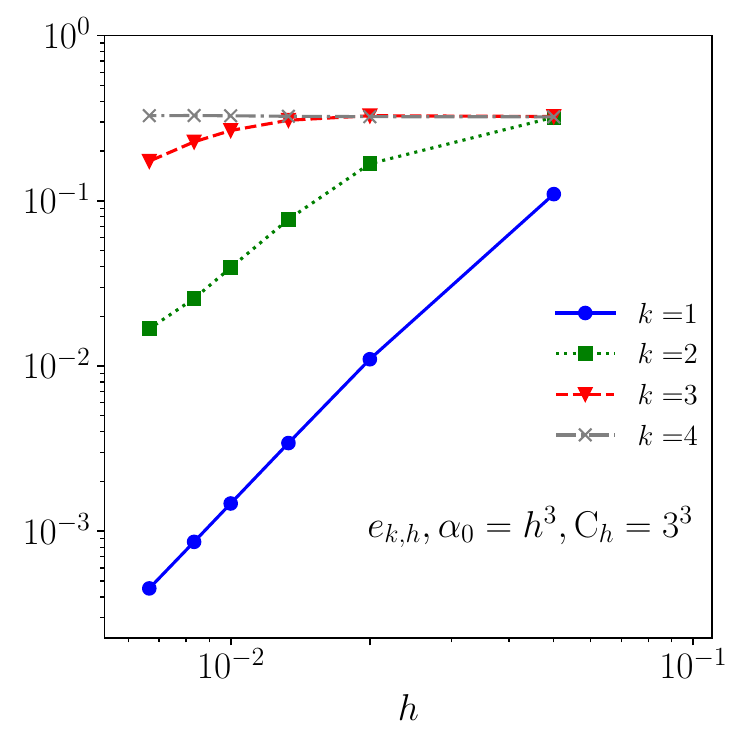}
\includegraphics[scale=0.35]{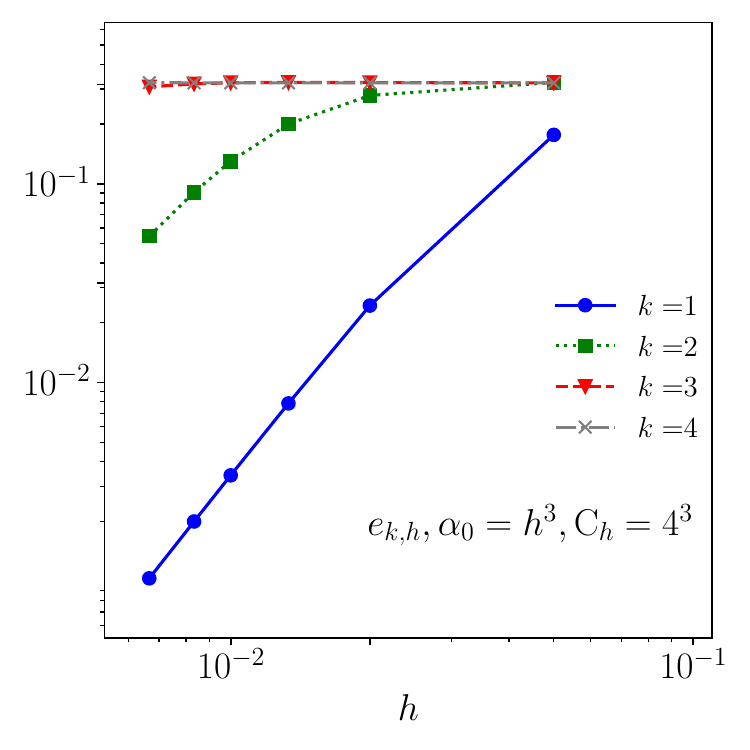}

\includegraphics[scale=0.35]{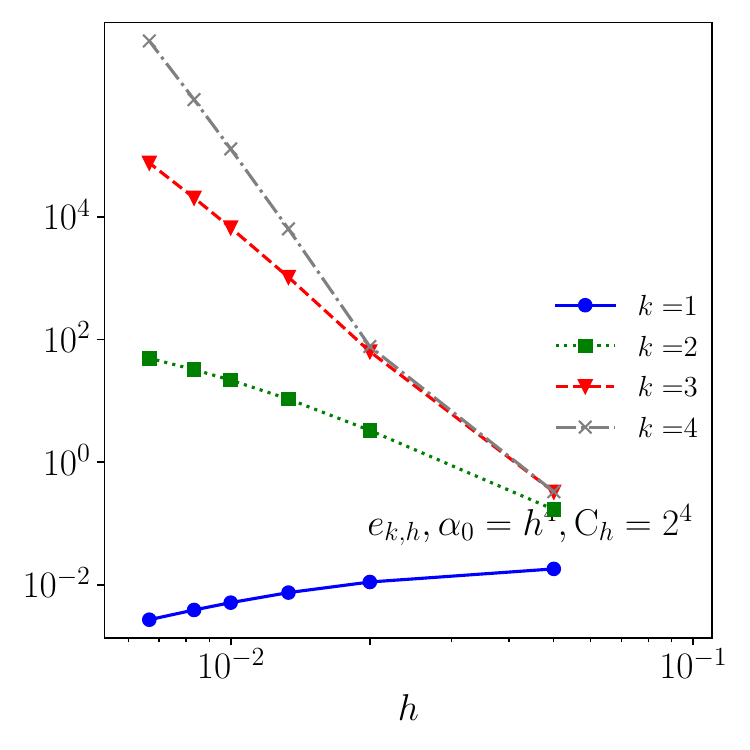}
\includegraphics[scale=0.35]{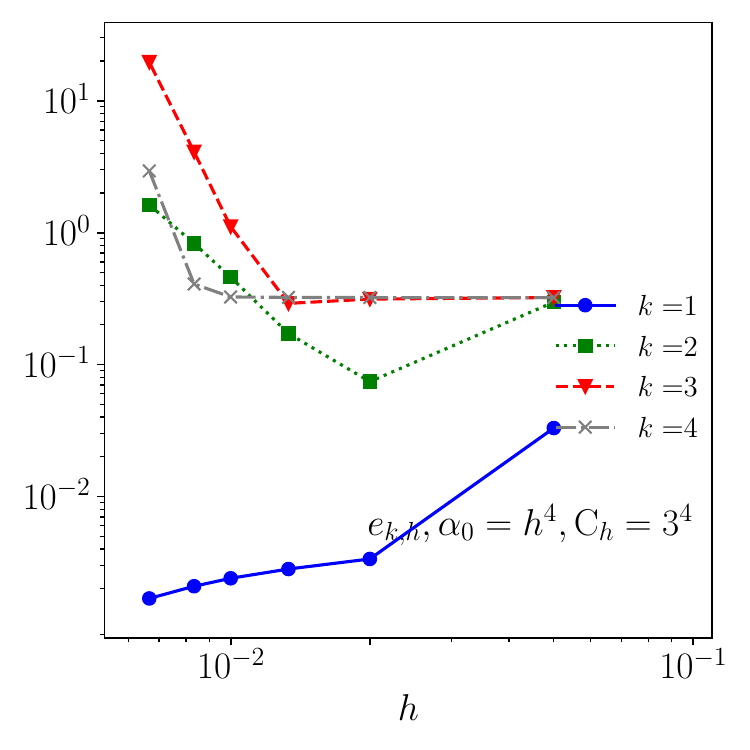}
\includegraphics[scale=0.35]{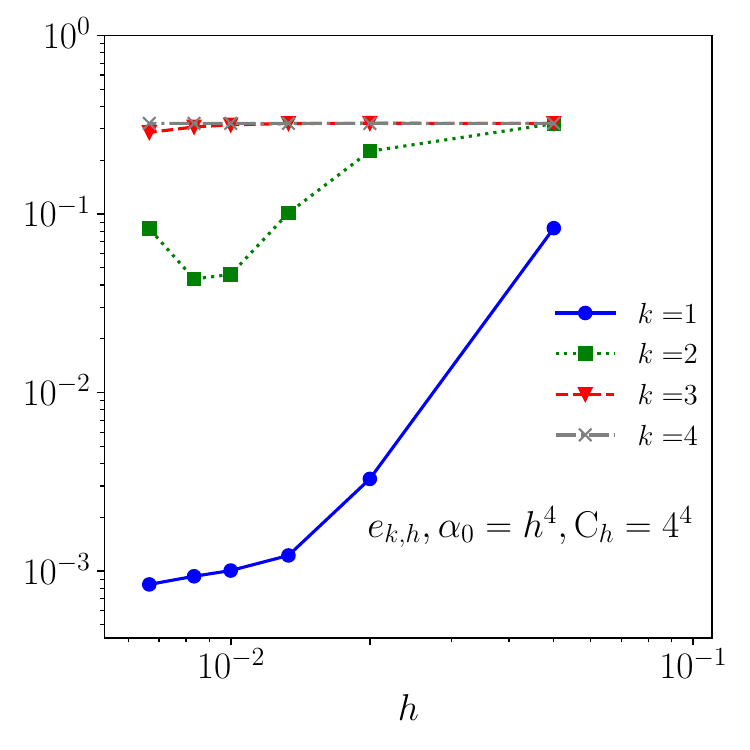}
 \end{center}
\caption{Convergence error of the first four terms of the series solution of Taylor Green for different stabilization features: $\alpha_0 = h^m$ for $m=2,3,4$ from top to bottom and $\C{h} = 2^m, 3^m,4^m$ from left to right. \label{Fig62}}
\end{figure}

To check the convergence of the stabilization process when $h$ decreases, we plot in \cref{Fig62} the evolution of the error  $e_{\k,h}$, for the first four terms of the series solution of the Taylor Green problem, with respect to $h$ for nine features of the stabilization technique: the cross of three values of $\alpha_0 = h^m$ for $m=2,3,4$  with three values of $\C{h}=2^m,3^m,4^m$. Every row in the plot represents the evolution of the error when the stabilization technique is applied using a fixed value of $\alpha_0$, while every column represents the evolution of the error for a fixed value of $\C{h}$. The first row of the plot presents the stabilization technique with a fixed value of $\alpha_0 = h^2$, the second row represent the cases for $\alpha_0=h^3$ and the third are for $\alpha_0=h^4$. Compared to \cref{Fig31}, all of the plots in this figure presents enhancing in the error when the stabilization is applied.
For all the studied cases, the error of the first terms converges to zero when $h\rightarrow 0$. This means that the stabilization does not deteriorate the quality of the first term that will be used to the computation of the following terms.

We can see that the plots in the upper triangular part of \cref{Fig62} presents the best results in convergence of the error and that for the first four terms. This aligns with the conclusion that the bigger the value of $\C{h}$ is, the better the stabilization is, and that is for the three cases of $\alpha_0$. This is compromised with a slight increase of the error when $\C{h}$ increases.  We can remark also that the plots in the lower triangular part presents the following evolution of the error for specific rank $\k$: The error decreases with $\h$ until reaching a minimum before increasing when $h$ decreases. The specific rank $\k$ increases when $\alpha_0$ decreases and when $\C{h}$ increases, which proves that the proposed stabilization technique, considering of choosing $\alpha_0 \sim h^m$ and $\alpha_\k = \C{h} \alpha_{\k-1}$ with $\C{h} \sim c^\k$, presents a convergence in the computation of terms of the series when $h\rightarrow 0$.

We end this convergence graphical study by showing in the error of the approximation obtained by the numerical method in the framework of \ac{TSE}-\ac{FEM}  with stabilization technique for three features: $\alpha_0 = h^m$ associated with $\C{h} = m^m$ for three values of $m = 2,3,4$. The error is calculated and plotted w.r.t to the time step $\Delta \t$ and mesh size $h$.

\begin{figure}[!ht]
 \begin{center}
  \includegraphics[scale=0.35]{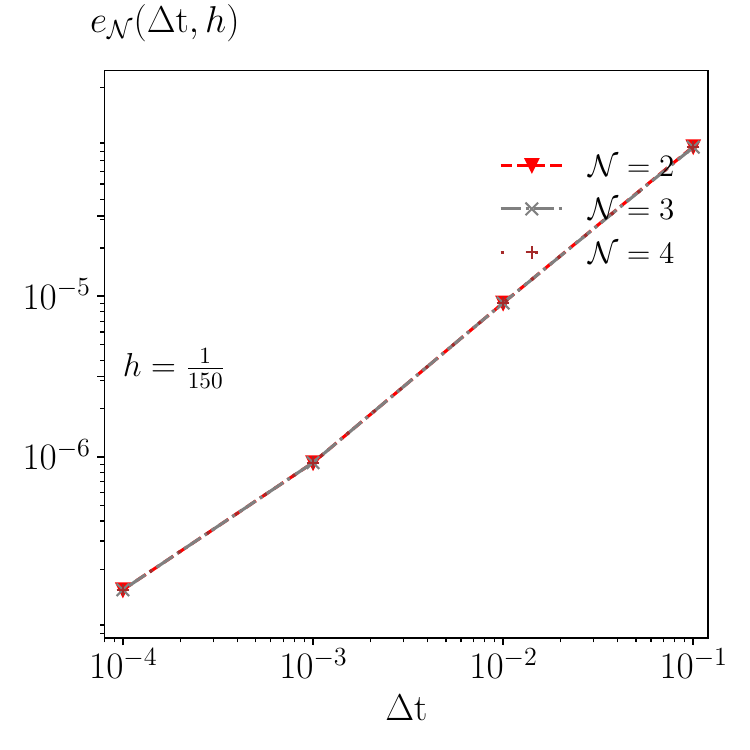}
  \includegraphics[scale=0.35]{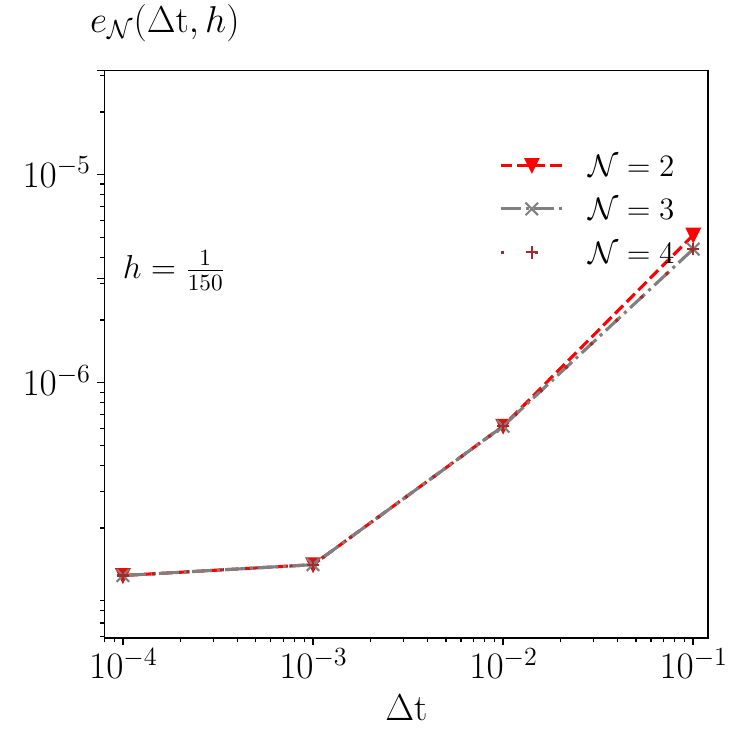}
  \includegraphics[scale=0.35]{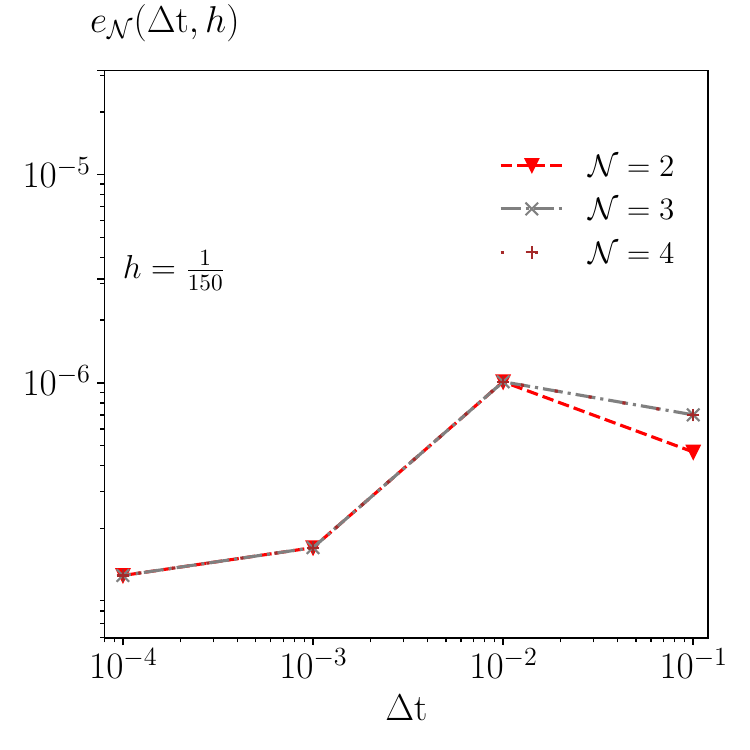}
 \end{center}
\caption{Plot of the error of the \ac{TSE}-\ac{FEM} with Stabilization when $h = 1/150$, $\Delta \t \rightarrow 0$ and for different truncation rank $\m = {2,3,4}$. The left panel represents the stabilization with $\alpha_0 = h^2, \C{h} = 2^2$, the middle panel is for  $\alpha_0 = h^3, \C{h} = 3^3$ and the right panel is for $\alpha_0 = h^4, \C{h} = 4^4$}
\label{Fig63}
\end{figure}

In \cref{Fig63}, the error w.r.t $\Delta \t$ presents a better stability with stabilization compared to plot in \cref{Fig32} regarding that the error does not increases when $\m$ increases. this aligns with the fact the the computation of terms does not explodes.

\begin{figure}[!ht]
 \begin{center}
  \includegraphics[scale=0.35]{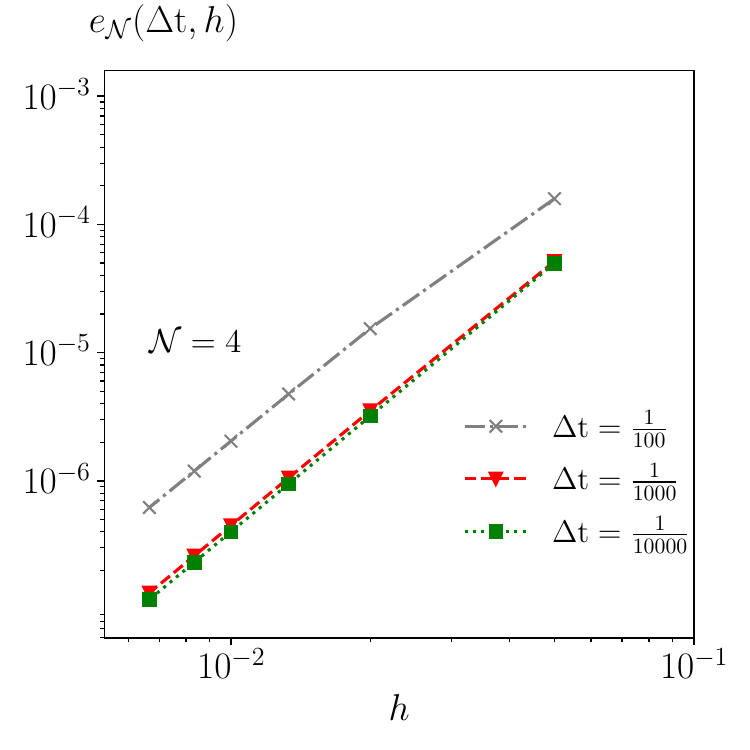}
  \includegraphics[scale=0.35]{FigStokesErr_hvsDt_stabTrue_2_N4_3_3_Nx150S.pdf}
  \includegraphics[scale=0.35]{FigStokesErr_hvsDt_stabTrue_2_N4_3_3_Nx150S.pdf}
 \end{center}
\caption{Plot of the error of the \ac{TSE}-\ac{FEM} with Stabilization when $\Delta \t$ is fixed and $\h \rightarrow 0$, for  $\m=4$ and for three features of stabilization: $\alpha_0 = h^2, \C{h} = 2^2$ (left panel),  $\alpha_0 = h^3, \C{h} = 3^3$ (middle panel) and $\alpha_0 = h^4, \C{h} = 4^4$ in the right panel.}
\label{Fig64}
\end{figure}

In \cref{Fig64}, the error decreases when $h$ decreases and that for all the values of the time step that were considered. This improve the stability condition of the scheme such that the one with $\Delta \t = 1/100$ has kept a stable simulation, when the one without stabilization had to have a time step that is 100 time smaller. This is promising in accelerating the simulation as the time step gets bigger.

\begin{figure}[!ht]
 \begin{center}
\includegraphics[scale=0.35]{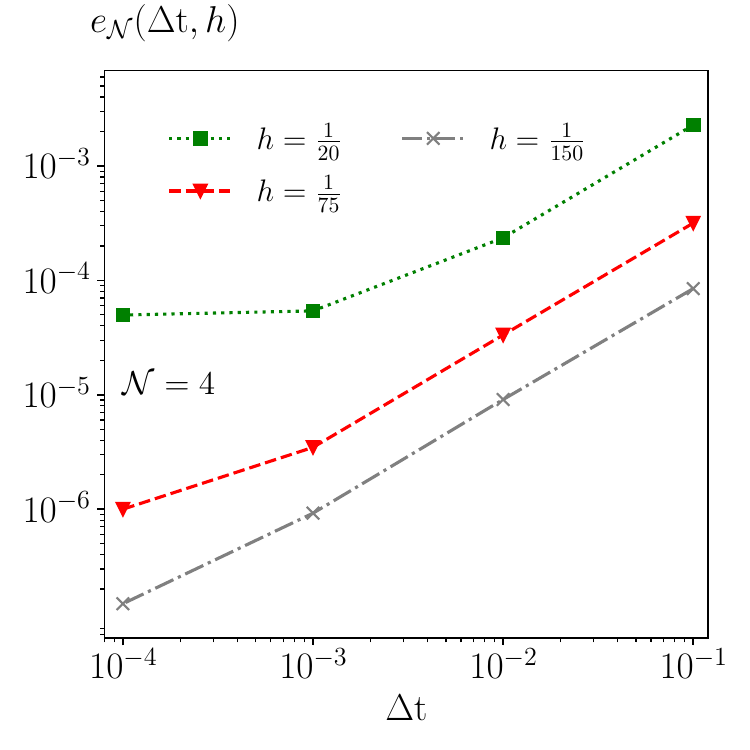}
\includegraphics[scale=0.35]{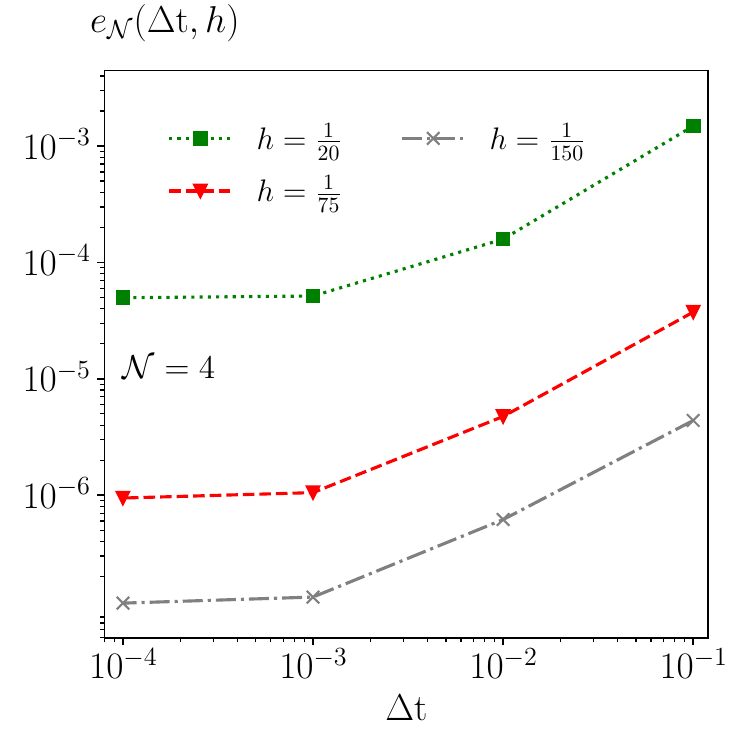}
\includegraphics[scale=0.35]{FigStokesErr_Dtvsh_stabTrue_2_N4_3_3_Nx150S}
\end{center}
\caption{Plot of the error of the \ac{TSE}-\ac{FEM} with Stabilization when $\h$ is fixed and $\Delta \t \rightarrow 0$, for  $\m=4$ and for three features of stabilization: $\alpha_0 = h^2, \C{h} = 2^2$ (left panel),  $\alpha_0 = h^3, \C{h} = 3^3$ (middle panel) and $\alpha_0 = h^4, \C{h} = 4^4$ in the right panel.\label{Fig65}}
\end{figure}

Although in \cref{Fig65}, when $h$ is fixed the error w.r.t $\Delta \t$ and with stabilization presents a similar decreasing  behaviour as of the errors without stabilization, as it was shown in \cref{Fig33}, but it is clear that the error is more consistent for smaller values of $h$ and for all the range values of $\Delta \t $ that were taken. This respond to the proposal work of the necessity of the stabilization technique  in the \ac{TSE}-\ac{FEM}.

\subsection{Continuation}
Having the initial condition $\ru_0$, the computation of the terms of the series is done up to a finite rank $\k=1,\ldots,\m$ ($\ru_0,\ldots,\ru_\m$), which means that the resummation of the truncated series verifies the equation over a finite time domain $\mathsf{T}_{0}^{1}\coloneqq [0,t_1]\subset [0,\T{}[$ up to a residual error $\tol$. This means that, for a fixed tolerance $\tol$ and a given rank $\m$, the associated numerical flow, denoted by $\Psi^{[\m]}$, approximates the solution for every $\tau \in \mathsf{T}_{0}^{1}$ as follows:
\begin{equation}
 \label{flow_resum}
 \bu_h(\tau,\cdot) \coloneqq \Psi_{\tau}^{[\m]}\left(\ru_0(\cdot)\right),
\end{equation}
while ensuring the following inequality:
\begin{equation}
 \label{residual_limit}
  \Res(\tau)(\ru_0,\rp_0)\leqslant \tol, \quad \forall \tau \in \mathsf{T}_{0}^{1}.
\end{equation}
Once the maximum value of $t_1$ is found, the solution will be approximated at this value $t_1$. This approximation is considered to be the initial condition for a new problem shifted by $t_1$ and the terms are computed again: this is the continuation process. The maximum $t_1$ could be obtained by a priori estimate study of the convergence of the numerical flow associated to the \ac{FS} scheme. This is done through a variable time step to follow the dynamics of the flow. This is not in the scope of this paper. We restrict our simulation in this work for a fixed time step during the simulation.

\section{Numerical tests}
\label{sec7}
In the following, a classical numerical test, the flow past a cylinder is to be performed for solving several settings for \ac{NS} equations using \ac{TSE} with stabilization techniques for the example of flow past a cylinder.
The triangulation of the geometry is built using the \textit{mshr} package in FEniCS software for different resolutions $N$.
Another option to generate meshes-- \textit{Gmsh} \cite{gmsh-09}-- a free open-source is available. It gives us the possibility to refine the mesh behind the cylinder. For the space approximation, we consider the mixed space $\V_{h}^2\times\P_h^1$ built on the $\Th$.

Regarding the stabilization of \ac{TSE} algorithm, we first construct the mesh size function, represented by the following discontinuous function $h(\bx)$ that belongs to the space of 0-discontinuous Galerkin element:
\begin{equation}
\label{function-mesh}
\h(\bx) = \sqrt{\int_{\K_j}1\, \d\bx}, \quad \forall \bx \in \K_j.
\end{equation}
Then the stabilization coefficient $\alpha_0$ is considered as follows:
\begin{equation}
\label{function-stab}
 \alpha_0(\bx) \coloneqq \h(\bx)^m.
\end{equation}
Note that this function could be amended by replacing it with a constant that is equal to the mesh size $h \equiv h_{\min} \coloneqq \min\limits_{\x \in \Th} \h(\bx) $, or $h \equiv h_{\max} \coloneqq  \max\limits_{\x \in \Th} \h(\bx)$, if the mesh is uniformly refined. The recurrence formula of $\alpha_\k$ follows several cases: $\alpha_\k = \C{h}\alpha_{\k-1}$, where $\alpha_0 = \h(\bx)^m$ is considered for different values of $m$ and $\C{h} \in\{\frac{1}{\k},1,2^m$\}. The computation of terms is done up to $\m = 5$, and the diagonal Padé approximant is considered. The \acf{GL} quadrature is done for $N_G=20$ points. The continuation time step $\tau$ is fixed during the simulation and the simulation is ran for different values of $\tau$.

\subsection{Laminar flow past a cylinder}
\begin{figure}[!ht]
 \begin{center}
  \includegraphics[scale=0.55]{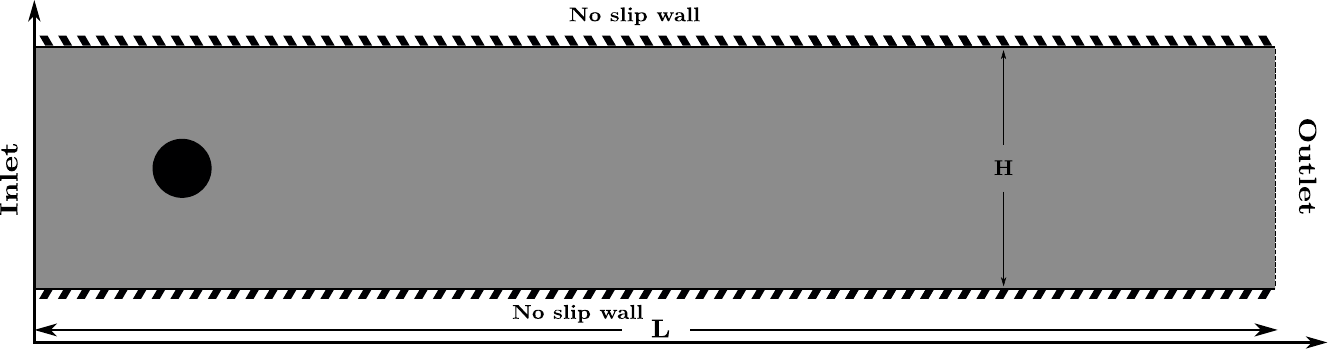}
 \end{center}
\caption{Schematic of a channel flow past cylinder.}
\label{Fig_channel_cylinder}
\end{figure}

We are considering the motion of a viscous, incompressible fluid in a rectangle of width $L=2.2$ and height $H=0.41$ around a cylindrical obstacle of center $(0.205,0.205)$ and radius $r_0=0.05$.
The initial condition is computed by solving the steady solution to the \ac{NS} equations, where the initial guess to the iterative solver will be the solution of the steady Stokes equation with the same Reynolds number.
Near the surface of the cylinder, the fluid velocity is zero due to the no-slip condition, as well as in the top and bottom boundaries of the rectangle. \cref{Fig_channel_cylinder} illustrates the geometry and the boundary conditions to be considered. On the left side of boundaries, we impose an inlet velocity condition:
\begin{equation}
 \bu_{\rm in} \coloneqq \bu(\t,\x=0,\y)= \left( 4\, \y\,(H-\y)/H^2,0 \right)^\top,  \x = 0,\, \forall \y\in [0,H],
\end{equation}
while on the right side a free stress tensor condition is imposed: $\sigma(\bx)\cdot \n = 0$.

For data post-processing, the drag $C_D$ and lift $C_L$ coefficients are evaluated:
\begin{align}
\label{drag-lift-coef}
 C_D &= \frac{F_D}{U_\infty^2\cdot r_0}, & C_L &= \frac{F_L}{U_\infty^2\cdot r_0}, & (F_L,F_D) &= \int_{\varGamma_c} \sigma(\bx)\cdot \n \,\d\bx
\end{align}
where $F_D$ and $F_L$ are the first and second component of the stress forces $\sigma(\bx)$ acting over the cylinder. First, we conduct a series of numerical simulations for a fixed Reynolds number $\Re = 400$ and different meshes generated by \textit{mshr}. We consider first the case without stabilisation and with a resolution $N = 100$. This means that the mesh size $h_{\min} \approx H/N $ and $h_{\max} \approx L/N$. If the algorithm is run with $\tau = 0.02$, the simulation was not able to reach $\t=0.1$ without exploding. Thus, we continue decreasing the time step to $\tau =0.01$ and the same blow-up happens so on until we prescribe $\tau = 0.0013$. We remind that we start with the initial condition that is the solution to the steady \ac{NS} problem, meaning that the terms of the series should be the null vector. We present in \cref{tab1} the values of the smallest time steps $\tau$ that should be prescribed to have a stable simulation along the time interval $[0,5]$ for the Reynolds number $\Re=400$ and for resolution $N$.

\begin{table}[!ht]
 \caption{Maximum values of the time step $\tau$ to ensure stable simulation reaching $\t=1$ for $\Re=400$ and without stabilization technique.\label{tab1}}
 \begin{center}
\begin{tabular}{|r|r|r|r|r|} \hline
$N$ & $35$ & $ 50$ & $75$ & $100$ \\\hline
$\tau$ & $0.002$ & $0.001$ & $0.0004$ & $0.0002$ \\ \hline
\end{tabular}
\end{center}
\end{table}
We can see in this table that the time step decreases as $N$ increases. This is because the approximation of the spatial modes in the time series expansion are not valid only for the first one. Thus, the \ac{TSE} algorithm without stabilization is equivalent to the first order explicit scheme.

\begin{figure}[!ht]
 \begin{center}
\includegraphics[scale=0.42]{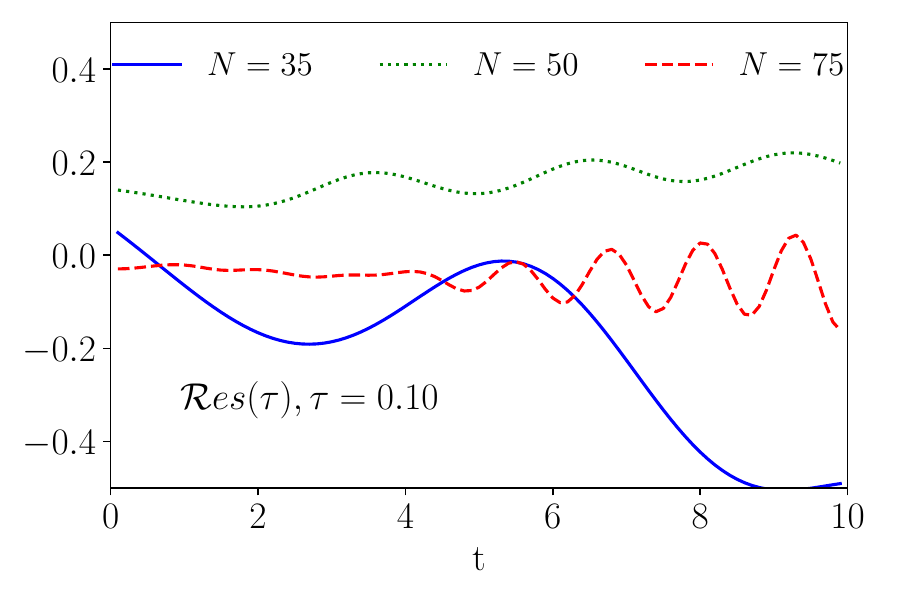}
\includegraphics[scale=0.42]{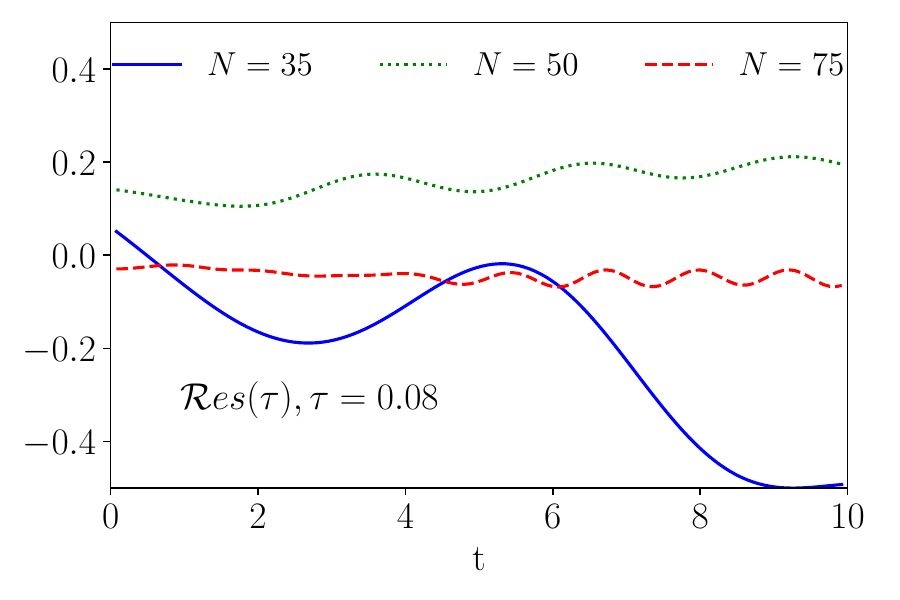}

\includegraphics[scale=0.42]{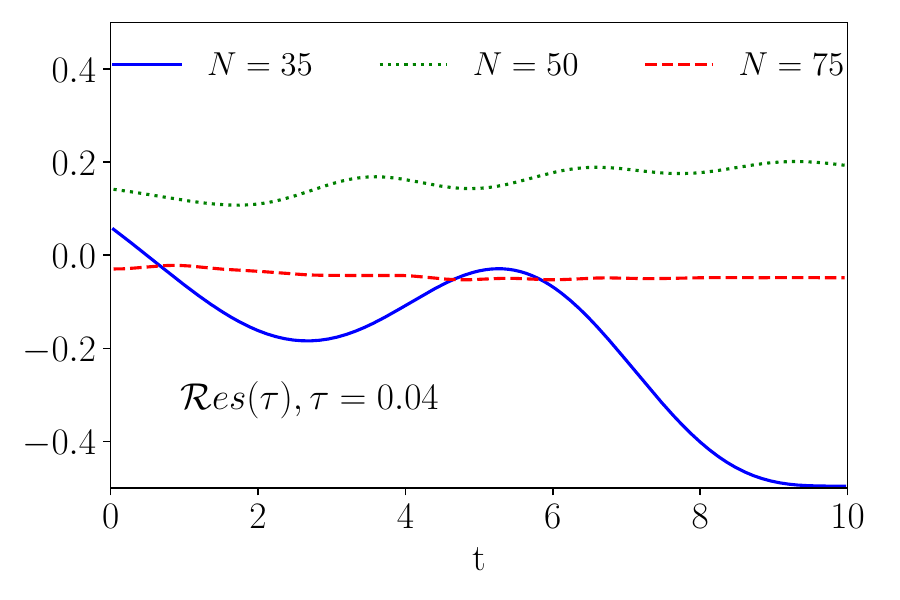}
\includegraphics[scale=0.42]{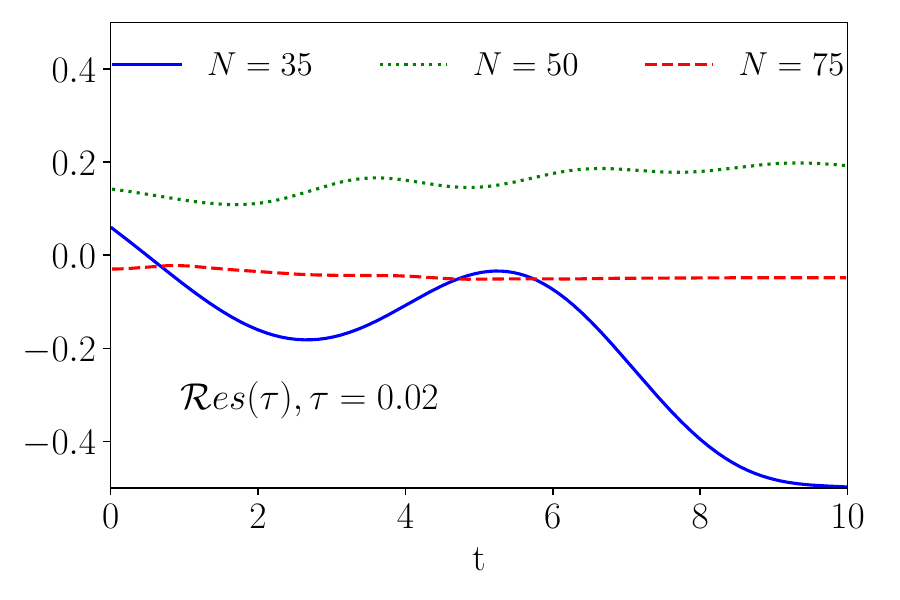}
 \end{center}
\caption{Evolution of the residual error during the simulation of the laminar flow past cylinder test using the \ac{FS} algorithm with stabilization technique: $m=1.8$, $\C{h} = 3^m$. \label{Fig7}}
\end{figure}

Nevertheless, if we consider the stabilization technique with $m=2$ and $\C{h} = 1$,  the simulation could reach $\t=10$ with $\tau = 0.04$ while keeping a completely stable computation. This means that the stabilization enhances the time step computation at least about $100$ times when $\Re = 400$ and the resolution is $N=75$.
Another simulation test is done with $\alpha_0 = h^m$ with $m=1.8$ and $\C{h} = 3^m$, with four time steps $\tau = 0.1,0.08,0.04,0.02$ and three mesh resolutions $N = {35},{50},{75}$. Results are plotted in \cref{Fig7}. We can see that the smaller the mesh size is or the time step is, the smaller the residual error absolute value is and its fluctuations.

In \cref{tab2}, we present the sum of the residual errors
\begin{equation}
 \label{residual_error}
 \Res(\tau) =\int_\gf \left[\left({(\bu_h)_t}+ (\bu_h\cdot \nabla)\bu_h\right)\cdot\bu_h +\frac{1}{\Re}\left(\nabla\bu_h:\nabla\bu_h\right) - \p_h \nabla\cdot \bu_h\right]\d\bx, \end{equation}
overall the time interval $[0,10]$ for different mesh resolutions $N$, different values of $m$ and $\C{h}$, and different time steps $\tau$.
\begin{table}[!ht]
 \caption{Values of the residual errors \cref{residual_error} for different space resolutions $N$. The double dashes -- means that the simulation explodes before reaching $\t=10$.\label{tab2}}
 \begin{center}

 \subfloat[$\C{h}=1$]{
\begin{tabular}{c|ccccccccc}
 \toprule
& \multicolumn{3}{c}{$m=1.8$} &\multicolumn{3}{c}{$m=2.0$} &\multicolumn{3}{c}{$m=2.2$} \\
& \multicolumn{3}{c}{$N$} &\multicolumn{3}{c}{$N$} &\multicolumn{3}{c}{$N$} \\
\cmidrule(lr){2-4}\cmidrule(lr){5-7}\cmidrule(lr){8-10}
  $\tau$ & $35$ & $50$ & $75$ & $35$ & $50$ & $75$ &$35$ & $50$ & $75$ \\ \midrule
  0.10 & 0.049 & -- & -- & 0.049 & -- & -- & 0.049& --& -- \\
 0.08 & 0.049 &0.037 & -- & 0.049 &0.037 & -- & 0.049& --& -- \\
 0.04 & 0.049 & 0.037 & 0.027 & 0.049 & 0.037&  0.027& 0.049 &0.037 &0.027 \\
 0.02 & 0.049 & 0.037 & 0.027 & 0.049 &0.037 & 0.027& 0.049 &0.037 & 0.027 \\
 \bottomrule
\end{tabular}}

 \subfloat[$\C{h}=2^{m}$]{
\begin{tabular}{c|ccccccccc}
 \toprule
& \multicolumn{3}{c}{$m=1.8$} &\multicolumn{3}{c}{$m=2.0$} &\multicolumn{3}{c}{$m=2.2$} \\
& \multicolumn{3}{c}{$N$} &\multicolumn{3}{c}{$N$} &\multicolumn{3}{c}{$N$} \\
\cmidrule(lr){2-4}\cmidrule(lr){5-7}\cmidrule(lr){8-10}
  $\tau$ & $35$ & $50$ & $75$ & $35$ & $50$ & $75$ &$35$ & $50$ & $75$ \\ \midrule
 0.10 & 0.048 & 0.037 & -- & 0.049 & -- & -- & --& --& -- \\
 0.08 & 0.049 &0.037 & -- & 0.049 &-- & -- & 0.049& --& -- \\
 0.04 & 0.049 & 0.037 & 0.027 & 0.049 & 0.037&  --& 0.049 & -- & -- \\
 0.02 & 0.049 & 0.037 & 0.027 & 0.049 &0.037 & 0.027& 0.049 &0.037 & 0.027 \\
 \bottomrule
\end{tabular}}

 \subfloat[$\C{h}=3^{m}$]{
\begin{tabular}{c|ccccccccc}
 \toprule
& \multicolumn{3}{c}{$m=1.8$} &\multicolumn{3}{c}{$m=2.0$} &\multicolumn{3}{c}{$m=2.2$} \\
& \multicolumn{3}{c}{$N$} &\multicolumn{3}{c}{$N$} &\multicolumn{3}{c}{$N$} \\
\cmidrule(lr){2-4}\cmidrule(lr){5-7}\cmidrule(lr){8-10}
  $\tau$ & $35$ & $50$ & $75$ & $35$ & $50$ & $75$ &$35$ & $50$ & $75$ \\ \midrule
 0.10 & 0.050 & 0.038 & 0.029 & 0.049 & 0.038 & -- & 0.049& --& -- \\
 0.08 & 0.050 &0.038 & 0.027 & 0.049 &0.038 & -- & 0.049& --& -- \\
0.04 & 0.050 & 0.038 & 0.027 & 0.050 & 0.038&  0.027& 0.049 & 0.038 & -- \\
 0.02 & 0.051 & 0.038 & 0.027 & 0.050 &0.038 & 0.027& 0.049 &0.038 & 0.027 \\
 \bottomrule
\end{tabular}}

\end{center}
\end{table}
We can see through \cref{tab2} that stabilization technique, with $\alpha_0 = h^2$ and $\C{h}=1,3^2$, improves the simulation stability of the \ac{NS} equations using the \ac{FS} algorithm by providing a time step that is $100$ times bigger in the case for $N=75$ without stabilization.

\begin{figure}[htpb]
\begin{center}
 \includegraphics[scale = 0.245]{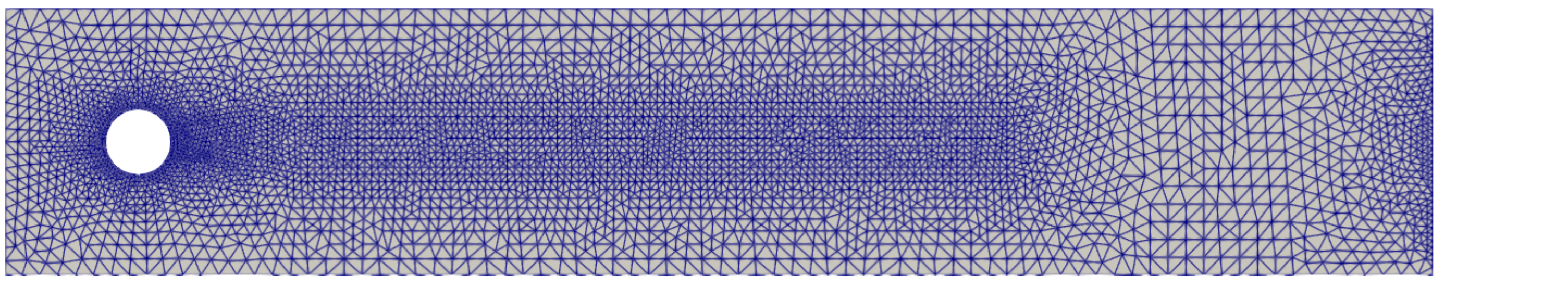}
 \caption{ Mesh discretization for the case when $\Re=4000$.\label{mesh-cyl_Re4000}}
\end{center}
\end{figure}

\begin{figure}[!ht]
 \begin{center}
\includegraphics[scale=0.245]{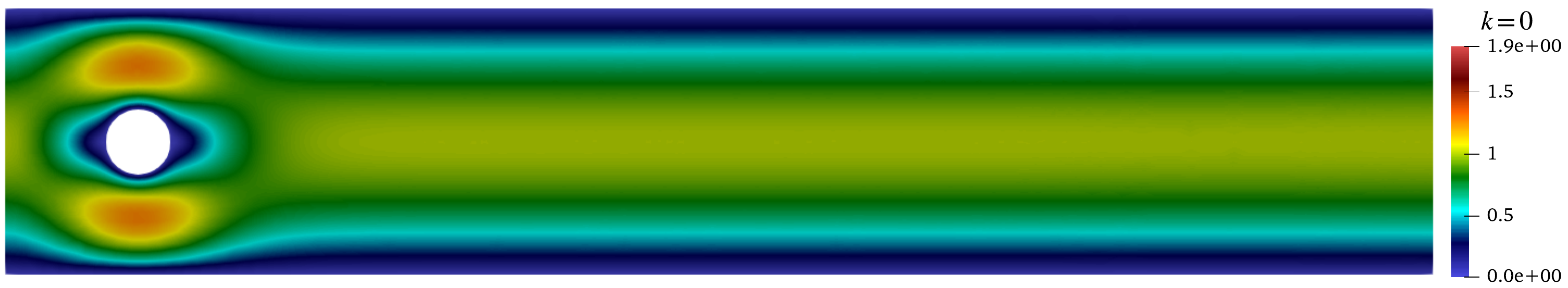}
\includegraphics[scale=0.245]{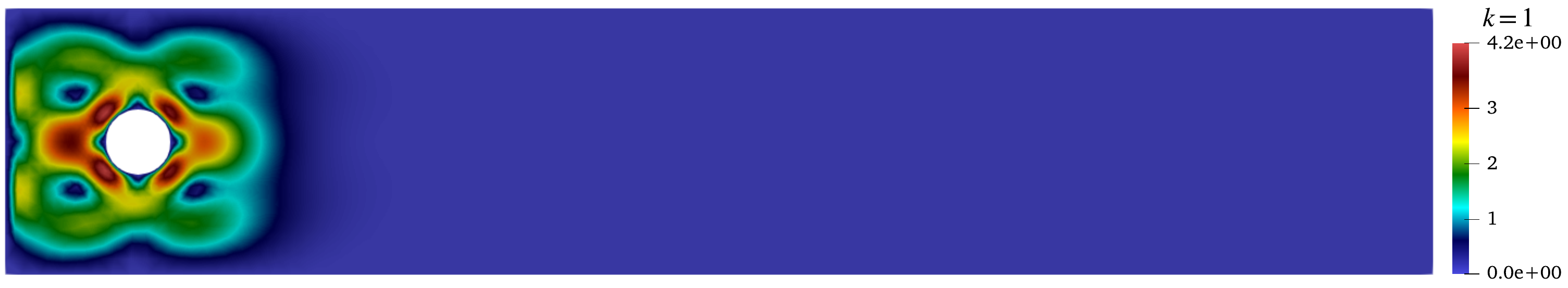}
\includegraphics[scale=0.245]{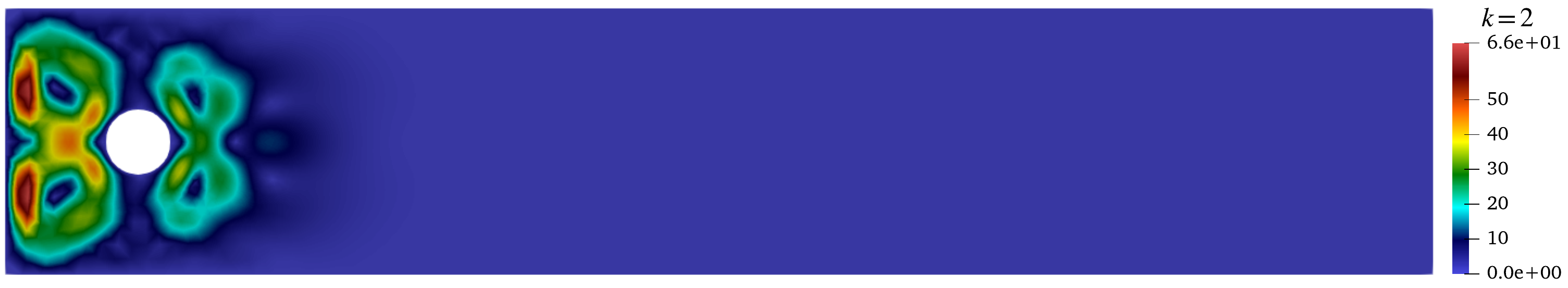}
\includegraphics[scale=0.245]{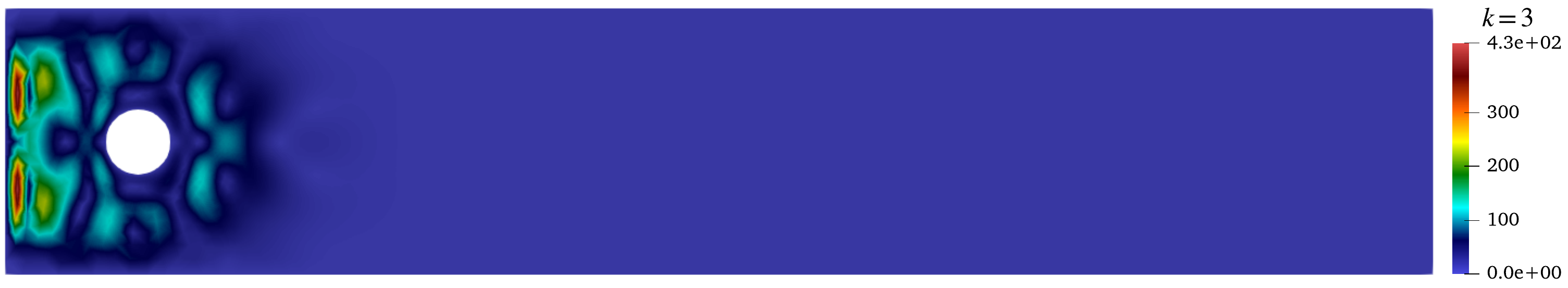}
  \caption{First four space modes of the time series expansion at $\t=0$ ($\Re=4000\,,m=2$).\label{figmodal4000-m2-0}}
 \end{center}
\end{figure}

We end by presenting the case when $\Re =4000$. Here the mesh is built via the \textit{Gmsh} software to take into consideration the turbulent wave behind the cylinder. We plot in \cref{mesh-cyl_Re4000} the vertices of the mesh. It contains 4\,347 vertices and thus a total of 38\,563 \ac{DOFs}.
\cref{figmodal4000-m2-0} presents the first four spatial modes of the time series velocity at $\t=0$ for the case $m=2$, $\C{h}=1$, and $N=75$. We can see in this figure that we have a representation of the velocity field spatial modes that could allow us to apply the \ac{FS} algorithm. \cref{figmodal4000-m2-0} presents the same spatial modes but rather at $\t=3$.

\begin{figure}[!ht]
 \begin{center}
\includegraphics[scale=0.245]{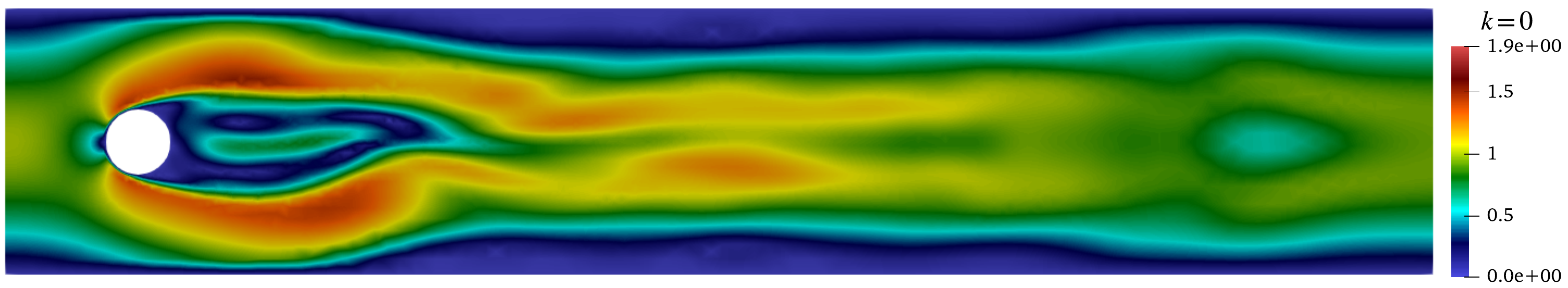}
\includegraphics[scale=0.245]{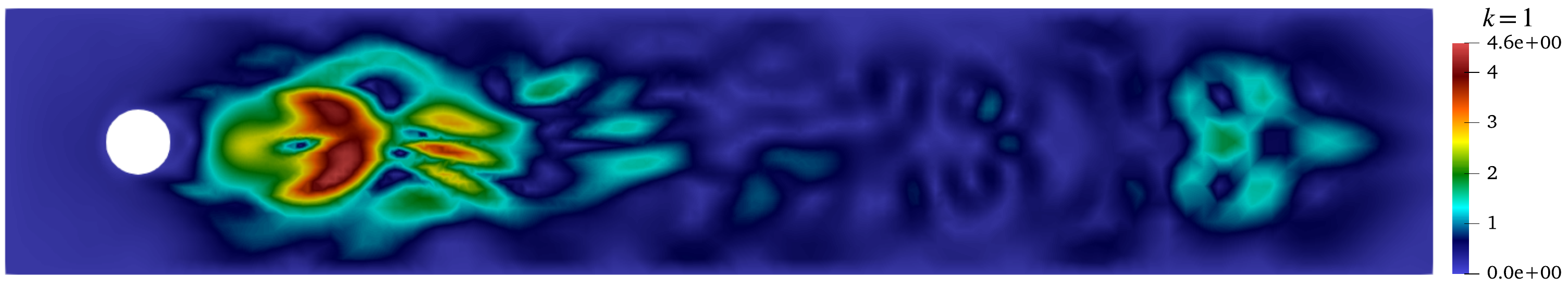}
\includegraphics[scale=0.245]{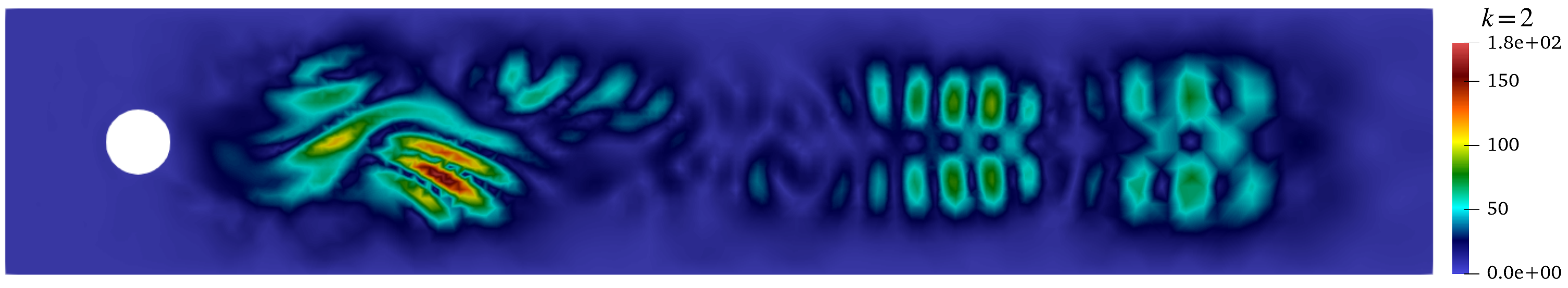}
\includegraphics[scale=0.245]{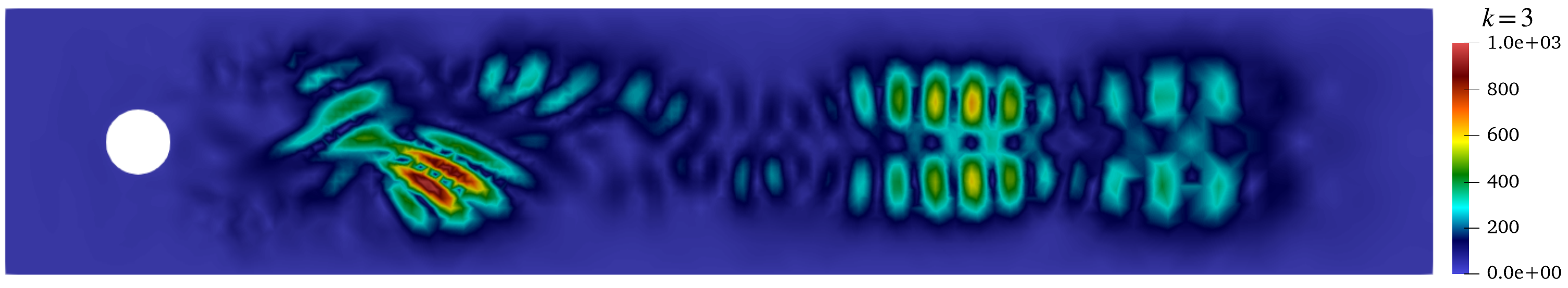}
  \caption{First four space modes of the time series expansion at $\t=3$ ($\Re=4000\,,m=2$).\label{figmodal4000-m2-3}}
 \end{center}
\end{figure}

\begin{figure}[!htp]
 \begin{center}
\includegraphics[scale=0.245]{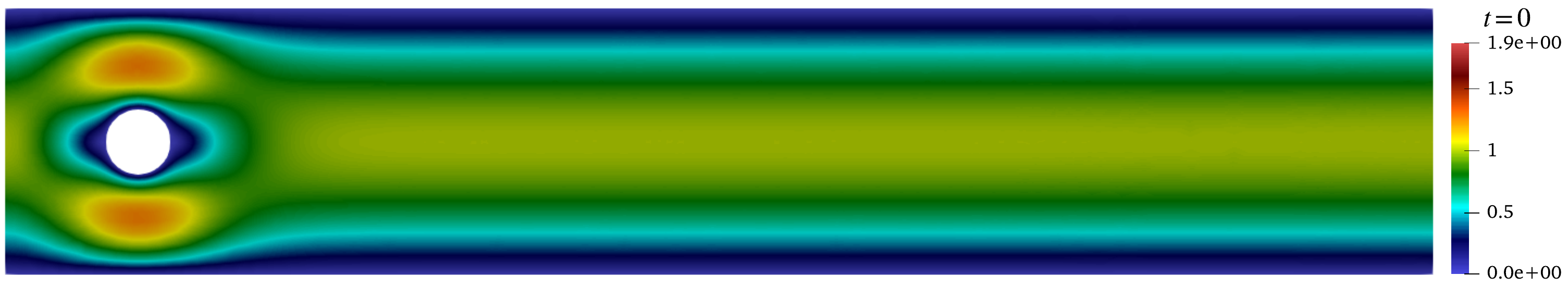}
\includegraphics[scale=0.245]{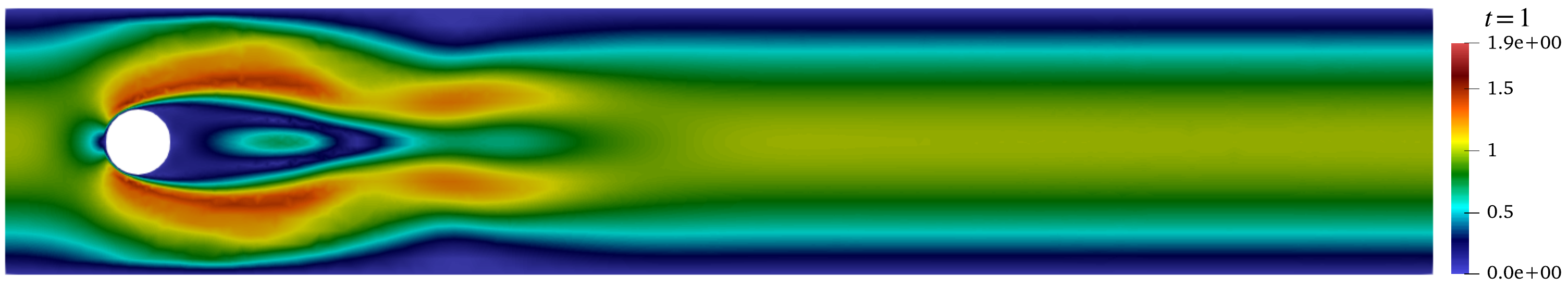}
\includegraphics[scale=0.245]{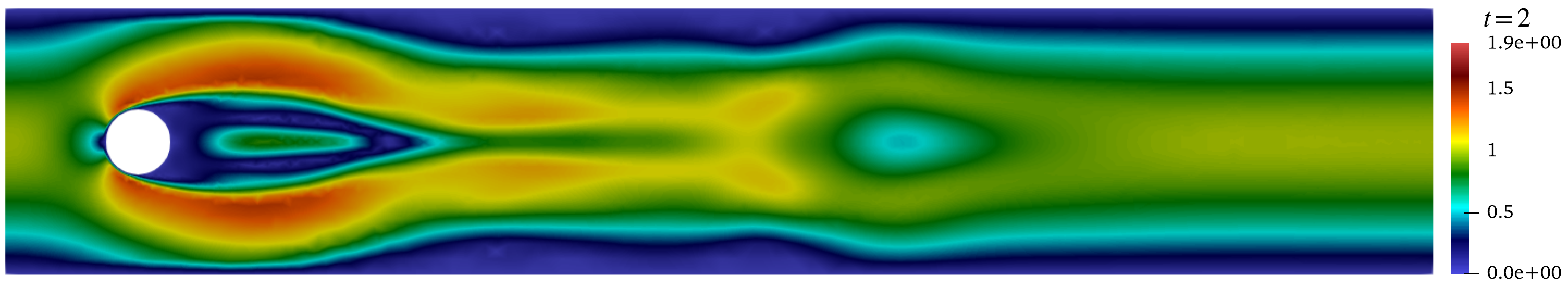}
\includegraphics[scale=0.245]{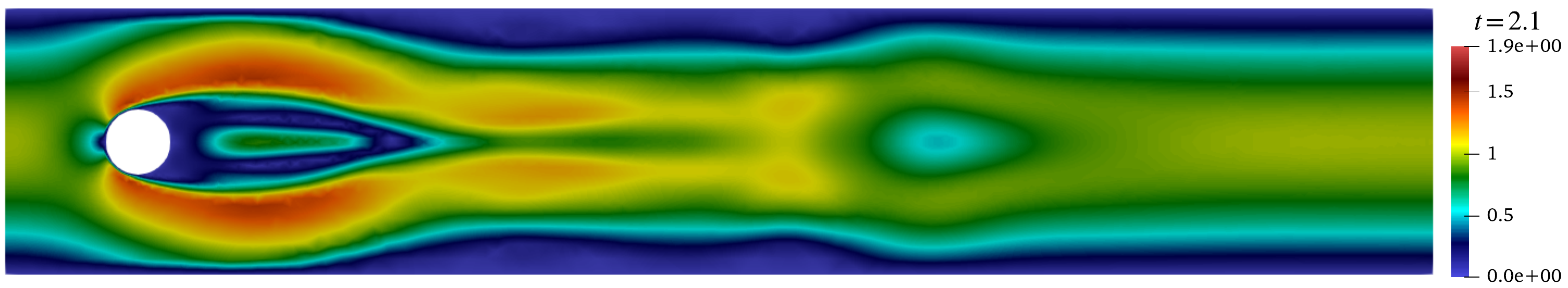}
\includegraphics[scale=0.245]{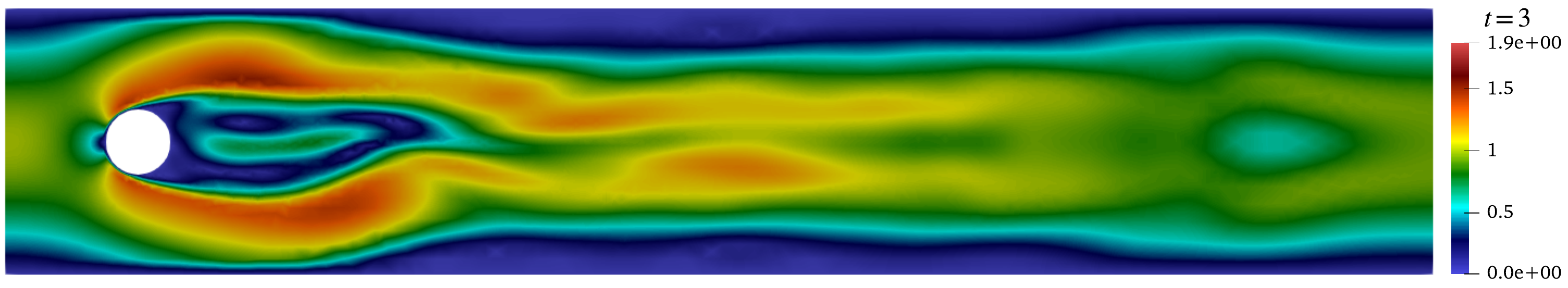}
\includegraphics[scale=0.245]{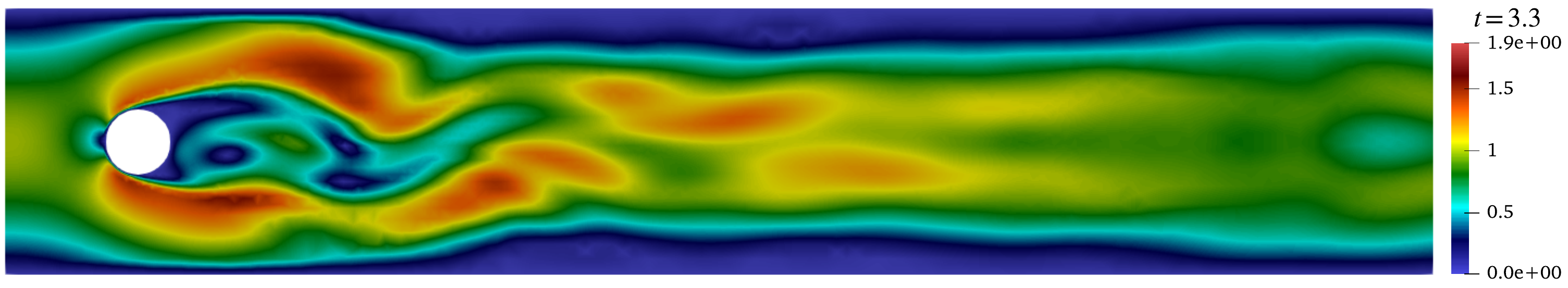}
\includegraphics[scale=0.245]{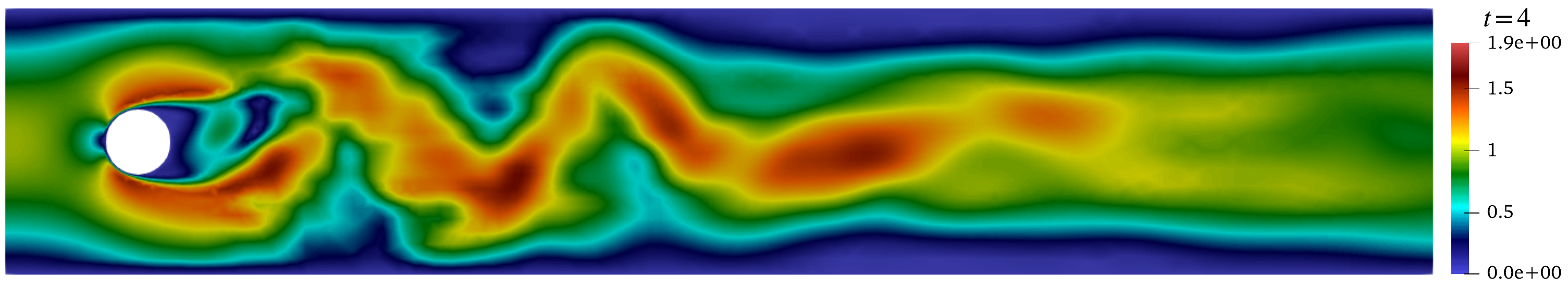}
\includegraphics[scale=0.245]{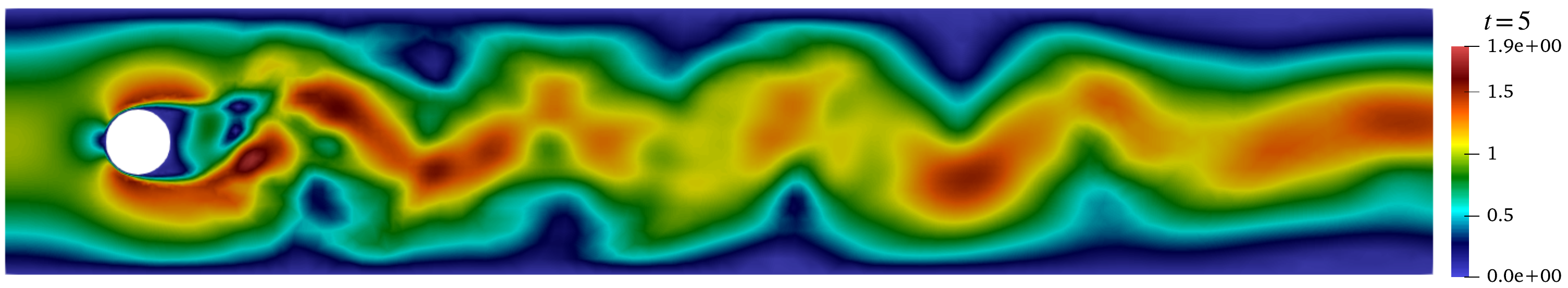}
  \caption{ Snapshots of the velocity field magnitude at different instants of time with stabilization ($\Re = 4000, \,\tau = 0.01\,, m = 2, \, \C{h} = 1$).\label{fig4000-m2-t}}
 \end{center}
\end{figure}

\begin{figure}[!htp]
 \begin{center}
\includegraphics[scale=0.245]{Fig1_cyl_4000_t0.pdf}
\includegraphics[scale=0.245]{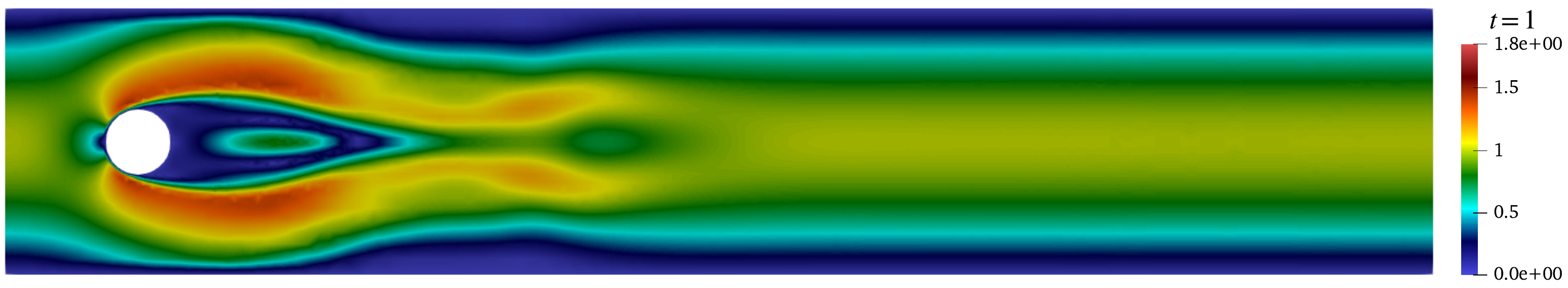}
\includegraphics[scale=0.245]{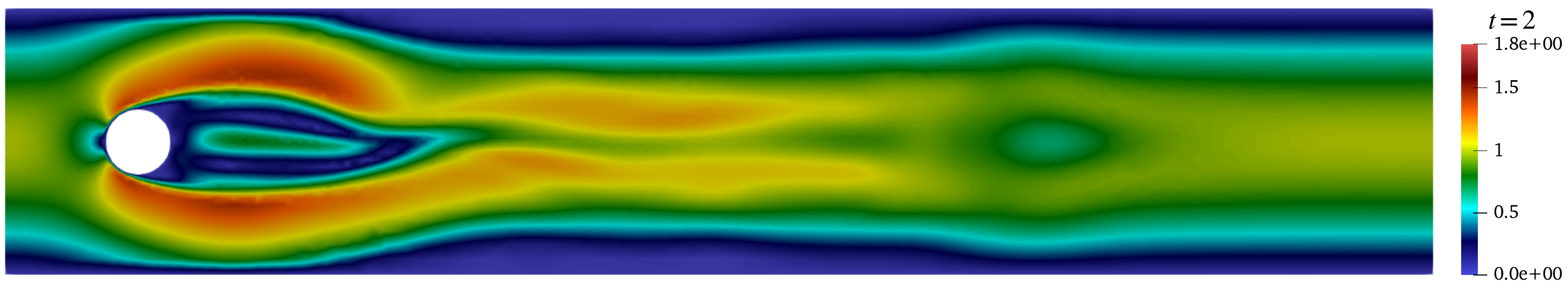}
\includegraphics[scale=0.245]{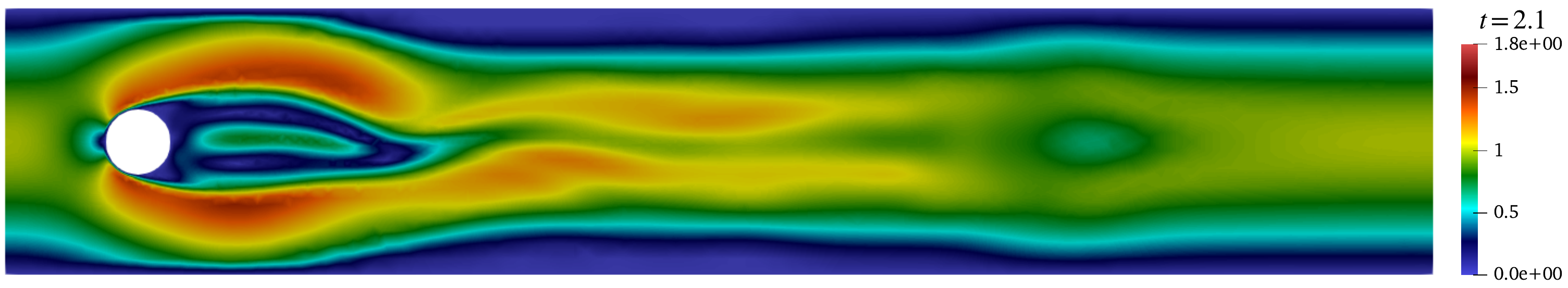}
\includegraphics[scale=0.245]{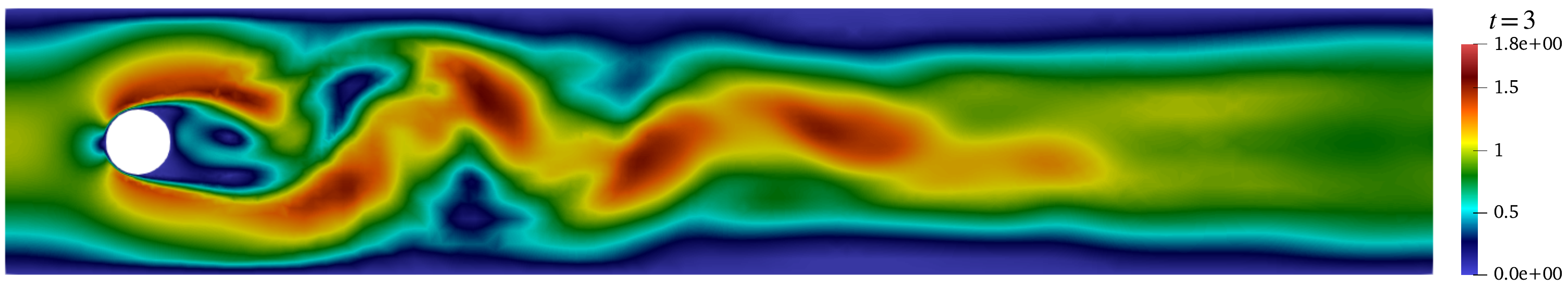}
\includegraphics[scale=0.245]{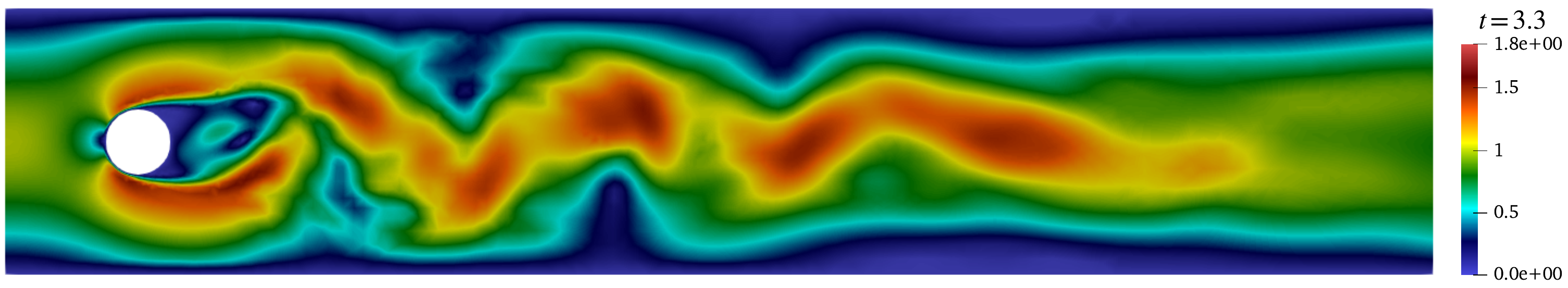}
\includegraphics[scale=0.245]{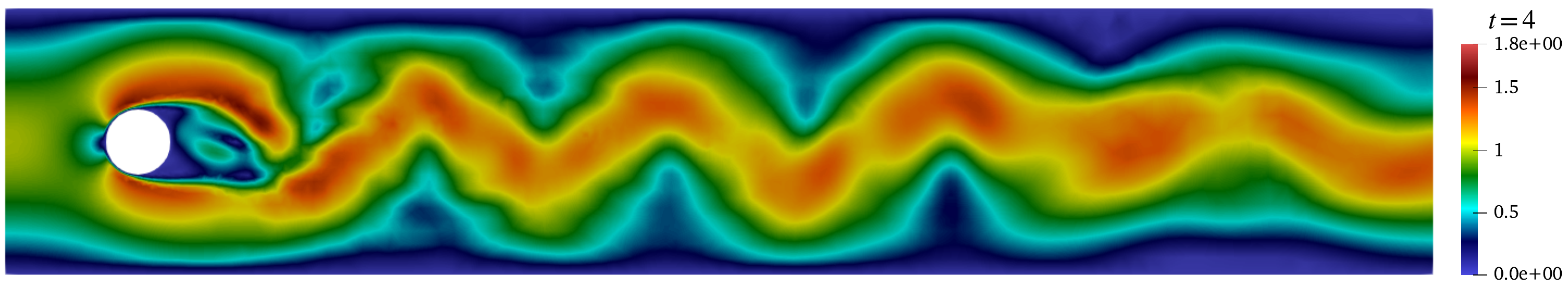}
\includegraphics[scale=0.245]{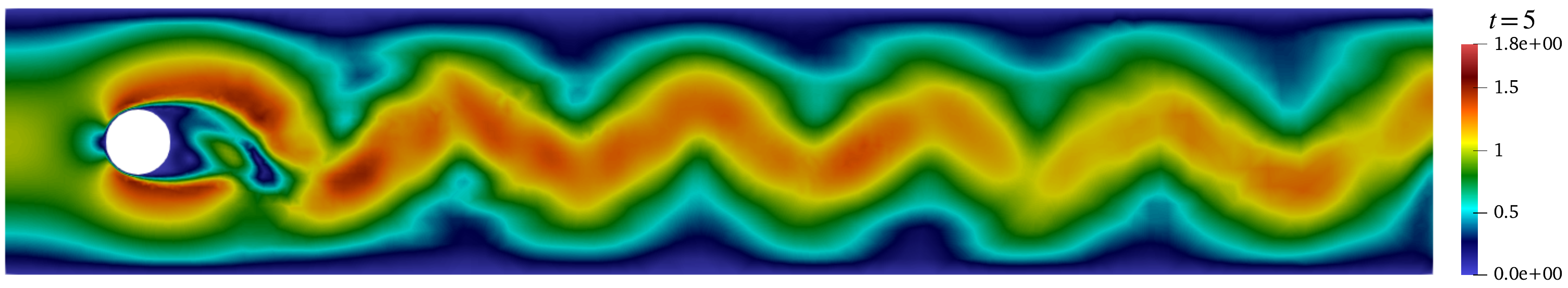}
  \caption{Snapshots of the velocity field magnitude at different instants of time with stabilization ($\Re = 4000, \,\tau = 0.01\,, m = 2.2, \, \C{h} = 1$). \label{fig4000-m2-2-t}}
 \end{center}
\end{figure}

\begin{figure}[!htp]
 \begin{center}
  \includegraphics[height=4cm]{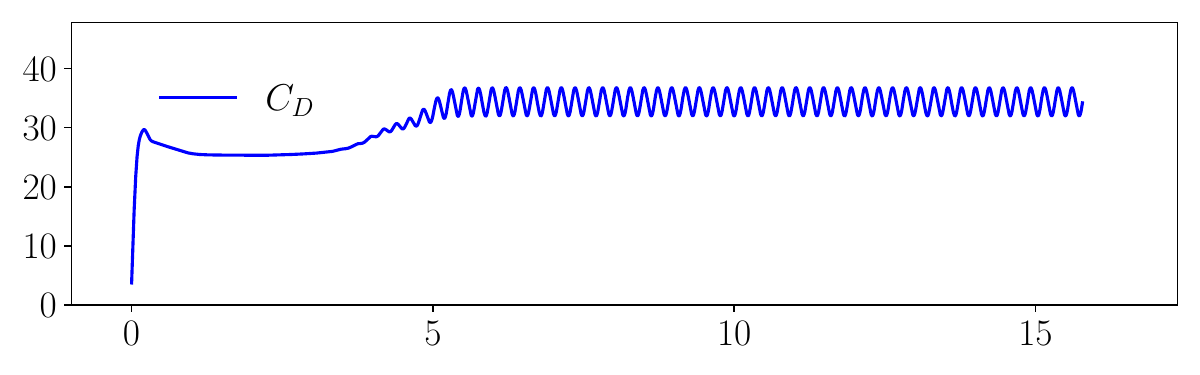}
    \includegraphics[height=4cm]{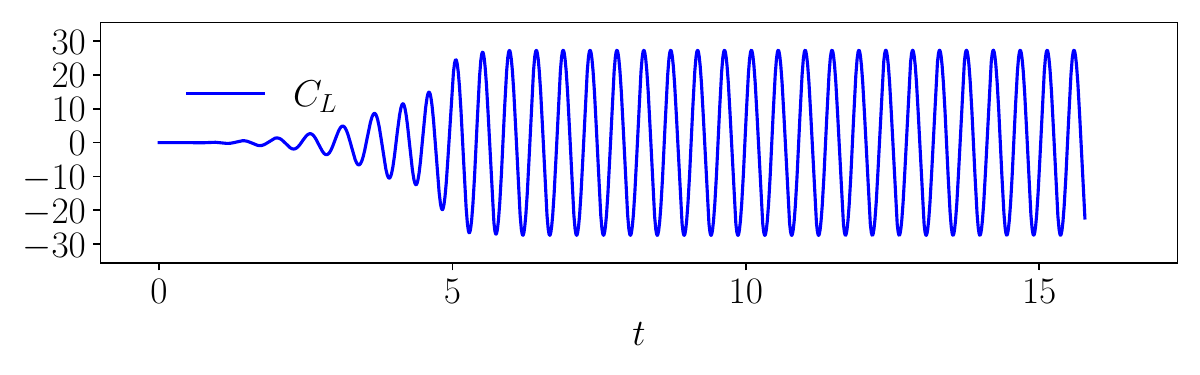}
  \caption{Drag and Lift coefficients, $C_D$ and $C_L$.\label{FigC_5_2_2_0}}
 \end{center}
\end{figure}

We apply now the \ac{FS} algorithm with a time step $\tau=0.01$, $\alpha_0(\bx) = \h(\bx)^m$ with two values of $m\in\{2, 2.2\}$ and $\C{h} = 1$. \cref{fig4000-m2-t} presents the snapshots for the case of $m=2$ and \cref{fig4000-m2-2-t} presents the snapshots when taking $m=2.2$ in the stabilization. The evolution of the drag and lift coefficients are presented in \Cref{FigC_5_2_2_0} in the case where the simulation is ran with the following parameters of the stabilization technique: $m=2$, $\alpha_0 = h^m$, $\alpha_\k =\alpha_0$ and a time step $\tau=0.02$.

\section{Conclusions}
\label{sec8}
In this paper, the \ac{FS} algorithm, which is based on divergent series resummation applied to the time series expansion, is tailored to integrate numerically in time the incompressible \ac{NS} equations, while the \ac{FEMF} is used for space approximation.
The solution is considered to be written in time series expansion form, that is injected in the system to generate a recurrence formula for velocity and pressure terms. The mixed weak formulation is adopted, where the conditions of existence and uniqueness of the solution at every rank of the recurrence formula follows from the Brezzi-conditions. Numerically, the \ac{FEM} is considered to approximate the solution in space variables where the Taylor-Hood \ac{FE} is employed for the velocity-pressure mixed \ac{FE} space.
A convergence study of the proposed method in the framework of \ac{TSE}-\ac{FEM} was presented for the evolutionary Stokes problem.
For simulation tests, the spatial modes of the series are computed up to the rank $\m=5$, where the \ac{FS} algorithm is applied to respect the functions space.

The stabilization technique is employed to enhance computing the terms of the series by adding an artificial diffusion to the left hand side of the recurrence formula, multiplied by a diffusion coefficient. This addition is seen as a regularization of terms by applying the operator $(\I\d-\alpha_0 \Delta)^{-1}$, or its approximation in the \ac{FEMF}. A numerical analysis relates this stabilization term to the enhancing of the condition number of the mass matrix, where the coefficient $\alpha_0$ is found by a minimization problem of the new mass matrix $(\mA+\alpha_0\mK)$. Several cases were considered regarding the dependence of $\alpha_\k$ on $\k$, where the case of $\alpha_\k = \alpha_{0}$ led to a new model resembling closely to \ac{NSa} model.

A series of numerical tests were conducted for the flow past a cylinder problem. First, the Reynolds number is considered to be equal to $400$ and the various space resolutions $N$, various values of the stabilization coefficients $\alpha_0$ and various time steps are considered conducting a practical stability analysis of the algorithm. This was compared with the algorithm without stabilization. The results show an enhancing of the computing of spatial modes of the time series allowing a time step $\mathcal{O}(100)$ bigger than the case without the stabilization.

\section{Perspectives}
\label{sec9}
In this work, the mixed formulation is employed, which involves a large dimensional sparse matrix, thus higher amount of CPU to achieve inversion of matrices. This could be reduced by developing the algorithm and combined with the Chorin projection. The time computation could be also reduced by applying the \ac{FS} on the vertices of the mesh and not only on every \ac{DOFs}. Then, a projection over the space of approximation could be done.
The fixed time step could be also replaced by a variable time step relative to the local dynamic of the problem. This could be done by conducting a priori estimate of the validity of the \ac{FS}-\ac{FEM} algorithm for every truncated series. A stability condition of the time step has also to be established.
Another perspective is the developing of the Padé-vector to produce a compact \ac{BPL} algorithm, allowing to run the full \ac{BPL} algorithm on \ac{HPC} -- as the use of \ac{API} in Python language programming. The application of the \ac{BPL} algorithm with stabilized \ac{TSE} to solve the Cahn-Hilliard equation will be of high interest when using high-order \ac{FEM}, as the equation contains fourth order space derivatives.

\section*{Acknowledgements}
This publication is based upon work supported by the Khalifa University of Science and Technology under Award No. FSU-2023-014.

\appendix

\section{Brezzi's conditions}
\label{append1}
We will state the Brezzi's conditions theorem. It was used in the error estimate of the approximation by the numerical method \ac{TSE}-\ac{FEM} in \cref{sec3}. The proof could be found in \cite{brenner-10}.
\begin{theorem}
Let $ a(\bu, \bv) $ be a continuous bilinear form defined on $ \V \times \V $, and $ b(\bv, \q) $ be a continuous bilinear form defined on $ \V \times \P $. Consider the variational problem for $ (\bu, \p) \in \V \times \P $,

\begin{equation}
\label{Mixed_problem}
\begin{aligned}
a(\bu, \bv) + b(\bv, \p) = \rF[\bv], \quad \forall \bv \in \V,\\
b(\bu, \q) = \rG[\q], \quad \forall \q \in \P,
\end{aligned}
\end{equation}
for $ \rF $ and $ \rG $ continuous linear forms on $ \V $ and $ \P $ respectively. Define the kernel $ \Z $ by
\begin{equation}
\Z = \{ \bu \in \V \ \vert\ b(\bu, \q) = 0 \ \forall \q \in \P \}.
\end{equation}
Assume the following conditions:
\begin{enumerate}
    \item $ a(\bu, \bv) $ is coercive on the kernel $ \Z $ with coercivity constant $ \delta $.
    \item There exists $ \beta > 0 $ such that the inf-sup condition for $ b(\bv, \q) $ holds.
\end{enumerate}
Then there exists a unique solution $ (\bu, \p) $ to the variational problem and we have the stability bound
\begin{equation}
\label{bound_brezzi}
\begin{aligned}
\| \bu \|_\V \leq \frac{1}{\delta} \| \rF \|_{\V'} + \frac{2M}{\delta \beta} \| \rG \|_{\P'},\\
\| \p \|_\P \leq \frac{2M}{\delta \beta} \| \rF \|_{\V'} + \frac{2M^2}{\delta \beta^2} \| \rG \|_{\P'},
\end{aligned}
\end{equation}
where $ M $ is the continuity constant of $ a $.
\end{theorem}

\section{Proof of bound error}
\label{append2}

In this section, we will show the proof the error bound already stated in \cref{error_Vh}. The proof is mostly based on a similar way of Céa Lemma. For more details please refer to the books of Brenner and Scott \cite{brenner-10} and Ciarlet \cite{ciarlet-02}.

\subsection{Proof of  bound errors in \cref{error_Vh}}
\label{append21}
Consider we have the mixed problem \eqref{Mixed_problem} defined in the mixed space $\Vr\times \P$, $a\equiv a_\k$ to be the bilinear form, $\rF \equiv \lins_{\ru_\k}$ with the exact weak solution $\ru_\k$ and  $\rG \equiv 0$. We have $(\ru_{\k+1},\rp_\k)$ verifies it for every $(\bv,\q) \in \Vr\times \P$. As $\Vrh \subset \Vr$ and $\P_h\subset \P$, we choose $(\bv,\q) \in \Vrh\times \P_h$. We denote by $(\reu_{\k+1,h},\rep_{\k,h})$ that verifies the mixed problem for using the exact $\ru_\k$ in the right hand side. By subtracting one from other, we obtain:
\begin{equation}
\begin{aligned}
 a_{\k+1}(\reu_{\k+1,h}-\ru_{\k+1},\bv) - b(\bv,\rep_{\k,h}-\rp_\k) &= 0 ,& \forall \bv \in \Vrh\\
 b(\reu_{\k+1,h}-\ru_{\k+1},\q) &=0,& \forall \q \in \P_h.
\end{aligned}
 \end{equation}
We replace $\ru_{\k+1} = \ru_{\k+1} - \ru_I + \ru_I$ and $\rp_{\k}  =\rp_{\k} - \rp_I + \rp_I$ for $\ru_I,\rp_I \in \Vrh\times \P_h$ and rearrange terms to find:
\begin{equation}
\begin{aligned}
 a_{\k+1}(\reu_{\k+1,h}-\ru_{I},\bv) - b(\bv,\rep_{\k,h} - \rp_{I}) &= \reF_{\k+1,\ru_I,\rp_I}[\bv] ,& \forall \bv \in \Vrh\\
 b(\reu_{\k+1,h}-\ru_{I},\q) &=\rG_{\ru_I}[\q],& \forall \q \in \P_h.
\end{aligned}
\end{equation}
We now use the stability bounds of Brezzi's conditions stated in \ref{append1} to get:
\begin{equation}
 \begin{aligned}
  \left\|\reu_{\k+1,h}-\ru_{I} \right\|_{\Vr} \leqslant \frac{1}{\delta_{\k+1}} \left\| \reF_{\k+1,\ru_I,\rp_I} \right\|_{(\Vr)^\prime} + \frac{2 M_{\k+1}}{\delta_{\k+1}\beta_h} \left\| \rG_{\ru_I}\right\|_{(\P_h)^{\prime}}\\
  \left\| \rep_{\k,h} - \rp_{I} \right\|_{\P_h} \leqslant \frac{2M_{\k+1}}{\delta_{\k+1}\beta_h} \left\| \reF_{\k+1,\ru_I,\rp_I} \right\|_{(\Vr)^\prime} + \frac{2 M_{\k+1}^2}{\delta_{\k+1}\beta_h^2} \left\| \rG_{\ru_I}\right\|_{(\P_h)^{\prime}}.
 \end{aligned}
\end{equation}
With
\begin{equation}
 \begin{aligned}
  \reF_{\k+1,\ru_I,\rp_I}[\bv]& = a_{\k+1}(\ru_{\k+1}-\ru_{I},\bv) - b(\bv,\rp_{\k} - \rp_{I}),& \forall \bv \in \Vrh\\
  \rG_{\ru_I}[\q] &= b(\ru_{\k+1}-\ru_{I},\q),& \forall \q \in \P_h,
 \end{aligned}
\end{equation}
and using continuity of $a_{\k+1}(\cdot,\cdot)$ and $b(\cdot,\cdot)$, we have:
\begin{align}
\begin{aligned}\left\|\reF_{\k+1,\ru_I,\rp_I}\right\|_{(\Vr)\prime} &\leqslant
 M_{\k+1}\|\ru_{\k+1}-u_I\|_{\Vr} + M_b\|\rp_\k-\rp_I\|_{\P},\\
 \|\rG_{u_I}\|_{(\P)\prime}
&\leqslant M_b\|\ru_{\k+1}-u_I\|_{\Vr}.
\end{aligned}
\end{align}
By substituting it and having also:
\begin{equation}
 \beta_h \leqslant \inf\limits_{0\neq \q \in \P} \sup\limits_{0\neq \bv \in \V_\k} \frac{b(\bv,\q)}{\|\bv\|_{\V_\k}\,\|\q\|_{\P}} \leqslant M_b \Longrightarrow \frac{M_b}{\beta_h}\geqslant 1,
\end{equation}
we obtain:
\begin{equation}
 \begin{aligned}
  \left\|\reu_{\k+1,h}-\ru_{I} \right\|_{\Vr} &\leqslant \frac{3MM_b}{\delta\beta_h} \|\ru_{\k+1}-u_I\|_{\Vr} + \frac{M_b}{\delta(\k+1)}\|\rp_\k-\rp_I\|_{\P}\\
  \left\| \rep_{\k,h} - \rp_{I} \right\|_{\P_h}& \leqslant
 \frac{4 M^2{(\k+1)}M_b}{\delta\beta_h^2} \|\ru_{\k+1}-u_I\|_{\Vr}+  \frac{2MM_b}{\delta\beta_h} \|\rp_\k-\rp_I\|_{\P} .
 \end{aligned}
\end{equation}
By using the triangle inequality, we have:
\begin{equation}
 \begin{aligned}
  \left\|\reu_{\k+1,h}-\ru_{I} \right\|_{\Vr} &\leqslant \frac{4MM_b}{\delta\beta_h} \|\ru_{\k+1}-u_I\|_{\Vr} + \frac{M_b}{\delta(\k+1)}\|\rp_\k-\rp_I\|_{\P}\\
  \left\| \rep_{\k,h} - \rp_{I} \right\|_{\P_h}& \leqslant
 \frac{4 M^2{(\k+1)}M_b}{\delta\beta_h^2} \|\ru_{\k+1}-u_I\|_{\Vr}+  \frac{3MM_b}{\delta\beta_h} \|\rp_\k-\rp_I\|_{\P} .
 \end{aligned}
\end{equation}
\subsection{Proof of  bound errors in \cref{error_Vhc}}
\label{append22}
We return to the original problem and denote by $(\ru_{\k+1},\rp_{\k,h})$ the numerical solution of the mixed problem when considering $\ru_{\k,h}$ in the right hand side; \emph{i.e.} using $\lins_{\ru_{\k,h}}$. In this case, the orthogonality of the error is represented by:
\begin{equation}
\begin{aligned}
 a_{\k+1}(\ru_{\k+1,h}-\ru_{I},\bv) - b(\bv,\rp_{\k,h} - \rp_{I}) &= \rF_{\k+1,\ru_I,\rp_I}[\bv] ,& \forall \bv \in \Vrh\\
 b(\ru_{\k+1,h}-\ru_{I},\q) &=\rG_{\ru_I}[\q],& \forall \q \in \P_h.
\end{aligned}
\end{equation}
With
\begin{equation}
 \begin{aligned}
  \rF_{\k+1,\ru_I,\rp_I}[\bv]& = a_{\k+1}(\ru_{\k+1}-\ru_{I},\bv) - b(\bv,\rp_{\k} - \rp_{I}) + \lins_{\reu_{k,h}-\ru_\k}(\bv),& \forall \bv \in \Vrh .\end{aligned}
\end{equation}
Following the same procedure we used to obtain the error bounds of $\|\reu_{k,h}-\ru_\k\|_{\Vr}$, using continuity of $a_{\k+1}$, $b$ and $\lins$,  we can obtain:
\begin{equation}
 \begin{aligned}
  \left\|\ru_{\k+1,h}-\ru_{I} \right\|_{\Vr} &\leqslant \frac{3MM_b}{\delta\beta_h} \|\ru_{\k+1}-u_I\|_{\Vr} + \frac{M_b}{\delta(\k+1)}\|\rp_\k-\rp_I\|_{\P}\\
   & \quad + \frac{1}{\delta_{\k+1}\Re} \left\|\ru_{\k}-\reu_{\k,h} \right\|_{\Vr} \\
  \left\| \rp_{\k,h} - \rp_{I} \right\|_{\P_h}& \leqslant
 \frac{4 M^2{(\k+1)}M_b}{\delta\beta_h^2} \|\ru_{\k+1}-u_I\|_{\Vr}+  \frac{2MM_b}{\delta\beta_h} \|\rp_\k-\rp_I\|_{\P} \\
 & \quad + \frac{2M^2_{\k+1}}{\delta_{\k+1} \beta} \frac{1}{\Re}  \left\|\ru_{\k}-\reu_{\k,h} \right\|_{\Vr}.
 \end{aligned}
\end{equation}
Now, we use the error bound for $\|\reu_{k,h}-\ru_\k\|_{\Vr}$ and the triangle inequality, we obtain:
\begin{equation}
 \begin{aligned}
  \left\|\ru_{\k+1}-\ru_{\k+1,h} \right\|_{\Vr} &\leqslant
  \frac{4MM_b}{\delta\beta}\left[1 + \frac{1}{(\k+1)\delta\Re} \right]O(h^s) \\
  &\quad+ \frac{M_b}{(\k+1)\delta} \left[ 1+\frac{1}{\k\delta \Re}\right] O(h^r)\\
  \left\|\rp_{\k}-\rp_{\k,h} \right\|_{\Vr} & \leqslant \frac{4M^2M_b}{\delta\beta^2}\left[\k+1 + \frac{2}{\delta\Re} \right]O(h^s) \\
  &\quad+ \frac{2MM_b}{\delta\beta} \left[ 1+\frac{1}{\k\delta \Re}\right] O(h^r)\\
 \end{aligned}
\end{equation}
which completes the proof of error bounds.

\section*{List of abbreviations}
\begin{acronym}[TDMA]
\acro{LBB}{Ladyzhenskaya–Babuška–Brezzi}
\acro{DMP}{Discrete Maximum Principle}
\acro{SUPG}{Streamline Upwind Petrov-Galerkin}
\acro{DOFs}{Degrees of Freedom}
\acro{ODE}{Ordinary Differential Equation}
\acro{PDE}{Partial Differential Equation}
\acro{CFD}{Computational Fluid Dynamics}
\acro{TSE}{Time Series Expansion}
\acro{BPL}{Borel-Padé-Laplace}
\acro{BL}{Borel-Laplace}
\acro{FE}{Finite Element}
\acro{FEM}{Finite Element Method}
\acro{FDM}{Finite Difference Methods}
\acro{FVM}{Finite Volume Methods}
\acro{PGD}{Proper Generalized Decomposition}
\acro{FEMF}{FEM framework}
\acro{NS}{Navier-Stokes}
\acro{NSa}[NS-$\alpha$]{Navier-Stokes-alpha}
\acro{IVP}{Initial Value Problem}
\acro{CP}{Cauchy Problem}
\acro{LSM}{Level Set Method}
\acro{Pa-ap}{Padé Approximant}
\acro{SVD}{Singular Value Decomposition}
\acro{GL}{Gau\ss\,\,-Laguerre}
\acro{LMS}{Linear Multi-Step}
\acro{BDF}{Backward Difference Formula}
\acro{BDFk}{\ac{BDF} of order $k$}
\acro{RK}{Runge-Kutta}
\acro{ERK}{Embedded-Runge-Kutta}
\acro{IRK}{Implicit \ac{RK}}
\acro{GRK}{Gau\ss-Runge-Kutta}
\acro{ETD}{Exponential-Time Difference}
\acro{DSR}{Divergent Series Resummation}
\acro{GFS}{Generalized Factorial Series}
\acro{ATS}{Adaptive Time Stepping}
\acro{HPC}{High Performance Computing}
\acro{API}{Application Programming Interface}
\acro{GFS}{Generalized Factorial Series}
\acro{FS}{Factorial Series}
\end{acronym}

 \bibliographystyle{elsarticle-num-names}
 \bibliography{references}

\begin{thebibliography}{73}
\expandafter\ifx\csname natexlab\endcsname\relax\def\natexlab#1{#1}\fi
\providecommand{\url}[1]{\texttt{#1}}
\providecommand{\href}[2]{#2}
\providecommand{\path}[1]{#1}
\providecommand{\DOIprefix}{doi:}
\providecommand{\ArXivprefix}{arXiv:}
\providecommand{\URLprefix}{URL: }
\providecommand{\Pubmedprefix}{pmid:}
\providecommand{\doi}[1]{\href{http://dx.doi.org/#1}{\path{#1}}}
\providecommand{\Pubmed}[1]{\href{pmid:#1}{\path{#1}}}
\providecommand{\bibinfo}[2]{#2}
\ifx\xfnm\relax \def\xfnm[#1]{\unskip,\space#1}\fi
\bibitem[{Sagot(1998)}]{CFD_sagot}
\bibinfo{author}{P.~Sagot}, \bibinfo{title}{Introduction a la simulation des
  grandes échelles pour les écoulements de fluide incompressible},
  volume~\bibinfo{volume}{30} of \textit{\bibinfo{series}{Series Title
  Mathématiques et Applications}}, \bibinfo{edition}{1} ed.,
  \bibinfo{publisher}{Springer Berlin}, \bibinfo{address}{Heidelberg},
  \bibinfo{year}{1998}.
\bibitem[{Wendt(2010)}]{CFD-book-09}
\bibinfo{editor}{J.~F. Wendt} (Ed.), \bibinfo{title}{Computational Fluid
  Dynamics: An Introduction}, \bibinfo{edition}{3} ed.,
  \bibinfo{publisher}{Springer}, \bibinfo{address}{Berlin, Heidelberg},
  \bibinfo{year}{2010}.
  \DOIprefix\doi{https://doi.org/10.1007/978-3-540-85056-4}.
\bibitem[{Bazilevs et~al.(2013)Bazilevs, Takizawa, and Tezduyar}]{book-CFSI-12}
\bibinfo{author}{Y.~Bazilevs}, \bibinfo{author}{K.~Takizawa},
  \bibinfo{author}{T.~E. Tezduyar}, \bibinfo{title}{Computational
  Fluid–Structure Interaction}, \bibinfo{publisher}{John Wiley \& Sons, Ltd},
  \bibinfo{address}{West Sussex, United Kingdom}, \bibinfo{year}{2013}.
  \DOIprefix\doi{https://doi.org/10.1002/9781118483565.fmatter}.
\bibitem[{Fambri(2020)}]{Archives-1}
\bibinfo{author}{F.~Fambri},
\newblock \bibinfo{title}{Discontinuous galerkin methods for compressible and
  incompressible flows on space–time adaptive meshes: Toward a novel family
  of efficient numerical methods for fluid dynamics},
\newblock \bibinfo{journal}{Archives of Computational Methods in Engineering}
  \bibinfo{volume}{27} (\bibinfo{year}{2020}) \bibinfo{pages}{199--283}.
  \DOIprefix\doi{https://doi.org/10.1007/s11831-018-09308-6}.
\bibitem[{Wang et~al.(2016)Wang, Chen, Wang, Liao, Zhu, and
  Li}]{SPH-archives-15}
\bibinfo{author}{Z.-B. Wang}, \bibinfo{author}{R.~Chen},
  \bibinfo{author}{H.~Wang}, \bibinfo{author}{Q.~Liao},
  \bibinfo{author}{X.~Zhu}, \bibinfo{author}{S.-Z. Li},
\newblock \bibinfo{title}{An overview of smoothed particle hydrodynamics for
  simulating multiphase flow},
\newblock \bibinfo{journal}{Applied Mathematical Modelling}
  \bibinfo{volume}{40} (\bibinfo{year}{2016}) \bibinfo{pages}{9625--9655}.
  \DOIprefix\doi{https://doi.org/10.1016/j.apm.2016.06.030}.
\bibitem[{Argyropoulos and Markatos(2015)}]{NM-TF-15}
\bibinfo{author}{C.~Argyropoulos}, \bibinfo{author}{N.~Markatos},
\newblock \bibinfo{title}{Recent advances on the numerical modelling of
  turbulent flows},
\newblock \bibinfo{journal}{Applied Mathematical Modelling}
  \bibinfo{volume}{39} (\bibinfo{year}{2015}) \bibinfo{pages}{693--732}.
  \DOIprefix\doi{https://doi.org/10.1016/j.apm.2014.07.001}.
\bibitem[{Hu et~al.(2019)Hu, Guo, Hu, Negrut, Xu, and Pan}]{ASPH-5}
\bibinfo{author}{W.~Hu}, \bibinfo{author}{G.~Guo}, \bibinfo{author}{X.~Hu},
  \bibinfo{author}{D.~Negrut}, \bibinfo{author}{Z.~Xu},
  \bibinfo{author}{W.~Pan},
\newblock \bibinfo{title}{A consistent spatially adaptive smoothed particle
  hydrodynamics method for fluid–structure interactions},
\newblock \bibinfo{journal}{Computer Methods in Applied Mechanics and
  Engineering} \bibinfo{volume}{347} (\bibinfo{year}{2019})
  \bibinfo{pages}{402--424}.
  \DOIprefix\doi{https://doi.org/10.1016/j.cma.2018.10.049}.
\bibitem[{Ashour et~al.(2023)Ashour, Valizadeh, and
  Rabczuk}]{phase-NS-hydro-23}
\bibinfo{author}{M.~Ashour}, \bibinfo{author}{N.~Valizadeh},
  \bibinfo{author}{T.~Rabczuk},
\newblock \bibinfo{title}{Phase-field navier--stokes model for vesicle doublets
  hydrodynamics in incompressible fluid flow},
\newblock \bibinfo{journal}{Computer Methods in Applied Mechanics and
  Engineering} \bibinfo{volume}{412} (\bibinfo{year}{2023})
  \bibinfo{pages}{116063}.
  \DOIprefix\doi{https://doi.org/10.1016/j.cma.2023.116063}.
\bibitem[{Barbeau et~al.(2024)Barbeau, Golshan, Deng, Étienne, Béguin, and
  Blais}]{IBM-CFD-24}
\bibinfo{author}{L.~Barbeau}, \bibinfo{author}{S.~Golshan},
  \bibinfo{author}{J.~Deng}, \bibinfo{author}{S.~Étienne},
  \bibinfo{author}{C.~Béguin}, \bibinfo{author}{B.~Blais},
\newblock \bibinfo{title}{High-order moving immersed boundary and its
  application to a resolved cfd-dem model},
\newblock \bibinfo{journal}{Computers \& Fluids} \bibinfo{volume}{268}
  (\bibinfo{year}{2024}) \bibinfo{pages}{106094}.
  \DOIprefix\doi{https://doi.org/10.1016/j.compfluid.2023.106094}.
\bibitem[{Eymard et~al.(2018)Eymard, Feron, and Guichard}]{NumScheme-INS-18}
\bibinfo{author}{R.~Eymard}, \bibinfo{author}{P.~Feron},
  \bibinfo{author}{C.~Guichard},
\newblock \bibinfo{title}{Family of convergent numerical schemes for the
  incompressible navier-stokes equations},
\newblock \bibinfo{journal}{Mathematics and Computers in Simulation}
  \bibinfo{volume}{144} (\bibinfo{year}{2018}) \bibinfo{pages}{196--218}.
  \DOIprefix\doi{https://doi.org/10.1016/j.matcom.2017.08.003}.
\bibitem[{Taylor and Hood(1973)}]{taylor-hood-73}
\bibinfo{author}{C.~Taylor}, \bibinfo{author}{P.~Hood},
\newblock \bibinfo{title}{A numerical solution of the navier-stokes equations
  using the finite element technique},
\newblock \bibinfo{journal}{Computers \& Fluids} \bibinfo{volume}{1}
  (\bibinfo{year}{1973}) \bibinfo{pages}{73--100}.
  \DOIprefix\doi{https://doi.org/10.1016/0045-7930(73)90027-3}.
\bibitem[{Cooke and Blanchard(1980)}]{FEM-HO-NS-80}
\bibinfo{author}{C.~H. Cooke}, \bibinfo{author}{D.~K. Blanchard},
\newblock \bibinfo{title}{A higher order finite element algorithm for the
  unsteady navier-stokes equations},
\newblock \bibinfo{journal}{Mathematics and Computers in Simulation}
  \bibinfo{volume}{22} (\bibinfo{year}{1980}) \bibinfo{pages}{127--132}.
  \DOIprefix\doi{https://doi.org/10.1016/0378-4754(80)90007-5}.
\bibitem[{Zienkiewicz et~al.(2013)Zienkiewicz, Taylor, and Zhu}]{Archives-4}
\bibinfo{author}{O.~Zienkiewicz}, \bibinfo{author}{R.~Taylor},
  \bibinfo{author}{J.~Zhu}, \bibinfo{title}{The Finite Element Method Set:Its
  Basis and Fundamentals}, \bibinfo{edition}{7} ed.,
  \bibinfo{publisher}{Elsevier-Butterworth-Heinemann},
  \bibinfo{address}{Oxford}, \bibinfo{year}{2013}.
  \DOIprefix\doi{https://doi.org/10.1016/C2009-0-24909-9}.
\bibitem[{Dettmer and Perić(2006)}]{FSI-06}
\bibinfo{author}{W.~Dettmer}, \bibinfo{author}{D.~Perić},
\newblock \bibinfo{title}{A computational framework for fluid–structure
  interaction: Finite element formulation and applications},
\newblock \bibinfo{journal}{Computer Methods in Applied Mechanics and
  Engineering} \bibinfo{volume}{195} (\bibinfo{year}{2006})
  \bibinfo{pages}{5754--5779}.
  \DOIprefix\doi{https://doi.org/10.1016/j.cma.2005.10.019},
  \bibinfo{note}{john H. Argyris Memorial Issue. Part II}.
\bibitem[{Logg et~al.(2012)Logg, Mardal, Wells et~al.}]{LoggEtal2012}
\bibinfo{author}{A.~Logg}, \bibinfo{author}{K.-A. Mardal},
  \bibinfo{author}{G.~N. Wells}, et~al., \bibinfo{title}{Automated Solution of
  Differential Equations by the Finite Element Method},
  \bibinfo{publisher}{Springer}, \bibinfo{year}{2012}.
  \DOIprefix\doi{10.1007/978-3-642-23099-8}.
\bibitem[{Kuzmin and Hämäläinen(2014)}]{FEM-CFD-14}
\bibinfo{author}{D.~Kuzmin}, \bibinfo{author}{J.~Hämäläinen},
  \bibinfo{title}{Finite Element Methods for Computational Fluid Dynamics: A
  Practical Guide}, \bibinfo{publisher}{Society for Industrial and Applied
  Mathematics}, \bibinfo{address}{Philadelphia, PA}, \bibinfo{year}{2014}.
  \DOIprefix\doi{10.1137/1.9781611973617}.
\bibitem[{Zienkiewicz et~al.(2014)Zienkiewicz, Taylor, and
  Nithiarasu}]{FEM-FD-14}
\bibinfo{author}{O.~Zienkiewicz}, \bibinfo{author}{R.~Taylor},
  \bibinfo{author}{P.~Nithiarasu}, \bibinfo{title}{The Finite Element Method
  for Fluid Dynamics}, \bibinfo{edition}{7} ed., \bibinfo{publisher}{Elsevier},
  \bibinfo{address}{Amsterdam}, \bibinfo{year}{2014}.
  \DOIprefix\doi{https://doi.org/10.1016/C2009-0-26328-8}.
\bibitem[{Fujita and Kato(1964)}]{NS-IVP-64}
\bibinfo{author}{H.~Fujita}, \bibinfo{author}{T.~Kato},
\newblock \bibinfo{title}{On the navier-stokes initial value problem. i},
\newblock \bibinfo{journal}{Archive for Rational Mechanics and Analysis}
  \bibinfo{volume}{16} (\bibinfo{year}{1964}) \bibinfo{pages}{269--315}.
  \DOIprefix\doi{10.1007/BF00276188}.
\bibitem[{Boukir et~al.(1994)Boukir, Maday, and Métivet}]{TIS-NS-CM-94}
\bibinfo{author}{K.~Boukir}, \bibinfo{author}{Y.~Maday},
  \bibinfo{author}{B.~Métivet},
\newblock \bibinfo{title}{A high order characteristics method for the
  incompressible navier—stokes equations},
\newblock \bibinfo{journal}{Computer Methods in Applied Mechanics and
  Engineering} \bibinfo{volume}{116} (\bibinfo{year}{1994})
  \bibinfo{pages}{211--218}.
  \DOIprefix\doi{https://doi.org/10.1016/S0045-7825(94)80025-1}.
\bibitem[{Si et~al.(2016)Si, Wang, and Sun}]{TIS-NS-CM-16}
\bibinfo{author}{Z.~Si}, \bibinfo{author}{J.~Wang}, \bibinfo{author}{W.~Sun},
\newblock \bibinfo{title}{Unconditional stability and error estimates of
  modified characteristics fems for the navier–stokes equations},
\newblock \bibinfo{journal}{Numerische Mathematik} \bibinfo{volume}{134}
  (\bibinfo{year}{2016}) \bibinfo{pages}{139--161}.
  \DOIprefix\doi{https://doi.org/10.1007/s00211-015-0767-9}.
\bibitem[{Bijl et~al.(2001)Bijl, Carpenter, and Vatsa}]{TIS-NS-01}
\bibinfo{author}{H.~Bijl}, \bibinfo{author}{M.~Carpenter},
  \bibinfo{author}{V.~Vatsa}, \bibinfo{title}{Time integration schemes for the
  unsteady Navier-Stokes equations}, \bibinfo{publisher}{AIAA},
  \bibinfo{address}{Anaheim, CA}, \bibinfo{year}{2001}.
  \DOIprefix\doi{10.2514/6.2001-2612}.
\bibitem[{Guesmi et~al.(2023)Guesmi, Grotteschi, and
  Stiller}]{Assessment-TIS-IMEX-23}
\bibinfo{author}{M.~Guesmi}, \bibinfo{author}{M.~Grotteschi},
  \bibinfo{author}{J.~Stiller},
\newblock \bibinfo{title}{Assessment of high-order imex methods for
  incompressible flow},
\newblock \bibinfo{journal}{International Journal for Numerical Methods in
  Fluids} \bibinfo{volume}{95} (\bibinfo{year}{2023})
  \bibinfo{pages}{954--978}. \DOIprefix\doi{https://doi.org/10.1002/fld.5177}.
\bibitem[{Kim and Bak(2021)}]{TIS-NS-24}
\bibinfo{author}{P.~Kim}, \bibinfo{author}{S.~Bak},
\newblock \bibinfo{title}{Algorithm for a cost-reducing time-integration scheme
  for solving incompressible navier–stokes equations},
\newblock \bibinfo{journal}{Computer Methods in Applied Mechanics and
  Engineering} \bibinfo{volume}{373} (\bibinfo{year}{2021})
  \bibinfo{pages}{113546}.
  \DOIprefix\doi{https://doi.org/10.1016/j.cma.2020.113546}.
\bibitem[{Feng et~al.(2011)Feng, He, and Liu}]{TIS-NS-11}
\bibinfo{author}{X.~Feng}, \bibinfo{author}{Y.~He}, \bibinfo{author}{D.~Liu},
\newblock \bibinfo{title}{Convergence analysis of an implicit fractional-step
  method for the incompressible navier–stokes equations},
\newblock \bibinfo{journal}{Applied Mathematical Modelling}
  \bibinfo{volume}{35} (\bibinfo{year}{2011}) \bibinfo{pages}{5856--5871}.
  \DOIprefix\doi{https://doi.org/10.1016/j.apm.2011.05.042}.
\bibitem[{Deteix et~al.(2022)Deteix, {Ndetchoua Kouamo}, and
  Yakoubi}]{TIS-NS-FTS-22}
\bibinfo{author}{J.~Deteix}, \bibinfo{author}{G.~{Ndetchoua Kouamo}},
  \bibinfo{author}{D.~Yakoubi},
\newblock \bibinfo{title}{A new energy stable fractional time stepping scheme
  for the navier--stokes/allen--cahn diffuse interface model},
\newblock \bibinfo{journal}{Computer Methods in Applied Mechanics and
  Engineering} \bibinfo{volume}{393} (\bibinfo{year}{2022})
  \bibinfo{pages}{114759}.
  \DOIprefix\doi{https://doi.org/10.1016/j.cma.2022.114759}.
\bibitem[{Anselmann and Bause(2021)}]{Nitsche-NS-21}
\bibinfo{author}{M.~Anselmann}, \bibinfo{author}{M.~Bause},
\newblock \bibinfo{title}{Higher order galerkin-collocation time discretization
  with nitsche’s method for the navier-stokes equations},
\newblock \bibinfo{journal}{Mathematics and Computers in Simulation}
  \bibinfo{volume}{189} (\bibinfo{year}{2021}) \bibinfo{pages}{141--162}.
  \DOIprefix\doi{https://doi.org/10.1016/j.matcom.2020.10.027},
  \bibinfo{note}{mATCOM Special Issue: Modelling 2019: The 6th International
  Conference on Mathematical Modelling and Computational Methods in Applied
  Sciences and Engineering}.
\bibitem[{Li et~al.(2020)Li, Ju, and Si}]{TIS-NS-20}
\bibinfo{author}{S.-J. Li}, \bibinfo{author}{L.~Ju}, \bibinfo{author}{H.~Si},
  \bibinfo{title}{Adaptive Exponential Time Integration of the Navier-Stokes
  Equations}, \bibinfo{publisher}{AIAA}, \bibinfo{address}{Orlando},
  \bibinfo{year}{2020}. \DOIprefix\doi{10.2514/6.2020-2033}.
\bibitem[{Cheng and Wang(2016)}]{TIS-NS-LMS-16}
\bibinfo{author}{K.~Cheng}, \bibinfo{author}{C.~Wang},
\newblock \bibinfo{title}{Long time stability of high order multistep numerical
  schemes for two-dimensional incompressible navier--stokes equations},
\newblock \bibinfo{journal}{SIAM Journal on Numerical Analysis}
  \bibinfo{volume}{54} (\bibinfo{year}{2016}) \bibinfo{pages}{3123--3144}.
  \DOIprefix\doi{10.1137/16M1061588}.
\bibitem[{Kim et~al.(2020)Kim, Kim, Piao, and Bak}]{TIS-NS-BDF-20}
\bibinfo{author}{P.~Kim}, \bibinfo{author}{D.~Kim}, \bibinfo{author}{X.~Piao},
  \bibinfo{author}{S.~Bak},
\newblock \bibinfo{title}{A completely explicit scheme of cauchy problem in
  bslm for solving the navier–stokes equations},
\newblock \bibinfo{journal}{Journal of Computational Physics}
  \bibinfo{volume}{401} (\bibinfo{year}{2020}) \bibinfo{pages}{109028}.
  \DOIprefix\doi{https://doi.org/10.1016/j.jcp.2019.109028}.
\bibitem[{Breckling et~al.(2024)Breckling, Fiordilino, Reyes, and
  Shields}]{TIS-NS-BDF-24}
\bibinfo{author}{S.~Breckling}, \bibinfo{author}{J.~Fiordilino},
  \bibinfo{author}{J.~Reyes}, \bibinfo{author}{S.~Shields},
\newblock \bibinfo{title}{A note on the long-time stability of pressure
  solutions to the 2d navier stokes equations},
\newblock \bibinfo{journal}{Applied Mathematics and Computation}
  \bibinfo{volume}{478} (\bibinfo{year}{2024}) \bibinfo{pages}{128839}.
  \DOIprefix\doi{https://doi.org/10.1016/j.amc.2024.128839}.
\bibitem[{Ji(2024)}]{TIS-NS-IMEX-BDF-24}
\bibinfo{author}{B.~Ji},
\newblock \bibinfo{title}{Convergence analysis of high-order imex-bdf schemes
  for the incompressible navier–stokes equations},
\newblock \bibinfo{journal}{Computers \& Fluids} \bibinfo{volume}{275}
  (\bibinfo{year}{2024}) \bibinfo{pages}{106251}.
  \DOIprefix\doi{https://doi.org/10.1016/j.compfluid.2024.106251}.
\bibitem[{DeCaria et~al.(2022)DeCaria, Gottlieb, Grant, and Layton}]{Layton-22}
\bibinfo{author}{V.~DeCaria}, \bibinfo{author}{S.~Gottlieb},
  \bibinfo{author}{Z.~J. Grant}, \bibinfo{author}{W.~J. Layton},
\newblock \bibinfo{title}{A general linear method approach to the design and
  optimization of efficient, accurate, and easily implemented time-stepping
  methods in cfd},
\newblock \bibinfo{journal}{Journal of Computational Physics}
  \bibinfo{volume}{455} (\bibinfo{year}{2022}) \bibinfo{pages}{110927}.
  \DOIprefix\doi{https://doi.org/10.1016/j.jcp.2021.110927}.
\bibitem[{Layton et~al.(2023)Layton, Pei, and Trenchea}]{Layton-23}
\bibinfo{author}{W.~Layton}, \bibinfo{author}{W.~Pei},
  \bibinfo{author}{C.~Trenchea},
\newblock \bibinfo{title}{Time step adaptivity in the method of dahlquist,
  liniger and nevanlinna},
\newblock \bibinfo{journal}{Advances in Computational Science and Engineering}
  \bibinfo{volume}{1} (\bibinfo{year}{2023}) \bibinfo{pages}{320--350}.
  \DOIprefix\doi{10.3934/acse.2023014}.
\bibitem[{Liao and Chen(2019)}]{TIS-RK-NS-19}
\bibinfo{author}{Y.~Liao, ShaokaiZhang}, \bibinfo{author}{D.~Chen},
\newblock \bibinfo{title}{Runge-kutta finite element method based on the
  characteristic for the incompressible navier-stokes equations},
\newblock \bibinfo{journal}{Advances in Applied Mathematics and Mechanics}
  \bibinfo{volume}{11} (\bibinfo{year}{2019}) \bibinfo{pages}{1415--1435}.
  \DOIprefix\doi{https://doi.org/10.4208/aamm.OA-2018-0150}.
\bibitem[{Franca and Frey(1992)}]{stab-FEM-NS-92}
\bibinfo{author}{L.~P. Franca}, \bibinfo{author}{S.~L. Frey},
\newblock \bibinfo{title}{Stabilized finite element methods: Ii. the
  incompressible navier-stokes equations},
\newblock \bibinfo{journal}{Computer Methods in Applied Mechanics and
  Engineering} \bibinfo{volume}{99} (\bibinfo{year}{1992})
  \bibinfo{pages}{209--233}.
  \DOIprefix\doi{https://doi.org/10.1016/0045-7825(92)90041-H}.
\bibitem[{Codina(2001)}]{stab-FEM-NS-01}
\bibinfo{author}{R.~Codina},
\newblock \bibinfo{title}{A stabilized finite element method for generalized
  stationary incompressible flows},
\newblock \bibinfo{journal}{Computer Methods in Applied Mechanics and
  Engineering} \bibinfo{volume}{190} (\bibinfo{year}{2001})
  \bibinfo{pages}{2681--2706}.
  \DOIprefix\doi{https://doi.org/10.1016/S0045-7825(00)00260-7}.
\bibitem[{Zhu and Chen(2015)}]{stab-FEM-NS-15}
\bibinfo{author}{L.~Zhu}, \bibinfo{author}{Z.~Chen},
\newblock \bibinfo{title}{A two-level stabilized nonconforming finite element
  method for the stationary navier-stokes equations},
\newblock \bibinfo{journal}{Mathematics and Computers in Simulation}
  \bibinfo{volume}{114} (\bibinfo{year}{2015}) \bibinfo{pages}{37--48}.
  \DOIprefix\doi{https://doi.org/10.1016/j.matcom.2011.02.015},
  \bibinfo{note}{sI:CBSACM 2009}.
\bibitem[{Zheng and Shang(2020)}]{stab-FEA-NS-20}
\bibinfo{author}{B.~Zheng}, \bibinfo{author}{Y.~Shang},
\newblock \bibinfo{title}{Local and parallel stabilized finite element
  algorithms based on the lowest equal-order elements for the steady
  navier-stokes equations},
\newblock \bibinfo{journal}{Mathematics and Computers in Simulation}
  \bibinfo{volume}{178} (\bibinfo{year}{2020}) \bibinfo{pages}{464--484}.
  \DOIprefix\doi{https://doi.org/10.1016/j.matcom.2020.07.010}.
\bibitem[{Chalot et~al.(2023)Chalot, Johan, Mallet, Billard, Martin, and
  Barré}]{stab-FEM-NS-23}
\bibinfo{author}{F.~Chalot}, \bibinfo{author}{Z.~Johan},
  \bibinfo{author}{M.~Mallet}, \bibinfo{author}{F.~Billard},
  \bibinfo{author}{L.~Martin}, \bibinfo{author}{S.~Barré},
\newblock \bibinfo{title}{Extension of methods based on the stabilized finite
  element formulation for the solution of the navier–stokes equations and
  application to aerodynamic design},
\newblock \bibinfo{journal}{Computer Methods in Applied Mechanics and
  Engineering} \bibinfo{volume}{417} (\bibinfo{year}{2023})
  \bibinfo{pages}{116425}.
  \DOIprefix\doi{https://doi.org/10.1016/j.cma.2023.116425}, \bibinfo{note}{a
  Special Issue in Honor of the Lifetime Achievements of T. J. R. Hughes}.
\bibitem[{Codina and Blasco(2000)}]{stab-FEM-NS-PG-00}
\bibinfo{author}{R.~Codina}, \bibinfo{author}{J.~Blasco},
\newblock \bibinfo{title}{Stabilized finite element method for the transient
  navier–stokes equations based on a pressure gradient projection},
\newblock \bibinfo{journal}{Computer Methods in Applied Mechanics and
  Engineering} \bibinfo{volume}{182} (\bibinfo{year}{2000})
  \bibinfo{pages}{277--300}.
  \DOIprefix\doi{https://doi.org/10.1016/S0045-7825(99)00194-2}.
\bibitem[{Zhang and He(2012)}]{Press-stab-stokes-12}
\bibinfo{author}{T.~Zhang}, \bibinfo{author}{Y.~He},
\newblock \bibinfo{title}{Fully discrete finite element method based on
  pressure stabilization for the transient stokes equations},
\newblock \bibinfo{journal}{Mathematics and Computers in Simulation}
  \bibinfo{volume}{82} (\bibinfo{year}{2012}) \bibinfo{pages}{1496--1515}.
  \DOIprefix\doi{https://doi.org/10.1016/j.matcom.2012.02.007}.
\bibitem[{Hachem et~al.(2010)Hachem, Rivaux, Kloczko, Digonnet, and
  Coupez}]{stab-VMS-FEM-NS-10}
\bibinfo{author}{E.~Hachem}, \bibinfo{author}{B.~Rivaux},
  \bibinfo{author}{T.~Kloczko}, \bibinfo{author}{H.~Digonnet},
  \bibinfo{author}{T.~Coupez},
\newblock \bibinfo{title}{Stabilized finite element method for incompressible
  flows with high reynolds number},
\newblock \bibinfo{journal}{Journal of Computational Physics}
  \bibinfo{volume}{229} (\bibinfo{year}{2010}) \bibinfo{pages}{8643--8665}.
  \DOIprefix\doi{https://doi.org/10.1016/j.jcp.2010.07.030}.
\bibitem[{Coupez and Hachem(2013)}]{stab-VMS-FEM-NS-13}
\bibinfo{author}{T.~Coupez}, \bibinfo{author}{E.~Hachem},
\newblock \bibinfo{title}{Solution of high-reynolds incompressible flow with
  stabilized finite element and adaptive anisotropic meshing},
\newblock \bibinfo{journal}{Computer Methods in Applied Mechanics and
  Engineering} \bibinfo{volume}{267} (\bibinfo{year}{2013})
  \bibinfo{pages}{65--85}.
  \DOIprefix\doi{https://doi.org/10.1016/j.cma.2013.08.004}.
\bibitem[{Ahmed et~al.(2016)Ahmed, Rebollo, John, and
  Rubino}]{stab-FEM-localProj-16}
\bibinfo{author}{N.~Ahmed}, \bibinfo{author}{T.~C. Rebollo},
  \bibinfo{author}{V.~John}, \bibinfo{author}{S.~Rubino},
\newblock \bibinfo{title}{Analysis of a full space–time discretization of the
  navier–stokes equations by a local projection stabilization method},
\newblock \bibinfo{journal}{IMA Journal of Numerical Analysis}
  \bibinfo{volume}{37} (\bibinfo{year}{2016}) \bibinfo{pages}{1437--1467}.
  \DOIprefix\doi{10.1093/imanum/drw048}.
\bibitem[{Baptista et~al.(2019)Baptista, Bento, Lima, Santos, Valli, and
  Catabriga}]{stab-FEM-subgrid-19}
\bibinfo{author}{R.~Baptista}, \bibinfo{author}{S.~S. Bento},
  \bibinfo{author}{L.~M. Lima}, \bibinfo{author}{I.~P. Santos},
  \bibinfo{author}{A.~M.~P. Valli}, \bibinfo{author}{L.~Catabriga},
\newblock \bibinfo{title}{A nonlinear subgrid stabilization parameter-free
  method to solve incompressible navier-stokes equations at high reynolds
  numbers},
\newblock in: \bibinfo{editor}{S.~Misra}, \bibinfo{editor}{O.~Gervasi},
  \bibinfo{editor}{B.~Murgante}, \bibinfo{editor}{E.~Stankova},
  \bibinfo{editor}{V.~Korkhov}, \bibinfo{editor}{C.~Torre},
  \bibinfo{editor}{A.~M.~A. Rocha}, \bibinfo{editor}{D.~Taniar},
  \bibinfo{editor}{B.~O. Apduhan}, \bibinfo{editor}{E.~Tarantino} (Eds.),
  \bibinfo{booktitle}{Computational Science and Its Applications -- ICCSA
  2019}, \bibinfo{publisher}{Springer International Publishing},
  \bibinfo{address}{Cham}, \bibinfo{year}{2019}, pp. \bibinfo{pages}{134--148}.
\bibitem[{Choi et~al.(1997)Choi, Choi, and Yoo}]{TIS-FEM-SUPG-NS-97}
\bibinfo{author}{H.~Choi}, \bibinfo{author}{H.~Choi}, \bibinfo{author}{J.~Yoo},
\newblock \bibinfo{title}{A fractional four-step finite element formulation of
  the unsteady incompressible navier-stokes equations using supg and linear
  equal-order element methods},
\newblock \bibinfo{journal}{Computer Methods in Applied Mechanics and
  Engineering} \bibinfo{volume}{143} (\bibinfo{year}{1997})
  \bibinfo{pages}{333--348}.
  \DOIprefix\doi{https://doi.org/10.1016/S0045-7825(96)01156-5}.
\bibitem[{{do Carmo} and Alvarez(2003)}]{stab-SUPG-FEM-03}
\bibinfo{author}{E.~G.~D. {do Carmo}}, \bibinfo{author}{G.~B. Alvarez},
\newblock \bibinfo{title}{A new stabilized finite element formulation for
  scalar convection--diffusion problems: the streamline and approximate
  upwind/petrov--galerkin method},
\newblock \bibinfo{journal}{Computer Methods in Applied Mechanics and
  Engineering} \bibinfo{volume}{192} (\bibinfo{year}{2003})
  \bibinfo{pages}{3379--3396}.
  \DOIprefix\doi{https://doi.org/10.1016/S0045-7825(03)00292-5}.
\bibitem[{Liu et~al.(2021)Liu, Gao, Zhu, and Jiang}]{satb-FEM-SUPG-NS-21}
\bibinfo{author}{M.~Liu}, \bibinfo{author}{G.~Gao}, \bibinfo{author}{H.~Zhu},
  \bibinfo{author}{C.~Jiang},
\newblock \bibinfo{title}{A cell-based smoothed finite element method
  stabilized by implicit supg/spgp/fractional step method for incompressible
  flow},
\newblock \bibinfo{journal}{Engineering Analysis with Boundary Elements}
  \bibinfo{volume}{124} (\bibinfo{year}{2021}) \bibinfo{pages}{194--210}.
  \DOIprefix\doi{https://doi.org/10.1016/j.enganabound.2020.12.018}.
\bibitem[{Razafindralandy(2005)}]{dina-thesis}
\bibinfo{author}{D.~Razafindralandy}, \bibinfo{title}{Contribution à
  l'{É}tude {M}athématique et {N}umérique de la {S}imulation des {G}randes
  {É}chelles}, Ph.D. thesis, La Rochelle University, \bibinfo{year}{2005}.
\bibitem[{Deeb(2015)}]{deeb-thesis}
\bibinfo{author}{A.~Deeb}, \bibinfo{title}{Intégrateurs {T}emporels {B}asés
  sur la {R}esommation des {S}éries {D}ivergentes. Applications en
  {M}écanique}, Ph.D. thesis, La Rochelle University, \bibinfo{year}{2015}.
\bibitem[{{B}orel(1901)}]{borel-1901}
\bibinfo{author}{E.~{B}orel}, \bibinfo{title}{Lecons sur les s{\'e}ries
  divergentes}, Collection de monographies sur la th{\'e}orie des fonctions,
  \bibinfo{publisher}{Gauthier-Villars}, \bibinfo{year}{1901}.
\bibitem[{Razafindralandy and Hamdouni(2013)}]{dina-2012}
\bibinfo{author}{D.~Razafindralandy}, \bibinfo{author}{A.~Hamdouni},
\newblock \bibinfo{title}{Time integration algorithm based on divergent series
  resummation, for ordinary and partial differential equations},
\newblock \bibinfo{journal}{Journal of Computational Physics}
  \bibinfo{volume}{236} (\bibinfo{year}{2013}) \bibinfo{pages}{56--73}.
  \DOIprefix\doi{https://doi.org/10.1016/j.jcp.2012.10.022}.
\bibitem[{Deeb et~al.(2014)Deeb, Hamdouni, Liberge, and
  Razafindralandy}]{ahmad_bpl_2014}
\bibinfo{author}{A.~Deeb}, \bibinfo{author}{A.~Hamdouni},
  \bibinfo{author}{E.~Liberge}, \bibinfo{author}{D.~Razafindralandy},
\newblock \bibinfo{title}{Borel-{L}aplace summation method used as time
  integration scheme},
\newblock \bibinfo{journal}{ESAIM: Procedings and Surveys} \bibinfo{volume}{45}
  (\bibinfo{year}{2014}) \bibinfo{pages}{318--327}.
  \DOIprefix\doi{https://doi.org/10.1051/proc/201445033}.
\bibitem[{Deeb et~al.(2016)Deeb, Razafindralandy, and
  Hamdouni}]{ahmad_comp_bpl_sfg_2015}
\bibinfo{author}{A.~Deeb}, \bibinfo{author}{D.~Razafindralandy},
  \bibinfo{author}{A.~Hamdouni},
\newblock \bibinfo{title}{Comparison between {B}orel-{P}ad{\'e} summation and
  factorial series, as time integration methods},
\newblock \bibinfo{journal}{Discrete and Continuous Dynamical Systems - Series
  S} \bibinfo{volume}{9} (\bibinfo{year}{2016}) \bibinfo{pages}{393--408}.
  \DOIprefix\doi{10.3934/dcdss.2016003}.
\bibitem[{Razafindralandy et~al.(2017)Razafindralandy, Hamdouni, and
  Deeb}]{ahmad_icnpaa_2016}
\bibinfo{author}{D.~Razafindralandy}, \bibinfo{author}{A.~Hamdouni},
  \bibinfo{author}{A.~Deeb},
\newblock \bibinfo{title}{Considering factorial series as time integration
  method},
\newblock in: \bibinfo{booktitle}{11th International Conference on Mathematical
  Problems in Engineering, Aerospace and Sciences}, volume
  \bibinfo{volume}{1798} of \textit{\bibinfo{series}{AIP Conference
  Proceedings}}, \bibinfo{publisher}{American Institute of Physics},
  \bibinfo{year}{2017}. \DOIprefix\doi{10.1063/1.4972721}.
\bibitem[{Razafindralandy et~al.(2019)Razafindralandy, Salnikov, Hamdouni, and
  Deeb}]{ahmad_robust_integrators_2019}
\bibinfo{author}{D.~Razafindralandy}, \bibinfo{author}{V.~Salnikov},
  \bibinfo{author}{A.~Hamdouni}, \bibinfo{author}{A.~Deeb},
\newblock \bibinfo{title}{Some robust integrators for large time dynamics},
\newblock \bibinfo{journal}{Advanced Modeling and Simulation in Engineering
  Sciences} \bibinfo{volume}{6} (\bibinfo{year}{2019}).
  \DOIprefix\doi{10.1186/s40323-019-0130-2}.
\bibitem[{Tayeh et~al.(2021)Tayeh, Girault, Guevel, and Cadou}]{tayeh-21}
\bibinfo{author}{C.~Tayeh}, \bibinfo{author}{G.~Girault},
  \bibinfo{author}{Y.~Guevel}, \bibinfo{author}{J.~Cadou},
\newblock \bibinfo{title}{Numerical time perturbation and resummation methods
  for nonlinear ode},
\newblock \bibinfo{journal}{Nonlinear Dyn} \bibinfo{volume}{103}
  (\bibinfo{year}{2021}) \bibinfo{pages}{617--642}.
  \DOIprefix\doi{https://doi.org/10.1007/s11071-020-06137-w}.
\bibitem[{Deeb et~al.(2022)Deeb, Hamdouni, and Razafindralandy}]{DEEB_2022_bpl}
\bibinfo{author}{A.~Deeb}, \bibinfo{author}{A.~Hamdouni},
  \bibinfo{author}{D.~Razafindralandy},
\newblock \bibinfo{title}{Performance of {B}orel-{P}adé-{L}aplace integrator
  for the solution of stiff and non-stiff problems},
\newblock \bibinfo{journal}{Applied Mathematics and Computation}
  \bibinfo{volume}{426} (\bibinfo{year}{2022}).
  \DOIprefix\doi{10.1016/j.amc.2022.127118}.
\bibitem[{Deeb and Dutykh(2024)}]{deeb:stab-serie}
\bibinfo{author}{A.~Deeb}, \bibinfo{author}{D.~Dutykh},
\newblock \bibinfo{title}{Stabilized time series expansions for high-order
  finite element solutions of partial differential equation},
\newblock \bibinfo{journal}{Studies in Applied Mathematics}
  \bibinfo{volume}{153} (\bibinfo{year}{2024}).
  \DOIprefix\doi{https://doi.org/10.1111/sapm.12708}.
\bibitem[{Brezzi(1974)}]{brezzi-75}
\bibinfo{author}{F.~Brezzi},
\newblock \bibinfo{title}{On the uniqueness, existence and approximation of
  saddle point problem arising from {L}agrangian multipliers},
\newblock \bibinfo{journal}{RAIRO Anal. Numer} \bibinfo{volume}{8}
  (\bibinfo{year}{1974}) \bibinfo{pages}{129--151}.
  \DOIprefix\doi{10.1051/m2an/197408R201291}.
\bibitem[{Brenner and Scott(2008)}]{brenner-10}
\bibinfo{author}{S.~C. Brenner}, \bibinfo{author}{L.~R. Scott},
  \bibinfo{title}{The Mathematical Theory of Finite Element Methods}, Texts in
  Applied Mathematics, \bibinfo{edition}{3} ed., \bibinfo{publisher}{Springer
  New York}, \bibinfo{address}{New York, NY}, \bibinfo{year}{2008}.
  \DOIprefix\doi{/10.1007/978-0-387-75934-0}.
\bibitem[{Céa(1964)}]{Cea-64}
\bibinfo{author}{J.~Céa},
\newblock \bibinfo{title}{Approximation variationnelle des problèmes aux
  limites},
\newblock \bibinfo{journal}{Annales de l'institut Fourier} \bibinfo{volume}{2}
  (\bibinfo{year}{1964}) \bibinfo{pages}{345--444}. \bibinfo{note}{PhD Thesis}.
\bibitem[{Balser(2000)}]{balser-99}
\bibinfo{author}{W.~Balser}, \bibinfo{title}{Formal Power Series and Linear
  Systems of Meromorphic Ordinary Differential Equations}, Universitext,
  \bibinfo{publisher}{Springer New York, NY}, \bibinfo{year}{2000}.
  \DOIprefix\doi{https://doi.org/10.1007/b97608}.
\bibitem[{Delabaere and Rasoamanana(2007)}]{Delabaere-07}
\bibinfo{author}{E.~Delabaere}, \bibinfo{author}{J.-M. Rasoamanana},
\newblock \bibinfo{title}{Sommation effective d{\textquoteright}une somme de
  {Borel} par s\'eries de factorielles},
\newblock \bibinfo{journal}{Annales de l'Institut Fourier} \bibinfo{volume}{57}
  (\bibinfo{year}{2007}) \bibinfo{pages}{421--456}. \URLprefix
  \url{https://aif.centre-mersenne.org/articles/10.5802/aif.2263/}.
  \DOIprefix\doi{10.5802/aif.2263}.
\bibitem[{Thomann(9091)}]{Thomann-90}
\bibinfo{author}{J.~Thomann},
\newblock \bibinfo{title}{Resommation des series formelles. solutions
  d'équations différentielles linéaires ordinaires du second ordre dans le
  champ complexe au voisinage de singularités irrégulières.},
\newblock \bibinfo{journal}{Numerische Mathematik} \bibinfo{volume}{58}
  (\bibinfo{year}{1990/91}) \bibinfo{pages}{503--536}. \URLprefix
  \url{http://eudml.org/doc/133514}.
\bibitem[{Baratta et~al.(2023)Baratta, Dean, Dokken, Habera, Hale, Richardson,
  Rognes, Scroggs, Sime, and Wells}]{BarattaEtal2023}
\bibinfo{author}{I.~A. Baratta}, \bibinfo{author}{J.~P. Dean},
  \bibinfo{author}{J.~S. Dokken}, \bibinfo{author}{M.~Habera},
  \bibinfo{author}{J.~S. Hale}, \bibinfo{author}{C.~N. Richardson},
  \bibinfo{author}{M.~E. Rognes}, \bibinfo{author}{M.~W. Scroggs},
  \bibinfo{author}{N.~Sime}, \bibinfo{author}{G.~N. Wells},
  \bibinfo{title}{{DOLFINx}: the next generation {FEniCS} problem solving
  environment}, \bibinfo{howpublished}{preprint}, \bibinfo{year}{2023}.
  \DOIprefix\doi{10.5281/zenodo.10447666}.
\bibitem[{Scroggs et~al.(2022)Scroggs, Dokken, Richardson, and
  Wells}]{ScroggsEtal2022}
\bibinfo{author}{M.~W. Scroggs}, \bibinfo{author}{J.~S. Dokken},
  \bibinfo{author}{C.~N. Richardson}, \bibinfo{author}{G.~N. Wells},
\newblock \bibinfo{title}{Construction of arbitrary order finite element
  degree-of-freedom maps on polygonal and polyhedral cell meshes},
\newblock \bibinfo{journal}{ACM Transactions on Mathematical Software}
  \bibinfo{volume}{48} (\bibinfo{year}{2022}) \bibinfo{pages}{{18:1--18:23}}.
  \DOIprefix\doi{10.1145/3524456}.
\bibitem[{Guermond et~al.(2003)Guermond, Oden, and Prudhomme}]{nsalpha-03}
\bibinfo{author}{J.~Guermond}, \bibinfo{author}{J.~Oden},
  \bibinfo{author}{S.~Prudhomme},
\newblock \bibinfo{title}{An interpretation of the {N}avier--{S}tokes-alpha
  model as a frame-indifferent {L}eray regularization},
\newblock \bibinfo{journal}{Physica D: Nonlinear Phenomena}
  \bibinfo{volume}{177} (\bibinfo{year}{2003}) \bibinfo{pages}{23--30}.
  \DOIprefix\doi{https://doi.org/10.1016/S0167-2789(02)00748-0}.
\bibitem[{Leray(1934)}]{leray-34}
\bibinfo{author}{J.~Leray},
\newblock \bibinfo{title}{Sur le mouvement d'un liquide visqueux emplissant
  l'espace.},
\newblock \bibinfo{journal}{Acta Math.} \bibinfo{volume}{63}
  (\bibinfo{year}{1934}) \bibinfo{pages}{193--248}.
  \DOIprefix\doi{https://doi.org/10.1007/BF02547354}.
\bibitem[{Hager(1984)}]{hager-84}
\bibinfo{author}{W.~W. Hager},
\newblock \bibinfo{title}{Condition estimates},
\newblock \bibinfo{journal}{SIAM Journal on Scientific and Statistical
  Computing} \bibinfo{volume}{5} (\bibinfo{year}{1984})
  \bibinfo{pages}{311--316}. \DOIprefix\doi{10.1137/0905023}.
\bibitem[{Farag\'{o} and Horv\'{a}th(2006)}]{DMP_Farago}
\bibinfo{author}{I.~Farag\'{o}}, \bibinfo{author}{R.~Horv\'{a}th},
\newblock \bibinfo{title}{Discrete maximum principle and adequate
  discretizations of linear parabolic problems},
\newblock \bibinfo{journal}{SIAM Journal on Scientific Computing}
  \bibinfo{volume}{28} (\bibinfo{year}{2006}) \bibinfo{pages}{2313--2336}.
  \DOIprefix\doi{10.1137/050627241}.
\bibitem[{Geuzaine and Remacle(2009)}]{gmsh-09}
\bibinfo{author}{C.~Geuzaine}, \bibinfo{author}{J.~F. Remacle},
\newblock \bibinfo{title}{Gmsh: a three-dimensional finite element mesh
  generator with built-in pre- and post-processing facilities},
\newblock \bibinfo{journal}{International Journal for Numerical Methods in
  Engineerins} \bibinfo{volume}{79} (\bibinfo{year}{2009})
  \bibinfo{pages}{1309--1331}.
  \DOIprefix\doi{https://doi.org/10.1002/nme.2579}.
\bibitem[{Ciarlet(2002)}]{ciarlet-02}
\bibinfo{author}{P.~G. Ciarlet}, \bibinfo{title}{The Finite Element Method for
  Elliptic Problems}, Classics in Applied Mathematics,
  \bibinfo{publisher}{Society for Industrial and Applied Mathematics},
  \bibinfo{address}{Philadelphia}, \bibinfo{year}{2002}.
  \DOIprefix\doi{10.1137/1.9780898719208}.

\end{thebibliography}

\end{document}